\documentclass[english]{article}
\usepackage[latin1]{inputenc}
\usepackage[T1]{fontenc}
\usepackage{amsmath,amssymb,amscd}
\usepackage{epsfig}
\usepackage{mathrsfs}
\usepackage{euscript}
\usepackage{pifont}
\usepackage[all]{xy}

\newcommand{\smind}[1]{{\!\mbox{\fontsize{2}{2}\selectfont #1}}}

\newcommand{\id}{\mathit{id}}

\newcommand{\Ho}{\operatorname{Ho}}
\newcommand{\AdEq}{\operatorname{AdEq}}
\newcommand{\diag}{\operatorname{diag}}

\newcommand{\op}{\operatorname{op}}
\newcommand{\B}{\mathcal{B}}
\newcommand{\C}{\mathcal{C}}
\newcommand{\D}{\mathcal{D}}
\newcommand{\E}{\mathcal{E}}
\newcommand{\F}{\mathcal{F}}
\newcommand{\G}{\mathcal{G}}
\newcommand{\I}{\mathcal{I}}

\newcommand{\R}{\mathbb{R}}
\newcommand{\T}{\mathcal{T}}
\newcommand{\U}{\mathcal{U}}

\newcommand{\Z}{\mathbb{Z}}
\newcommand{\bbC}{\mathbb{C}}

\newcommand{\Set}{\mathcal{S}\!\mathit{et}}
\newcommand{\Top}{\mathcal{T}\!\!\mathit{op}}

\newcommand{\Con}{\operatorname{Con}}

\newcommand{\tr}{\operatorname{tr}}

\newcommand{\colim}{\operatornamewithlimits{colim}}

\newcommand{\Sing}{\operatorname{Sing}_{\bullet}}

\newcommand{\abs}[1]{\lvert {#1} \rvert}

\newcommand{\co}{\colon\thinspace}    

\newtheorem{Thm}{Theorem}[section]

\newtheorem{Prop}[Thm]{Proposition}
\newtheorem{Lem}[Thm]{Lemma}

\newtheorem{Cor}[Thm]{Corollary}

\newtheorem{Concl}[Thm]{Conclusion}

\newtheorem{Def}[Thm]{Definition}
\newtheorem{Exa}[Thm]{Example}
\newtheorem{Rem}[Thm]{Remark}

\newenvironment{proof}{\mbox{}\newline\textbf{Proof:} }{\nopagebreak\hfill$\square$\\}

\title{Two-Categorical Bundles and Their Classifying Spaces}

\author{Nils A. Baas, Marcel B\"o{}kstedt and Tore August Kro\footnote{The first and third author gratefully acknowledge support
by the Institute Mittag-Leffler (Djursholm, Sweeden).} }

\begin{document}
\maketitle

\begin{abstract}
For a $2$-category $2\C$ we associate a notion of a principal
$2\C$-bundle. In case of the $2$-category of $2$-vector spaces in
the sense of M.M. Kapranov and V.A. Voevodsky this gives the the
$2$-vector bundles of N.A. Baas, B.I. Dundas and J. Rognes. Our main
result says that the geometric nerve of a good $2$-category is a
classifying space for the associated principal $2$-bundles. In the
process of proving this we develop a lot of powerful machinery which
may be useful in further studies of $2$-categorical topology. As a
corollary we get a new proof of the classification of principal
bundles. A calculation based on the main theorem shows that the
principal $2$-bundles associated to the $2$-category of $2$-vector
spaces in the sense of J.C. Baez and A.S. Crans split, up to
concordance, as two copies of ordinary vector bundles. When $2\C$ is
a cobordism type $2$-category we get a new notion of cobordism-bundles
which turns out to be classified by the Madsen-Weiss spaces.
\end{abstract}

\section{Introduction and main result}

The main purpose of this paper is to introduce a general notion of
bundles associated to topological $2$-categories, and classify
these. This encompasses the notion of $2$-vector bundles developed
in~\cite{BaasDundasRognes:04}.

In this paper we mean by a $2$-category what is sometimes called a
bicategory~\cite{Benabou:67}, or a weak $2$-category, i.e. where the
associativity conditions of the horizontal composition is relaxed.
We specify strict $2$-categories when we need them.

The key idea behind principal $2\C$-bundles is to categorify the
transition data description of a principal $G$-bundle by formally
replacing the group $G$ with a $2$-category $2\C$. Given an ordered
open cover $\{U_\alpha\}$ of a space, we associate objects in $2\C$
to points in each $U_\alpha$ and $1$-morphisms to points in each
double intersection $U_\alpha\cap U_\beta$. Instead of the ordinary
cocycle condition, we associate $2$-morphisms to points in each
triple intersection. These are subject to a higher cocycle condition
on each quadruple intersection.

A $2$-category has a geometric nerve, $\Delta2\C$, which is a
simplicial set,~\cite{Duskin:02}. There are obvious generalizations
of these nerves to topological $2$-categories, which produce
simplicial spaces. A principal $2$-bundle is the same as a
simplicial map from the ordered \v{C}ech complex, $U_\bullet$, of
the ordered open cover to the geometric nerve of the $2$-category.
Alternatively, one may describe a principal $2$-bundle as a
$2$-functor from a certain $2$-category into $2\C$.

Two principal $2$-bundles over $X$ are concordant if they are the
restriction to $X\times\{0,1\}$ of some principal $2$-bundle over
$X\times I$. This is an equivalence relation, and our aim is to
determine the set of concordance classes. We restrict to base spaces
$X$ having the homotopy type of a CW complex. A space $B$ is said to
be a classifying space if there is a bijection between principal
$2$-bundles over such $X$ and homotopy classes of maps from $X$ into
$B$. Our main theorem, which is a $2$-categorical answer to M.
Weiss' question,~\cite{Weiss:05}, ``What does the classifying space
of a category classify?'', says:

\begin{Thm}\label{thm:main}
The realization of the geometric nerve of a good topological
$2$-category $2\C$ is a classifying space for principal
$2\C$-bundles.
\end{Thm}

We get the homotopy class in $[X,|\Delta2\C|]$ corresponding to a
principal $2\C$-bundle $U_\bullet\rightarrow\Delta2\C$ by taking
geometric realization. The adjective ``good'' concerns the topology
of $2\C$, see Definition~\ref{def:good}, and this condition is
satisfied whenever the total space of $2$-morphisms is a CW complex
and all source-target maps are Hurewicz fibrations, see
Theorem~\ref{Thm:good}.

In~\cite{BaasDundasRognes:04} Baas, Dundas and Rognes introduce
$2$-vector bundles and define an associated second order $K$-theory.
The representing spectrum of this cohomology theory is $K(ku)$,
see~\cite{BaasDundasRognes:04,BaasDundasRichterRognes:07}, and it
qualifies as a form of elliptic cohomology theory since the
chromatic filtration essentially is $2$, see~\cite{Ausoni:06}. We
recover in Example~\ref{exa:2KVn} the notion of $2$-vector bundles
as the principal $2\C$-bundles for a suitable $2\C$ defined in the
spirit of Kapranov and Voevodskys $2$-vector spaces.

A similar construction based on Baez and Crans version of $2$-vector
spaces~\cite{BaezCrans:04} leads to a cohomology theory which is two
copies of ordinary $K$-theory, and hence not a form of elliptic
cohomology, see Conclusion~\ref{concl:BC}.

The last four sections are more technical. They deal with facts
which are used in the proof of the main theorem. In
Section~\ref{sect:Kan} we study the following problem. Suppose that
$Z_\bullet$ is a simplicial space. By applying the simplicial
functor degreewise, we get a bisimplicial set. Let us diagonalize
this bisimplicial set. What are reasonable conditions on $Z_\bullet$
so that this diagonal simplicial set has the Kan property? The
answer to this question, Theorem~\ref{th:Kan}, is a topological Kan
condition.

In Section~\ref{sect:GSF} we provide criteria for what we think a
good and sufficiently fibrant topological $2$-category $2\C$ should
be, i.e. we ask when the geometric nerve $\Delta2\C$ is a good
simplicial space and when it satisfies the topological Kan
condition.

In Section~\ref{sect:concordance} we generalize the concept of
concordance classes by replacing the geometric nerve of a
topological category with an arbitrary simplicial space,
$Z_\bullet$. We give some general properties (homotopy invariance,
exact sequence, gluing) for such generalized concordance classes. We
also show how to replace arbitrary $Z_\bullet$ by another simplicial
space $\tilde{Z}_\bullet$, satisfying the topological Kan condition,
without changing the set of concordance classes. This relies on
Quillen's small object argument.

In Section~\ref{sect:proof} we generalize the statement of
Theorem~\ref{thm:main}, and prove that concordance classes of
$Z_\bullet$-bundles are classified by $|Z_\bullet|$ provided that
$Z_\bullet$ is a good simplicial space,
Theorem~\ref{thm:generalizationofmain}.

\section{Principal $2\C$-bundles}\label{sect:principal}

In this section we define what a principal $2\C$-bundle is. In our
definition we specify transition data. However, there are two other,
equivalent, definitions; either as a simplicial map from the ordered
\v{C}ech complex to the geometric nerve, or as a continuous functor
from a certain $2$-category into $2\C$. We end the section with the
definition of concordance and the restatement of our main theorem.

Our structure $2$-category $2\C$ will be \emph{topological}. This
may be understood in at least two different ways, either as a
$2$-category enriched in spaces or as a $2$-category internal in
spaces. We take the latter, and the most general, point of view.
Hence, misunderstanding our usage of the term topological
$2$-category will not lead to any mistakes. Discrete (or ordinary)
$2$-categories is the special case where the spaces of objects, $1$-
and $2$-morphisms have the discrete topology.

We apply the following notation: $2\C_0$, $2\C_1$, and $2\C_2$ are
the topological spaces of objects, $1$-morphisms and $2$-morphisms
respectively. If $x,y\in2\C_0$, then $2\C_1(x,y)$ is the space of
$1$-morphisms $x\rightarrow y$, while $2\C(x,y)$ is the topological
category with $1$-morphisms $x\rightarrow y$ as objects and
$2$-morphisms between these as morphisms. If $f,g\in2\C_1(x,y)$,
then $2\C_2(f,g)$ is the space of $2$-morphisms $f\Rightarrow g$. We
denote horizontal composition by $*$, whereas vertical composition
(of $2$-morphisms) is written by juxtaposition. The natural
associativity, left and right unit coherence isomorphisms are
denoted $\underline{\alpha}$, $\underline{\lambda}$, and
$\underline{\rho}$ respectively.

As defined in~\cite{BaasDundasRognes:04}, an \emph{ordered open
cover} $\U$ of a topological space $X$ consists of a family
$\{U_\alpha\}$ of open subsets $U_\alpha\subseteq X$ indexed over a
partially ordered set $\I$ such that the family cover $X$, i.e.
$\bigcup_{\alpha\in \I}U_\alpha = X$, and whenever a finite
intersection $U_{\alpha_0\cdots\alpha_k}=U_{\alpha_0}\cap\cdots\cap
U_{\alpha_k}$ is nonempty then the partial ordering of $\I$
restricts to a total ordering on $\{\alpha_0,\ldots,\alpha_k\}$.

\begin{Def}\label{def:principal2Cbundle}
Let $X$ be a topological space, and let $2\C$ be a topological
$2$-category. A \emph{principal $2\C$-bundle} $\E$ over~$X$ consists
of an ordered open cover $\U$, indexed by $\I$, together with
\begin{itemize}
\item[(1)] for each $\alpha$ in $\I$ a continuous map
\[
V_\alpha\co U_\alpha\rightarrow 2\C_0,
\]
and
\item[(2)] for each $\alpha<\beta$ a
continuous family of $1$-morphisms
\[
E_{\alpha\beta}\co U_{\alpha\beta}\rightarrow 2\C_1
\]
with source and target $V_{\alpha}\xrightarrow{E_{\alpha\beta}}
V_{\beta}$, and
\item[(3)] for each $\alpha<\beta<\gamma$ a
continuous family of 2-morphisms
\[
\phi_{\alpha\beta\gamma}\co U_{\alpha\beta\gamma}\rightarrow
2\C_2
\]
with source and target
$E_{\alpha\gamma}\overset{\phi_{\alpha\beta\gamma}}{\Rightarrow}
E_{\beta\gamma}* E_{\alpha\beta}$, such that
\item[(4)] the diagram
\[
\xymatrix{ E_{\gamma\delta} * (E_{\beta\gamma} * E_{\alpha\beta})
\ar@{=>}[rr]^{\underline{\alpha}} && (E_{\gamma\delta} *
E_{\beta\gamma}) * E_{\alpha\beta}\\
E_{\gamma\delta} * E_{\alpha\gamma}
\ar@{=>}[u]_{E_{\gamma\delta}*\phi_{\alpha\beta\gamma} } &
E_{\alpha\delta} \ar@{=>}[r]_{\phi_{\alpha\beta\delta}}
\ar@{=>}[l]^{\phi_{\alpha\gamma\delta}} & E_{\beta\delta} *
E_{\alpha\beta}
\ar@{=>}[u]^{\phi_{\beta\gamma\delta}*E_{\alpha\beta}} }
\]
commutes over~$U_{\alpha\beta\gamma\delta}$, for each chain $\alpha
< \beta < \gamma < \delta$ in $\I$. The $2$-morphism
$\underline{\alpha}$ is the natural associativity isomorphism of the
$2$-category $2\C$.
\end{itemize}
\end{Def}

In Section~\ref{sect:examples} we will give concrete examples of
structure $2$-categories giving rise to different types of principal
$2$-bundles. In particular, we recover the charted $2$-vector
bundles of~\cite{BaasDundasRognes:04} for a suitably chosen $2\C$.
In the forthcoming examples, it is often possible to interpret
$V_\alpha$, $E_{\alpha\beta}$, and $\phi_{\alpha\beta\gamma}$ as
bundles of objects, $1$-morphisms, and $2$-morphisms respectively.

Recall from~\cite{BullejosCegarra:03,Duskin:02,Street:87}, the
notion of the geometric nerve $\Delta2\C$ of a discrete $2$-category
$2\C$. The simplicial set $\Delta2\C$ has $0$-simplices the objects
$x_0$ of $2\C$. The $1$-simplices are the $1$-morphisms
$x_0\xrightarrow{x_{01}}x_1$ of $2\C$, and the $2$-simplices are
triangles
\[
\xymatrix{ & x_1 \ar[dr]^{x_{12}}  \\
x_0 \ar[ur]^{x_{01}} \ar[rr]_{x_{02}}& \ar@{=>}[u]_(.3){x_{012}} &
x_2 }
\]
where $x_{012}$ is a $2$-morphism $x_{02}\Rightarrow x_{12}*
x_{01}$. For $n\geq3$ the $n$-simplices are built from
$2$-simplices, such that for each subtetrahedron the obvious
coherence condition for $2$-morphisms is satisfied. The geometric
nerve is $3$-coskeletal, i.e. any simplex in the geometric nerve is
uniquely determined by its $3$-skeleton,
see~\cite[Proposition~3.1]{Beke:04}. The definition of the geometric
nerve extends, in the obvious way, to the cases where $2\C$ is
topological or simplicial, and then $\Delta2\C$ becomes a simplicial
space or a bisimplicial set respectively.

To an ordered cover $\U$ of $X$ we associate the \v{C}ech complex
$\check{C}_{\bullet}$ and the ordered \v{C}ech complex
$U_{\bullet}$. Their definitions are as follows:
$\check{C}_{\bullet}$ and $U_{\bullet}$ are both $1$-coskeletal
simplicial spaces (and thus determined by their $0$- and
$1$-simplices). We define $\check{C}_{0}$ and $\check{C}_{1}$ to be
\[
\coprod_{\alpha} U_{\alpha}\quad\text{and}\quad
\coprod_{\alpha,\beta} U_{\alpha\beta}
\]
respectively. The the ordered \v{C}ech complex $U_{\bullet}$ is
defined as the sub-simplicial space of $\check{C}_{\bullet}$ having
the same $0$-simplices, but the $1$-simplices being the disjoint
union of all $U_{\alpha\beta}$'s with $\alpha\leq\beta$.

Let $|-|$ denote geometric realization. From~\cite[Theorem~2.1 and
Proposition~2.6]{DuggerIsaksen:04} we have that the spaces
$|\check{C}_{\bullet}|$ and $|U_{\bullet}|$ are weakly equivalent to
$X$ via the natural maps.

With the geometric nerve, we get our first reformulation of the
definition of a principal $2\C$-bundle:

\begin{Prop}
There is a one-to-one correspondence between principal $2\C$-bundles
subordinate to the ordered cover $\U$ and simplicial maps
$U_\bullet\rightarrow\Delta2\C$.
\end{Prop}

\begin{proof}
Let $\tr_3$ denote the functor truncating a simplicial space at its
$3$-skeleton. Since the geometric nerve is $3$-coskeletal, a
simplicial map $U_\bullet\rightarrow \Delta2\C$ corresponds to its
truncation $\tr_3 U_\bullet \rightarrow \tr_3\Delta2\C$. The ordered
\v{C}ech complex has free
degeneracies,~\cite[Definition~A.4]{DuggerIsaksen:04}, hence the
truncated map is uniquely given as maps from the non-degenerate
$0$-, $1$-, $2$-, and $3$-simplices of $U_\bullet$ satisfying the
face relations. These data are precisely what is written out in
Definition~\ref{def:principal2Cbundle} above.
\end{proof}

Below we will give two additional reformulations of
Definition~\ref{def:principal2Cbundle}. The idea is to look at
$U_\bullet\rightarrow \Delta2\C$ and think about the geometric nerve
as a right adjoint. To illustrate this we consider the
$1$-categorical example $U_\bullet\rightarrow N_\bullet \C$, where
$N_\bullet \C$ is the nerve of the category $\C$.

Define $X_\U^{\text{ord}}$ to be the topological category whose
spaces of objects and morphisms are the disjoint unions
\[
\coprod_{\alpha} U_\alpha\quad\text{and}\quad
\coprod_{\alpha\leq\beta} U_{\alpha\beta}
\]
respectively. Whenever some underlying point $x\in X$ is contained
in $U_\alpha$ we write $x_\alpha$ for the object of
$X_\U^{\text{ord}}$ corresponding to this copy of $x$. Similarly, we
denote a typical morphism of $X_\U^{\text{ord}}$ by
$x_{\alpha\beta}$. The source and target of $x_{\alpha\beta}$ are
$x_\alpha$ and $x_\beta$ respectively, and composition is given by
$x_{\beta\gamma}x_{\alpha\beta}=x_{\alpha\gamma}$.

Intuitively, $X_\U^{\text{ord}}$ is the left adjoint of $N_\bullet$
applied to the ordered \v{C}ech complex $U_\bullet$. Hence, there is
a one-to-one correspondence between simplicial maps
$U_\bullet\rightarrow N_\bullet\C$ and continuous functors
$X_\U^{\text{ord}}\rightarrow \C$. Also observe that $N_\bullet
X_\U^{\text{ord}}$ is isomorphic to $U_\bullet$, so it follows that
$|N_\bullet X_\U^{\text{ord}}|$ is weakly equivalent to $X$ via the
natural map.

\begin{Rem}
The construction above is related to topological category $X_\U$
defined by Segal~\cite{Segal:68} as follows: To an open cover $\U$
indexed by an unordered set $\I$ we can associate a new open cover
cover $\{U_S\}$, indexed over all finite non-empty subsets $S$ of
$\I$. There is a partial ordering of finite subsets given by
inclusion $R\subseteq S$. Even though $\{U_S\}$ does not satisfy the
conditions defining ``an ordered open covering'' the construction of
$X_{\{U_S\}}^{\text{ord}}$ is still well-defined. Moreover, the
category $X_{\{U_S\}}^{\text{ord}}$ is equal to $X_\U$. Explicitly;
the space of objects of $X_{\U}$ is the disjoint union of all $U_S$
and the space of morphisms is the disjoint union of $U_S$ over all
pairs $R\subseteq S$. The source map forgets $R$, and the target map
is induced from the inclusions $U_S\subseteq U_R$.
By~\cite[Proposition~2.7]{DuggerIsaksen:04}, $|N_\bullet X_\U|$ is
weakly equivalent to $X$ via the natural map.
\end{Rem}

\begin{Def}
Let $2X_\U^{\text{norm}}$ be the topological $2$-category
constructed from $X_\U^{\text{ord}}$ by adding identity
$2$-morphisms only.
\end{Def}

We now get the following reformulation of the definition of a
principal $2\C$-bundle subordinate to an ordered cover $\U$:

\begin{Prop}
There is a one-to-one correspondence between simplicial maps
$U_\bullet\rightarrow \Delta2\C$ and continuous normalized colax
$2$-functors $F\co 2X_\U^{\text{norm}}\rightarrow 2\C$.
\end{Prop}

\begin{proof}
Observe that $\Delta 2X_\U^{\text{norm}}\cong N_\bullet
X_\U^{\text{ord}}\cong U_\bullet$. Hence, a continuous normalized
colax $2$-functor $F$ yields a simplicial map $U_\bullet\rightarrow
\Delta2\C$ by functorality of the geometric nerve construction.

To reconstruct $F$ from a simplicial map
$U_\bullet\rightarrow\Delta2\C$, notice that the restrictions to
$0$- and $1$-morphisms, $U_0\rightarrow(\Delta2\C)_0$ and
$U_1\rightarrow(\Delta2\C)_1$, determine $F$ on objects and
$1$-morphisms respectively. Since $F$ is normalized, it takes
identities to identities, and there are no non-identity
$2$-morphisms in $2X_\U^{\text{norm}}$. The remaining data of the
colax $2$-functor $F$ are $2$-cells
\[
\phi\co F(x_{\alpha\gamma})\Rightarrow
F(x_{\beta\gamma})*F(x_{\alpha\beta}),
\]
and we define these to be the value of $U_2\rightarrow
(\Delta2\C)_2$ on the $2$-simplices $x_{\alpha\beta\gamma}$.

The verification of the cocylce condition on $3$-cells is omitted.
\end{proof}

If we restrict our attention to strict $2$-functors, the left
adjoint of $\Delta$ yields more complicated $2$-categories. We
define:

\begin{Def}\label{def:twoXU}
Define the topological $2$-category $2X_\U$ as follows: Let the
space of objects, $(2X_\U)_0$, be the disjoint union $\coprod
U_{\alpha}$. A typical object is denoted $x_\alpha$. Let the space
of $1$-morphisms be non-associatively freely generated by the
disjoint union $\coprod_{\alpha<\beta}U_{\alpha\beta}$. A typical
$1$-morphism is thus a parenthesized sequence of
$x_{\alpha_{i-1}\alpha_i}$'s. For example if
$\alpha<\beta<\gamma<\delta<\epsilon$ and $x\in
U_{\alpha\beta\gamma\delta\epsilon}$, then all the following
expressions are different $1$-morphisms:
\begin{multline*}
x_{\alpha\beta},\quad x_{\alpha\gamma},\quad
x_{\beta\gamma}*x_{\alpha\beta},\quad
(x_{\gamma\delta}*x_{\beta\gamma})*x_{\alpha\beta},\quad
x_{\gamma\delta}*(x_{\beta\gamma}*x_{\alpha\beta}),\quad
(x_{\delta\epsilon}*x_{\gamma\delta})*(x_{\beta\gamma}*x_{\alpha\beta}).
\end{multline*}
Observe that any $1$-morphism $f\co x_{\alpha_0}\rightarrow
x_{\alpha_k}$ has an associated chain
$\alpha_0<\alpha_1<\cdots<\alpha_k$ of indices, thus $f$ is some
parenthesizing of $x_{\alpha_{k-1}\alpha_k}* \cdots *
x_{\alpha_1\alpha_2} *x_{\alpha_0\alpha_{1}}$.

The space of $2$-morphisms $(2X_\U)_2$ is defined as follows: for
arbitrary $1$-morphisms $f\co x_{\alpha_0}\rightarrow x_{\alpha_k}$
and $g\co y_{\beta_0}\rightarrow y_{\beta_l}$ there are either one
or zero $2$-morphisms $f\Rightarrow g$. If $x=y$, $\alpha_0=\beta_0$
and $\alpha_k=\beta_l$, and the chain
$\alpha_0<\alpha_1<\cdots<\alpha_k$ is a refinement of
$\beta_0<\beta_1<\cdots<\beta_l$, then there exists a unique
$2$-morphism $f\Rightarrow g$. Otherwise, there is no $2$-morphism
$f\Rightarrow g$. Let $x_{\alpha\beta\gamma}$ denote the unique
$2$-morphism $x_{\alpha\gamma}\Rightarrow x_{\beta\gamma}*
x_{\alpha\beta}$.
\end{Def}

This leads to the following reformulation of
Definition~\ref{def:principal2Cbundle}:

\begin{Prop}
There is a one-to-one correspondence between simplicial maps
$U_\bullet\rightarrow \Delta2\C$ and continuous strict $2$-functors
$F\co 2X_\U\rightarrow 2\C$.
\end{Prop}

\begin{proof}
The proof is to inspect the definitions. Suppose first that we are
given such a strict $2$-functor $F$. Since $\Delta2\C$ is
$3$-coskeletal, a map $U_{\bullet}\rightarrow\Delta2\C$ is uniquely
determined by the maps of $0$-, $1$-, $2$-, and $3$-simplices.
Restricting $F$ to the objects, $1$-morphisms, and $2$-morphisms of
the form $x_{\alpha}$, $x_{\alpha\beta}$, and
$x_{\alpha\beta\gamma}$ we recover the maps $U_0\rightarrow (\Delta
2\C)_0$, $U_1\rightarrow (\Delta 2\C)_1$, and $U_2\rightarrow
(\Delta 2\C)_2$ respectively. To get the map on $3$-simplices, we
need to verify the tetrahedron coherence condition of $\Delta2\C$.
Take a $3$-simplex of $U_\bullet$. It is represented by a point $x$
in $U_{\alpha\beta\gamma\delta}$. Now inspect how $F$, in the
following diagram, transforms the top, consisting of $2$-morphisms
in $2X_{\U}$, into the bottom, consisting of $2$-morphisms in $2\C$:
\[
\xymatrix@C=7pt{ &
x_{\gamma\delta}*(x_{\beta\gamma}*x_{\alpha\beta})
\ar@{~>}[ddd]^(.4)F |!{[dl];[ddr]}\hole
\ar@{=>}[rr]^{\underline{\alpha}} &&
(x_{\gamma\delta}*x_{\beta\gamma})*x_{\alpha\beta} \ar@{~>}[ddd]_(.4)F |!{[ddl];[dr]}\hole &\\
x_{\gamma\delta}*x_{\alpha\gamma} \ar@{~>}[ddd]^F
\ar@{=>}[ur]^{1*x_{\alpha\beta\gamma}} &&&&
x_{\beta\delta}*x_{\alpha\beta} \ar@{~>}[ddd]_F
\ar@{=>}[ul]_{x_{\beta\gamma\delta}*1} \\
&& x_{\alpha\delta} \ar@{~>}[ddd]_F \ar@{=>}[ull]^{x_{\alpha\gamma\delta}} \ar@{=>}[urr]_{x_{\alpha\beta\delta}} \\
& E_{\gamma\delta}*(E_{\beta\gamma}*E_{\alpha\beta})
\ar@{=>}'[r][rr]^(.1){\underline{\alpha}} &&
(E_{\gamma\delta}*E_{\beta\gamma})*E_{\alpha\beta} &\\
E_{\gamma\delta}*E_{\alpha\gamma}
\ar@{=>}[ur]_{1*\phi_{\alpha\beta\gamma}} &&&&
E_{\beta\delta}*E_{\alpha\beta}.
\ar@{=>}[ul]^{\phi_{\beta\gamma\delta}*1} \\
&& E_{\alpha\delta} \ar@{=>}[ull]_{\phi_{\alpha\gamma\delta}}
\ar@{=>}[urr]^{\phi_{\alpha\beta\delta}} }
\]
This verifies the tetrahedron coherence condition.

Next, we consider the problem of reconstructing a strict $2$-functor
$F\co 2X_\U\rightarrow 2\C$ from a simplical map
$U_\bullet\rightarrow\Delta 2\C$. The map on $0$-, $1$-, and
$2$-simplices determine $F$ on objects, $1$-morphisms, and
$2$-morphisms of the form $x_{\alpha}$, $x_{\alpha\beta}$, and
$x_{\alpha\beta\gamma}$. Now observe that $2X_{\U}$ is generated by
such $1$- and $2$-morphisms, and the associators. Thus the strict
functor $F$ is uniquely determined by the data above. We omit
further details.
\end{proof}

\begin{Prop}
Let $\U$ be an ordered open cover of $X$. The corresponding ordered
\v{C}ech complex $U_{\bullet}$ embeds as a simplicial deformation
retract of $\Delta2X_{\U}$, and consequently we get a homotopy
equivalence $|\Delta 2X_\U|\simeq X$.
\end{Prop}

We know that $U_\bullet\cong\Delta2X_\U^{\text{norm}}$. If
$2X_\U^{\text{norm}}$ and $2X_{\U}$ are biequivalent as
$2$-categories, then the result would follow
from~\cite[Section~10]{Street:96}. However, a direct argument is not
difficult:

\begin{proof}
Observe that $\Delta2X_{\U}$ actually is $2$-coskeletal, and we can
write any $n$-simplex as a flag
\[
\begin{pmatrix}
x_{\alpha_0} & f_{01} & f_{02} & f_{03} & \cdots & f_{0n} \\
& x_{\alpha_1} & f_{12} & f_{13} & \cdots & f_{1n}\\
&& x_{\alpha_2} & f_{23} & \cdots & f_{2n}\\
&&& x_{\alpha_3} & \cdots & f_{3n}\\
&&&& \ddots & \vdots \\
&&&&& x_{\alpha_n}
\end{pmatrix}_,
\]
where $\alpha_0<\alpha_1<\ldots<\alpha_n$, $x\in
U_{\alpha_0\cdots\alpha_n}$ and $f_{ij}\co x_{\alpha_i}\rightarrow
x_{\alpha_j}$ are $1$-morphisms in $2X_{\U}$ such that for every
$i<j<k$ there exists a $2$-morphism $f_{ik}\Rightarrow
f_{jk}*f_{ij}$. Define $U_{\bullet}\rightarrow\Delta2X_{\U}$ by
sending the $n$-simplex $x\in U_{\alpha_0\alpha_1\cdots\alpha_n}$ to
the $n$-simplex in $\Delta2X_{\U}$ given by
$f_{ij}=x_{\alpha_i\alpha_j}$. Observe that this gives an inclusion,
and we view $U_{\bullet}$ as a sub-simplicial space of
$\Delta2X_{\U}$.

Now we define a simplicial deformation retraction $h$ from
$\Delta2X_{\U}$ down to $U_{\bullet}$ by
\[
h_k(f_{ij})=
\begin{pmatrix}
x_{\alpha_0} & \cdots & x_{\alpha_0\alpha_{k-1}} &
x_{\alpha_0\alpha_{k}}
    & x_{\alpha_0\alpha_{k}} & x_{\alpha_0\alpha_{k+1}} &\cdots & x_{\alpha_0\alpha_{n}}\\
& \ddots & \vdots & \vdots & \vdots & \vdots & & \vdots \\
&& x_{\alpha_{k-1}} & x_{\alpha_{k-1}\alpha_k} &
x_{\alpha_{k-1}\alpha_k} & x_{\alpha_{k-1}\alpha_{k+1}}
    & \cdots & x_{\alpha_{k-1}\alpha_n} \\
&&& x_{\alpha_k} & 1_{x_{\alpha_k}} & x_{\alpha_k\alpha_{k+1}} & \cdots & x_{\alpha_k\alpha_n} \\
&&&& x_{\alpha_k} & f_{k\,k+1} & \cdots & f_{kn}\\
&&&&& x_{\alpha_{k+1}} & \cdots & f_{k+1\,n}\\
&&&&&& \ddots & \vdots\\
&&&&&&& x_{\alpha_n}
\end{pmatrix}_.
\]
These maps $h_k\co (\Delta2X_{\U})_n\rightarrow
(\Delta2X_{\U})_{n+1}$, for $0\leq k \leq n$, verify the usual
identities for a simplicial deformation retraction.
\end{proof}

Given a principal $2\C$-bundle $\E$ over $X$ and a subspace $A$ of
$X$ it is clear that we may restrict $\E$ to $A$. We denote the
restriction by $\E|_A$. More generally, there is a pullback $f^*\E$
over $Y$, for every continuous map $f\co Y\rightarrow X$.

We define equivalence of principal $2\C$-bundles by concordance:

\begin{Def}
Two principal $2\C$-bundles $\E_0$ and $\E_1$ over $X$ are
\emph{concordant} if there exists a principal $2\C$-bundle $\E$ over
$X\times I$ such that the restrictions $\E|_{X\times\{0\}}$ and
$\E|_{X\times\{0\}}$ are $\E_0$ and $\E_1$ respectively.
\end{Def}

Concordance is an equivalence relation. The main theorem,
Theorem~\ref{thm:main}, classifies the set of concordance classes of
principal $2\C$-bundles over $X$: under mild assumptions on $X$ and
$2\C$, there is a bijection between concordance classes and homotopy
classes of maps $X\rightarrow |\Delta 2\C|$.

\begin{Rem}
There are other ways to define an equivalence relation between
principal $2\C$-bundles. In all the following cases one should allow
elementary refinements of the cover,
compare~\cite[Definition~2.4]{BaasDundasRognes:04}. In addition one
declares two principal $2\C$-bundles $\E_0$ and $\E_1$, subordinate
to the same ordered open cover $\U$, to be equivalent if either
\begin{itemize}
\item[i)] there exists a natural transformation between the $2$-functors
$2X_\U\rightarrow 2\C$ corresponding to $\E_0$ and $\E_1$,
\item[ii)] there is an adjoint equivalence relating $\E_0$ and
$\E_1$ in the $2$-category of $2$-functors, natural transformations
and modifications, or
\item[iii)] there is a simplical homotopy $U_\bullet\times
\Delta^1_\bullet\rightarrow \Delta2\C$ between the simplicial maps
$U_\bullet\rightarrow\Delta2\C$ representing $\E_0$ and $\E_1$
respectively.
\end{itemize}
This list of possible equivalence relations is not exhaustive. Under
what circumstances will this actually lead to a set of equivalence
classes different from the concordance classes? Discussing this in
detail would be beyond the scope of this article, but we believe
that for compact $X$ and those $2\C$ that occur most commonly, the
equivalences of type i)---iii) should all coincide with concordance.
\end{Rem}

\section{Examples}\label{sect:examples}

In this section we give several examples of structure
$2$-categories, their principal $2$-bundles, and the corresponding
classifying space. Let us start off by discussing a simple example
which illuminates the fact that different categories can have the
same concordance classes of principal bundles:

\begin{Exa}[(Complex line bundles)]
In this example we consider the more elementary setting of
topological categories. In this case the nerve, $N_\bullet$, and the
geometric nerve, $\Delta$, coincide, and we may use Quillen's
Theorem~A,~\cite{Quillen:73}, for calculations.

Let $U(1)$ be the topological group of complex numbers with modulus
$1$. Needless to say, a principal bundle with structure group $U(1)$
corresponds to a complex line bundle. If we want to think about
$U(1)$ as a topological category, it has only one object, $*$, and
$U(1)$ is the space of morphisms. We may also think about the space
$\mathbb{CP}^\infty$ as a category: The space of objects is
$\mathbb{CP}^\infty$, and we add identity morphisms only. It is
classical that maps $X\rightarrow \mathbb{CP}^\infty$ corresponds to
a complex line bundle by pulling back the canonical line bundle.

There is an intermediate topological category, which we will denote
by $\mathcal{CU}$. The space of objects is $\mathbb{CP}^\infty$. The
morphisms $\mathcal{CU}(L,L')$ between two complex lines in
$\mathbb{C}^\infty$ are the unitary linear transformations $f\co
L\rightarrow L'$. To topologize the total space of all morphisms we
consider basis vectors $b$ and $b'$ of length $1$ for $L$ and $L'$
respectively, together with a point $c\in U(1)$. To the triple
$(b,b',c)$ we associate the transformation $f\co zb\mapsto zcb'$.
Therefore $\mathcal{CU}_1$ is defined as the quotient of all such
triples, where we identify $(b,b',c)$ with $(\theta b,\theta
b',\theta^{-1}\theta'c)$ for all $\theta,\theta'\in U(1)$.

Fix a complex line $L_0$ in $\mathcal{CU}_0$. This gives an
inclusion functor $i\co U(1)\rightarrow\mathcal{CU}$ by viewing
$c\in U(1)$ as an operator on $L_0$. For every object $L$ of
$\mathcal{CU}$, the comma category $L\backslash i$ consists of
terminal objects only, and hence $L\backslash i$ is contractible. By
Quillen's Theorem~A it follows that $i$ induces a weak equivalence
\[
BU(1)=|N_\bullet
U(1)|\xrightarrow{\simeq}|N_\bullet\mathcal{CU}|=B\mathcal{CU}.
\]

There is also an inclusion $j\co \mathbb{CP}^\infty\rightarrow
\mathcal{CU}$ given as the identity map on objects. Let $L$ be an
arbitrary object of $\mathcal{CU}$. Observe that all morphisms of
the comma category $L\backslash j$ are identity morphisms, and that
the space of objects can be identified with the unit sphere in
$\mathbb{C}^\infty$. Consequently, $B(L\backslash j)$, which in this
case is equal to the space of objects, is contractible, and by
Quillen's Theorem~A it follows that $j$ induces a weak equivalence
\[
\mathbb{CP}^\infty=|N_\bullet\mathbb{CP}^\infty|
\xrightarrow{\simeq} |N_\bullet\mathcal{CU}|=B\mathcal{CU}.
\]
\end{Exa}

\begin{Exa}[(Gerbes)]
Recall from~\cite[Section~3.1]{MurrayStevenson:00} the notion of a
\emph{local bundle gerbe}. This fits into our framework as follows.
Let $2\G$ be the topological $2$-category with $2\G_0=*$,
$2\G_1=\mathbb{CP}^{\infty}$, and interpreting
$L,L'\in\mathbb{CP}^\infty$ as complex lines in $\mathbb{C}^\infty$,
a $2$-morphism $f\co L\Rightarrow L'$ is the same as a unitary
linear transformation between these lines. Represent the lines as
$L=[z_1,z_2,z_3,\ldots]$ and $L'=[z'_1,z'_2,z'_3,\ldots]$. Then we
may write their tensor product as line represented by the collection
of all products $z_iz'_j$, i.e. $L\otimes
L'=[z_1z'_1,z_2z'_1,z_1z'_2,z_1z'_3,z_2z'_2,z_3z'_1,\ldots]$. This
gives a map
$\mathbb{CP}^{\infty}\times\mathbb{CP}^{\infty}\rightarrow\mathbb{CP}^{\infty}$,
and which we take as the definition of horizontal composition of
$1$-morphisms.

Now consider a principal $2\G$-bundle. Since $\mathbb{CP}^{\infty}$
classifies complex line bundles, we interpret $E_{\alpha\beta}$ as a
complex line bundle over $U_{\alpha\beta}$. From this point of view,
$\phi_{\alpha\beta\gamma}$ is a bundle isomorphism from
$E_{\alpha\gamma}$ to $E_{\beta\gamma}\otimes E_{\alpha\beta}$, and
it corresponds to the trivialization $\theta_{\alpha\beta\gamma}$
considered by~\cite{MurrayStevenson:00}. There is a one-to-one
correspondence between stable isomorphism classes of gerbes over $X$
and elements of the cohomology group $H^3(X;\Z)$,
see~\cite[Theorem~3.3]{MurrayStevenson:00}. Thus we write
\[
\operatorname{Gerbe}(X)=[X,K(\Z,3)].
\]
\end{Exa}

\begin{Exa}[(BDR $2$-vector bundles)]\label{exa:2KVn}
In order to classify the concordance classes of $2$-vector bundles
it is sufficient to describe the corresponding structure
$2$-category. By Theorem~\ref{thm:main} the geometric nerve of the
structure $2$-category classifies the concordance classes. Since
this result depends only on transition data, the description of the
fiber plays no role here.

Based on the $2$-vector spaces of~\cite{KapranovVoevodsky:94}, we
define a topological $2$-category $2\mathcal{KV}(n)$ as
follows: Take a single point $*$ as the space of objects, and let
the $1$-morphisms be the space of $n\times n$-matrices of complex
Grassmannians
\[
\coprod \big( \operatorname{Gr}(a_{ij})\big)_{i,j=1}^n,
\]
where the disjoint union runs over all invertible $n\times
n$-matrices $(a_{ij})$ consisting of non-negative integers. Thus a
map $E_{\alpha\beta}\co U_{\alpha\beta}\rightarrow
(2\mathcal{KV}(n))_1$ defines an $n\times n$-matrix
$(E_{ij})_{\alpha\beta}$ of complex vector bundles such that
$\det(\dim (E_{ij})_{\alpha\beta})=\pm 1$. The $2$-morphisms are
defined as matrices of linear isomorphisms, and the composition,
$*$, of $1$-morphisms is defined as reverse matrix multiplication
using the exterior $\otimes$ of Grassmanians as product between
matrix entries, and $\oplus$ as sum. This recovers the definition of
\emph{charted $2$-vector bundles} given
in~\cite{BaasDundasRognes:04}. Their composition $\cdot$ is dual to
our $*$, and our $\phi_{\alpha\beta\gamma}$ is the inverse of their
$\phi^{\alpha\beta\gamma}$. We write
\[
2\operatorname{Vect}_n(X)=[X,|\Delta2\mathcal{KV}(n)|].
\]

This leads to a geometrically defined form of elliptic cohomology:
Baas, Dundas and Rognes show in~\cite{BaasDundasRognes:04} how to
define a kind of second order $K$-theory, namely the $2K$-theory of
$2$-vector bundles. They conjectured that the representing spectrum
was $K(ku)$, and this has been proved by Baas, Dundas, Richter
and Rognes in~\cite{BaasDundasRichterRognes:07}. Calculations by
Ausoni and Rognes~\cite{Ausoni:06} show that this spectrum has
telescopic complexity $2$, and hence it qualifies as a form of elliptic
cohomology.
\end{Exa}

\begin{Exa}[(String bundles)]
Let $2\mathscr{S}$ denote the surface $2$-category defined by U.
Tillmann in~\cite[Section~2]{Tillmann:97}. Its objects,
$1$-morphisms, and $2$-morphisms are closed strings, oriented
cobordisms, and their diffeomorphisms respectively. Define
\emph{string bundles} as principal $2\mathscr{S}$-bundles. The
classifying space $|\Delta2\mathscr{S}|$ is an infinite loop space,
and it is related to the stabilized mapping class group of oriented
surfaces, $\Gamma_\infty$, by the formula:
\[
\Omega|N_\bullet2\mathscr{S}|\cong\Z\times B\Gamma_\infty^+,
\]
see~\cite[Theorem~3.1]{Tillmann:97}.

Let $\mathbb{CP}_{-1}^\infty$ denote the spectrum whose $2k$th space
is the Thom space of the complex bundle complementary to the
canonical complex line bundle over $\mathbb{CP}^{k-1}$. By I. Madsen
and M. Weiss' proof of the generalized Mumford
conjecture~\cite{MadsenWeiss:07}, the classifying space of string
bundles is $\Omega^{\infty-1}\mathbb{CP}_{-1}^\infty$.
In~\cite[Definition~5.1]{GalatiusMadsenTillmannWeiss:06} S.
Galatius, I. Madsen, U. Tillmann, and M. Weiss introduce a
topological category $\EuScript{C}_2^+$ of embedded oriented
surfaces. By their new proof of the generalized Mumford conjecture,
the notion of principal $\EuScript{C}_2^+$-bundles coincides with
the notion of string bundles, up to concordance. They also introduce
the notation $MTSO(2)$ for the spectrum $\mathbb{CP}_{-1}^{\infty}$.
To summarize, the concordance classes of string bundles are
classified by $\Omega^{\infty-1}MTSO(2)=
\Omega^{\infty-1}\mathbb{CP}_{-1}^\infty$, i.e.
\[
\operatorname{String}(X)=[X,\Omega^{\infty-1}MTSO(2)]
=[X,\Omega^{\infty-1}\mathbb{CP}_{-1}^\infty].
\]
\end{Exa}

\begin{Exa}[(Cobordism bundles)]
Extending the previous example, we now define \emph{$d$-cobordism
bundles} as principal $\EuScript{C}_d$-bundles, where
$\EuScript{C}_d$ is the topological category of embedded
$(d-1)$-manifolds and embedded cobordisms introduced
in~\cite[Section~2.1]{GalatiusMadsenTillmannWeiss:06}. Furthermore,
let $MTO(d)$ be the spectrum whose $(n+d)$th space is the Thom
spectrum of the bundle $U_{d,n}^\perp$ complementary to the
canonical $d$-plane bundle over the Grassmannian of $d$-planes in
$\R^{n+d}$. The main theorem of Galatius, Madsen, Tillmann, and
Weiss shows that the classifying space of $d$-cobordism bundles is
$\Omega^{\infty-1}MTO(d)$. There are variations of these
constructions by placing various types of tangential structure on
the cobordism category, e.g. orientation, etc. We write
\[
\operatorname{Cob}^d(X)=[X,\Omega^{\infty-1}MTO(d)].
\]
\end{Exa}

\begin{Exa}[(Graph bundles)]
Let $\mathcal{C}_\infty$ denote the graph cobordism category defined
by Galatius in~\cite[Section~3]{Galatius:06}. We define \emph{graph
bundles} as principal $\mathcal{C}_\infty$-bundles. On one hand, the
classifying space $B\mathcal{C}_\infty$ is weakly equivalent to
$\Omega^{\infty-1}S$, where $S$ is the sphere spectrum. On the other
hand, there is a weak equivalence $\Z\times
B\!\operatorname{Aut}^+_\infty\simeq \Omega B\mathcal{C}_\infty$,
where $\operatorname{Aut}_\infty$ is the automorphism group of the
free group on $n$ generators as $n$ tends to infinity. For the
concordance classes of graph bundles over $X$ we write
\[
\operatorname{Graph}(X)=[X,\Omega^{\infty-1}S].
\]
\end{Exa}

\begin{Exa}[(BC $2$-vector bundles)]\label{exa:BC2cat}
Let us now discuss choices of $2\C$ based on Baez and Crans'
$2$-vector spaces. There are numerous variations depending on how
strictly the $1$-morphisms should be equivalences, how much
structure the $2$-morphisms should contain, and how to choose the
objects from $2$-vector spaces within the same weak equivalence
class. Let $b_1$ and $b_0$ be non-negative integers and let $C_*$ be
a $2$-term chain complex of vector spaces.

\begin{description}
\item[i)] Let $2\mathcal{B}_{\text{strict}}(C_*)$ be the $2$-category
with a single object $C_*$, the $1$-morphisms are the chain
isomorphisms $f\co C_*\rightarrow C_*$, and the $2$-morphisms
$f\Rightarrow g$ are chain homotopies $\phi$ from $f$ to $g$.
\item[ii)] Let $2\mathcal{B}_{\text{weak}}(C_*)$ be the $2$-category
with a single object $C_*$, the $1$-morphisms are the chain
equivalences $f\co C_*\rightarrow C_*$, and the $2$-morphisms
$f\Rightarrow g$ are chain homotopies $\phi$ from $f$ to $g$.
\item[iii)] Let $2\mathcal{B}_{\text{eq}}(C_*)$ be the $2$-category
with a single object $C_*$, the $1$-morphisms are tuples
$(f,\bar{f},\iota_f,\epsilon_f)$, where $f$ and $\bar{f}$ are chain
maps $C_*\rightarrow C_*$, and $\iota_f\co 1\Rightarrow f\bar{f}$
and $\epsilon_f\co \bar{f}f\Rightarrow 1$ are chain homotopies, and
the $2$-morphisms $(f,\bar{f},\iota_f,\epsilon_f)\Rightarrow
(g,\bar{g},\iota_g,\epsilon_g)$ are pairs of chain homotopies
$\phi\co f\Rightarrow g$ and $\bar{\phi}\co
\bar{f}\Rightarrow\bar{g}$ such that the following identities hold:
\[
\iota_g-\iota_f=f\bar{\phi}+\phi\bar{g} \quad\text{and}\quad
\epsilon_f-\epsilon_g=f\bar{\phi}+\phi\bar{g}.
\]
\item[iv)] Let $2\mathcal{B}_{\text{ad}}(C_*)$ be the $2$-category
with a single object $C_*$, the $1$-morphisms are tuples
$(f,\bar{f},\iota_f,\epsilon_f)$, where $f$ and $\bar{f}$ are chain
maps $C_*\rightarrow C_*$, and $\iota_f\co 1\Rightarrow f\bar{f}$
and $\epsilon_f\co \bar{f}f\Rightarrow 1$ are chain homotopies such
that the zigzag identities
\[
f\epsilon_f+\iota_f f=0 \quad\text{and}\quad
\epsilon_f\bar{f}+\bar{f}\epsilon_f=0
\]
hold, and the $2$-morphisms
$(f,\bar{f},\iota_f,\epsilon_f)\Rightarrow
(g,\bar{g},\iota_g,\epsilon_g)$ are pairs of chain homotopies
$\phi\co f\Rightarrow g$ and $\bar{\phi}\co
\bar{f}\Rightarrow\bar{g}$ such that the following identities hold:
\[
\iota_g-\iota_f=f\bar{\phi}+\phi\bar{g} \quad\text{and}\quad
\epsilon_f-\epsilon_g=f\bar{\phi}+\phi\bar{g}.
\]
\item[v)]\label{item:Bweak} Let $2\mathcal{B}_{\text{weak}}(b_1,b_0)$ be the $2$-category
with objects all $2$-term chain complexes $C_*$ with Betti-numbers
$b_0,b_1$, the $1$-morphisms are the chain equivalences $f\co
C_*\rightarrow C'_*$, and the $2$-morphisms $f\Rightarrow g$ are
chain homotopies $\phi$ from $f$ to $g$.

\item[vi)] Let $2\mathcal{B}_{\text{eq}}(b_1,b_0)$ be the $2$-category
with objects all $2$-term chain complexes $C_*$ with Betti-numbers
$b_0,b_1$, the $1$-morphisms are tuples
$(f,\bar{f},\iota_f,\epsilon_f)$, where $f\co C_*\rightarrow C'_*$
and $\bar{f}\co C'_*\rightarrow C_*$ are chain maps, and $\iota_f\co
1\Rightarrow f\bar{f}$ and $\epsilon_f\co \bar{f}f\Rightarrow 1$ are
chain homotopies, and the $2$-morphisms
$(f,\bar{f},\iota_f,\epsilon_f)\Rightarrow
(g,\bar{g},\iota_g,\epsilon_g)$ are pairs of chain homotopies
$\phi\co f\Rightarrow g$ and $\bar{\phi}\co
\bar{f}\Rightarrow\bar{g}$ such that the following identities hold:
\[
\iota_g-\iota_f=f\bar{\phi}+\phi\bar{g} \quad\text{and}\quad
\epsilon_f-\epsilon_g=f\bar{\phi}+\phi\bar{g}.
\]
\item[vii)] Let $2\mathcal{B}_{\text{ad}}(b_1,b_0)$ be the $2$-category
with objects all $2$-term chain complexes $C_*$ with Betti-numbers
$b_0,b_1$, the $1$-morphisms are tuples
$(f,\bar{f},\iota_f,\epsilon_f)$, where $f\co C_*\rightarrow C'_*$
and $\bar{f}\co C'_*\rightarrow C_*$ are chain maps, and $\iota_f\co
1\Rightarrow f\bar{f}$ and $\epsilon_f\co \bar{f}f\Rightarrow 1$ are
chain homotopies such that the zigzag identities
\[
f\epsilon_f+\iota_f f=0 \quad\text{and}\quad
\epsilon_f\bar{f}+\bar{f}\epsilon_f=0
\]
hold, and the $2$-morphisms
$(f,\bar{f},\iota_f,\epsilon_f)\Rightarrow
(g,\bar{g},\iota_g,\epsilon_g)$ are pairs of chain homotopies
$\phi\co f\Rightarrow g$ and $\bar{\phi}\co
\bar{f}\Rightarrow\bar{g}$ such that the following identities hold:
\[
\iota_g-\iota_f=f\bar{\phi}+\phi\bar{g} \quad\text{and}\quad
\epsilon_f-\epsilon_g=f\bar{\phi}+\phi\bar{g}.
\]
\end{description}
\end{Exa}

We will now calculate the classifying space in the case~\textbf{v)}
of the example above. Set $2\B=2\B_{\text{weak}}(b_1,b_0)$. In order
to study $2\B$ in more detail, we choose for each object $C_*$
isomorphisms
\[
C_1\cong H_1\oplus B\quad\text{and}\quad C_0\cong B\oplus H_0,
\]
where $H_1$ and $H_0$ are complex vector spaces of complex dimension
$b_1$ and $b_0$ respectively, and $B$ is an arbitrary complex vector
space. Moreover, we can choose these isomorphisms such that the
differential $d\co C_1\rightarrow C_0$ is identified with the matrix
$\begin{pmatrix}0&I\\0&0\end{pmatrix}$. A $1$-morphism $f\co
C_*\rightarrow C'_*$ is then identified with a pair of matrices,
\[
\begin{pmatrix}f_{H'_1}^{H_1}&f_{H'_1}^{B}\\0&f_{B'}^B\end{pmatrix}
\quad\text{and}\quad
\begin{pmatrix}f_{B'}^B&f^{H_0}_{B'}\\0&f_{H'_0}^{H_0}\end{pmatrix}
\]
representing $f_1$ and $f_0$ respectively. Here the notation
$f_{B'}^{H_1}$ indicates a linear map from $H_1$ to $B'$, etc.
Observe that $f$ is a chain equivalence if and only if both
$f_{H'_1}^{H_1}$ and $f_{H'_0}^{H_0}$ are invertible. A chain
homotopy $\phi$ from $f$ to $g$ is a linear map $C_0=B\oplus
H_0\rightarrow H'_1\oplus B'=C'_1$, and thus represented by a matrix
\[
\begin{pmatrix}\phi^B_{H'_1}&\phi^{H_0}_{H'_1}\\
\phi^{B}_{B'}&\phi^{H_0}_{B'}\end{pmatrix},
\]
such that
\[
\phi^{B}_{H'_1}=f^{B}_{H'_1}-g^{B}_{H'_1},\quad
\phi^{B}_{B'}=f^{B}_{B'}-g^{B}_{B'},\quad\text{and}\quad
\phi^{H_0}_{B'}=f^{H_0}_{B'}-g^{H_0}_{B'}.
\]
Horizontal composition corresponds to matrix multiplication and
vertical composition to addition of matrices.

Because $\phi^{H_0}_{H'_1}$ is arbitrary the space of chain
homotopies between two fixed chain maps $f$ and $g$ admits an affine
structure. This indicates that there are too many $2$-morphisms in
$2\B$. Therefore, define $2\B^{\Ho}$ as the quotient of $2\B$, where
the objects and $1$-morphisms are the same, and two $2$-morphisms
from $2\B$ are identified in the quotient if their source and target
$1$-morphisms are pairwise equal. Explicitly, the class of $\phi$ in
$2\B^{\Ho}$ is represented by a matrix
\[
\begin{pmatrix}\phi^B_{H'_1}&*\\
\phi^{B}_{B'}&\phi^{H_0}_{B'}\end{pmatrix},
\]
where we do not care about the stuff in the upper right corner. Not
unexpectedly, we get:

\begin{Lem}
The geometric nerves $\Delta2\B$ and $\Delta2\B^{\Ho}$ are weakly
equivalent.
\end{Lem}

\begin{proof}
Here is one way to show this: both $2$-categories satisfies the
assumtions~\textbf{G1'}, \textbf{G2'}, \textbf{G3}, and~\textbf{G4}
of Section~\ref{sect:GSF}, and the canonical map
$2\B\rightarrow2\B^{Ho}$ is a continuous strict $2$-functor. In
order to apply Lemma~\ref{lem:DeltaSimplified}, we fix chain
compexes $C_*$ and $C'_*$, and compare nerves of morphism
categories:
\[
N_\bullet2\B(C_*,C'_*)\rightarrow N_\bullet2\B^{\Ho}(C_*,C'_*).
\]
At each simplicial degree $k$, this map is a fiber bundle with fiber
an affine space, namely a $(k-1)$-tuple of $\phi^{H_0}_{H_1}$s.
Since the fiber is contractible, the map is a weak equivalence.
\end{proof}

We finish the calculation by comparing $2\B^{\Ho}$ to the product of
the complex general linear groups $GL(b_1)$ and $GL(b_0)$. We have
$2$-functors
\[
GL(b_1)\times GL(b_0)\xrightarrow{i} 2\B^{\Ho}
\xrightarrow{H}GL(b_1)\times GL(b_0),
\]
where $i$ includes $GL(b_1)\times GL(b_0)$ as the automorphisms of
the chain complex $\bbC^{b_1} \xrightarrow{0} \bbC^{b_0}$, and $H$
is homology of the chain complex $C_*$, i.e. we choose bases for
$H_1$ and $H_0$ for each object of $2\B^{\Ho}$. For $\bbC^{b_1}$ and
$\bbC^{b_0}$ we choose the standard bases. Observe that the
composition $Hi$ is the identity on $GL(b_1)\times GL(b_0)$. We will
now construct a normalized colax transformation $\sigma$ from the
identity on $2\B^{\Ho}$ to the composition $iH$.

The choice of bases for $H_1$ and $H_0$ specifies projection maps
$C_1\cong H_1\oplus B\rightarrow \bbC^{b_1}$ and $C_0\cong B\oplus
H_0\rightarrow \bbC^{b_0}$. This yields for each chain complex $C_*$
a chain map $\sigma_{C_*}$ into
$\bbC^{b_1}\xrightarrow{0}\bbC^{b_0}$. For a chain equivalence $f\co
C_*\rightarrow C'_*$ we observe that
\[
\xymatrix@C=2.5cm{ H_1\oplus B \ar[d]^{\sigma_{C_*}}
\ar[r]^{\scriptsize
\begin{pmatrix}f_{H'_1}^{H_1}&f_{H'_1}^{B}\\0&f_{B'}^B\end{pmatrix}}
& H'_1\oplus B' \ar[d]^{\sigma_{C'_*}}\\
\bbC^{b_1} \ar[r]^{[f^{H_1}_{H'_1}]} & \bbC^{b_1}}
\]
does not commute. Therefore we define $\sigma_f$ as the class of
chain homotopies represented by the matrix
$\begin{pmatrix}f^B_{H'_1} & *\end{pmatrix}$. We leave it to the
reader to verify that $\sigma$ is a normalized colax transformation.

\begin{Rem}
To see why we need the intermediate $2$-category $2\B^{\Ho}$, try to
replace $2\B^{\Ho}$ with $2\B$, and then define $\sigma$ as above.
In that case take $\sigma_f=\begin{pmatrix} f^B_{H'_1} & 0
\end{pmatrix}$. Given a pair of composable chain maps
$C_*\xrightarrow{f}C'_*\xrightarrow{g}C''_*$, observe that, in
general, we have $\sigma_{gf}\neq g\circ \sigma_f + \sigma_g\circ f$
in $2\B$. Hence $\sigma$ is not well-defined on $2\B$.
\end{Rem}

We now apply the following standard technique:

\begin{Lem}
Let $F$ and $G$ be normalized colax $2$-functors $2\C\rightarrow
2\D$. A normalized colax transformation $\phi\co F\rightarrow G$
induces a homotopy between the induced maps $|\Delta F|$ and
$|\Delta G|$.
\end{Lem}

\begin{proof}
Let $\I$ be the category with two objects, $0$ and $1$, and one
non-identity morphisms, $0\rightarrow 1$. Our normalized colax
transformation $\phi$ defines a $2$-functor $2\C\times\I\rightarrow
2\D$, which restricts to $F$ and $G$ over $0$ and $1$ respectively.
Since the geometric nerve of $\I$ is an interval, we get our
homotopy.
\end{proof}

This lemma immediately implies that the geometric nerve of
$2\B^{\Ho}$ is weakly equivalent to the geometric nerve of
$GL(b_1)\times GL(b_0)$. We arrive at:

\begin{Concl}\label{concl:BC}
The geometric nerve of $2\B=2\B_{\text{weak}}(b_1,b_0)$ is weakly
equivalent to the nerve $N_\bullet GL(b_1)\times N_\bullet GL(b_0)$.
Hence, the corresponding Baez--Crans $2$-vector bundles splits, up
to concordance, as two copies of ordinary vector bundles. With
little or no modification, similar calculations can be carried out
for the other $2$-categories of Example~\ref{exa:BC2cat}. In all
cases the geometric nerve splits, up to weak equivalence, as
$N_\bullet GL(b_1)\times N_\bullet GL(b_0)$. We can therefore
conclude that $2K$-theories based on Baez--Crans $2$-vector bundles
are given as two copies of ordinary $K$-theory.
\end{Concl}





\section{The topological Kan condition}\label{sect:Kan}

We consider a simplicial space $Z_\bullet$ with a topological Kan
property. We apply the following notation: let $\Delta^n$ be the
topological $n$-simplex, i.e. the subspace of $\R^{n+1}$ consisting
of tuples $(t_0,t_1,\ldots,t_n)$ such that each $t_i\geq0$ and $\sum
t_i=1$. However, $\Delta^n_\bullet$ denotes the simplicial
$n$-simplex, i.e. the simplicial set with $q$-simplices all order
preserving maps $[q]\rightarrow[n]$. We define the \emph{space of
$k$-horns of dimension $n$} in the simplicial space $Z_\bullet$ to
be
\begin{gather}
\Lambda^n_k(Z_\bullet)=\nonumber\\
\{(z_0,\dots,z_{k-1},-,z_{k+1},\dots,z_n)\mid
d_i(z_j)=d_{j-1}(z_i),i<j,i\not=k\not=j\}\nonumber\\
\subset \underbrace{Z_{n-1}\times Z_{n-1}\times\cdots \times
Z_{n-1}} _{n\text{ factors}}\nonumber
\end{gather}
We give $\Lambda^n_k(Z_\bullet)$ the subspace topology. For our
purposes, it seems more natural to associate this horn to the
complementary set $I(k)=\{0,1,\dots,k-1,k+1,\dots n\}$. So we think
of a $k$-horn as an $I(k)$-cohorn. If $I\subset
[n]=\{0,1,\dots,n\}$, we define the \emph{space of $I$-cohorns} to
be
\begin{gather}
\Lambda_I^n(Z_\bullet)=\nonumber\\
\{(z_i)_{i\in I}\mid
d_i(z_j)=d_{j-1}(z_i),i<j,i\in I\ni j\}\nonumber\\
\subset \underbrace{ Z_{n-1}\times Z_{n-1}\times\cdots \times
Z_{n-1} } _{\abs I \text{ factors}}\nonumber
\end{gather}
Note that $\Lambda_k^n=\Lambda_{I(k)}^n$. The face maps $d_i\co
Z_n\rightarrow Z_{n-1}$, $i\neq k$ induces a canonical map
$c_{k}^n\colon Z_{n}\to \Lambda^n_k(Z_\bullet)$, and more generally,
for any $I\subset [n]$ we have a canonical map $c_I^n\colon Z_{n}\to
\Lambda_I^n(Z_\bullet)$.

\begin{Def}\label{Def:Kanconditions}
We say that $Z_\bullet$ satisfies the \emph{discrete Kan condition}
if the maps $c^n_k$ are surjective for all $n$ and $0\leq k\leq n$.
We say that $Z_\bullet$ satisfies the \emph{topological Kan
condition} if the maps $c^n_k$ are surjective Serre fibrations for
all $n$ and $0\leq k\leq n$.
\end{Def}

\begin{Rem}
E.H. Brown and R.H. Szczarba consider
in~\cite[Definition~2.1]{BrownSzczarba:95} a related notion of
fibrancy.
\end{Rem}

\begin{Lem}\label{le:Serre}
If the simplicial space $Z_\bullet$ satisfies the topological Kan
condition, then for all $n$ and every non-empty proper subset $I$ of
$[n]$ the map $c_I^n$ is a surjective Serre fibration.
\end{Lem}

\begin{proof}
The proof is by double induction. First by induction on $n$: assume
that all $c_{I'}^{n'}$ are surjective Serre fibrations for $n'<n$
and every non-empty proper subset $I'$ of $[n']$. Next by downwards
induction on $\abs I=m$: assume that all $c_J^n$ are surjective
Serre fibrations for all non-empty proper subsets $J$ of $[n]$ with
$\abs J>m$. This is true for $m=n-2$ by assumption, since the $J$'s
in this case have the form $I(k)$ and $c_J^n=c_{I(k)}^n=c_k^n$. For
the inductive step assume that $k\not\in I$ and $J=I\cup \{k\}$, so
that by the induction hypothesis $c_J^n$ is a surjective Serre
fibration. Consider the canonical map which forgets $z_k$
\[
r^J_I\colon \Lambda_J^n(Z_\bullet) \to \Lambda_I^n(Z_\bullet).
\]
Then $c^n_I=r^J_I \circ c^n_J$. So it is enough to show that the map
$r_I^J$ is a surjective Serre fibration.

By $\delta_k\colon [n-1]\to [n]$ we denote the injective map which
omits $k$, thus the face map $d_k\colon Z_n\to Z_{n-1}$ is the map
induced by $\delta_k$. Let the sequence $(z_i)_{i\in I}$ be a point
in $\Lambda_I^n(Z_\bullet)$. Let $I^\prime\subset [n-1]$ be the
subset determined by $I=\delta_k(I^\prime)$. There is such a subset,
since $k\not\in I$. Moreover, $\abs{I^\prime}=\abs{I}$ and the map
$\delta_k\colon I^\prime \to I$ is a bijection. For $i\in I^\prime$
we define elements
\[
y_i=
\begin{cases}
d_{k-1}(z_{\delta_k(i)})=d_{k-1}(z_i)& \text{ if $i \in I^\prime$, and $i<k$,}\\
d_{k}(z_{\delta_k(i)})=d_{k}(z_{i+1}) & \text{ if $i \in I^\prime$,
and $k \leq i$.}
\end{cases}
\]
We claim that the sequence $(y_i)_{i\in I^\prime}$ forms an
$I^\prime$-cohorn of dimension $(n-1)$ in $Z_{\bullet}$. This can be
seen either by meditating on the geometry of the simplex, or more
computationally in the following way: for $i,j\in I^\prime$ and
$i<j$, we have
\[
d_iy_j=
\begin{cases}
d_id_{k-1}(z_{j})=d_{k-2}d_i(z_j)         &i<j<k\\
d_id_{k}(z_{j+1})=d_{k-1}d_i(z_{j+1})         &i<k \leq j\\
d_{i}d_{k}(z_{j+1})=d_kd_{i+1}(z_{j+1})       &k \leq i<j\\
\end{cases}
\]
and also, using that $(z_i)_{i\in I}$ is in the space of
$I$-cohorns,
\[
d_{j-1}y_i=
\begin{cases}
d_{j-1}d_{k-1}(z_{i})=d_{k-2}d_{j-1}(z_i)=d_{k-2}d_i(z_j)    &i<j\leq k\\
d_{j-1}d_{k-1}(z_{i})=d_{k-1}d_j(z_{i})=d_{k-1}d_i(z_{j+1})    &i<k < j\\
d_{j-1}d_{k}(z_{i+1})=d_kd_{j}(z_{i+1})=d_{i+1}(z_{j+1})    &k \leq i<j.\\
\end{cases}
\]

The association $(z_i)_{i\in I} \mapsto (y_i)_{i\in I^\prime}$
defines a continuous map $P^I_k\colon \Lambda_I^n(Z_\bullet)\to
\Lambda_{I^\prime}^{n-1}(Z_\bullet)$. Since $k\in J$ there is also a
continuous map $p_k\colon \Lambda_{J}^n(Z_\bullet)\to Z_{n-1}$
defined by picking out the appropriate component, i.e. we have
$p_k((z_i)_{i\in J})=z_k$. These maps fit into a diagram of spaces
\[
\xymatrix@C=1.5cm{ \Lambda_J^n(Z_\bullet) \ar[r]^{p_k}\ar[d]^{r_I^J}
& Z_{n-1} \ar[d]^{c_{I^\prime}^{n-1}}\\
\Lambda_I^n(Z_\bullet)\ar[r]^{P^{I}_k}
&\Lambda^{n-1}_{I^\prime}(Z_\bullet) }
\]
This is actually a pullback diagram. The right map is a surjective
Serre fibration by the induction hypothesis of the first induction,
so it follows that $r^J_I$ is also a surjective Serre fibration.
\end{proof}

We can now prove the main theorem of this section:

\begin{Thm}
\label{th:Kan} Let $Z_{\bullet}$ be a simplicial space. Assume that
$Z_{\bullet}$ satisfies the topological Kan condition. Then the
diagonal simplicial set $D_\bullet=\diag(\Sing Z_\bullet)$ is a Kan
complex, i.e. it satisfies the discrete Kan condition.
\end{Thm}

\begin{proof}
We need to show that any horn in $D_\bullet$ can be filled. The Kan
condition is always satisfied for $n=1$, since a horn of dimension
$1$ is filled by the degeneracy map $s_0\colon D_0\to D_1$. So we
assume without restriction that $n\geq 2$. A $k$-horn of dimension
$n$ in $D_\bullet$ consists of a sequence of continuous maps
\[
z_i\colon \Delta^{n-1} \to Z_{n-1} \quad \text{ defined for $0\leq i
\leq n-1$, $i\not=k$,}
\]
with the property that if $z_i$ and $z_j$ are defined and $i<j$,
then the following two continuous maps agree:
\begin{align*}
  \Delta^{n-2} &\xrightarrow{\delta_i} \Delta^{n-1}
\xrightarrow{z_j} Z_{n-1} \xrightarrow{d_i} Z_{n-2},\\
  \Delta^{n-2} &\xrightarrow{\delta_{j-1}} \Delta^{n-1}
\xrightarrow{z_i} Z_{n-1} \xrightarrow{d_{j-1}} Z_{n-2}.\\
\end{align*}
Here, and in the rest of the proof, we abuse notation by not
distinguishing between a map $[m]\to [n]$ and its affine extension
to a continuous map $\Delta^m\to \Delta^n$ between the associated
topological standard simplices. The way to visualize these $z_i$'s,
is to consider the topological horn
\[
\Lambda^n_k=\cup_{i\not=k}\delta_i(\Delta^{n-1})\subset \Delta^{n}.
\]
The domain of $z_i$ can be identified with the subset
$\delta_i(\Delta^{n-1})$ of $\Lambda^n_k$. The maps corresponding to
$z_i$ and $z_j$ do not necessarily agree on the intersection
$\delta_i(\Delta^{n-1})\cap \delta_j(\Delta^{n-1})$ of these
subsets, but satisfy a more complicated compatibility condition.

\emph{Main step.} We define a continuous map
\[
y\colon \Lambda^n_k\to Z_{n}
\]
such that for every $i\not=k$ we have that $z_i=d_i \circ y\circ
\delta_i$. We define this map inductively. The induction is done
over a sequence of subspaces of the skeleton of $\Lambda^n_k$, which
we define now. By $v_i\in\Delta^n$ we denote the $i$'th vertex. Let
$W_r$ be the union of all $r$-dimensional subsimplices of
$\Delta^{n}$ which contain $v_k$. Then
\[
\{v_k\}=W_0\subset W_1\subset \dots W_{n-1}=\Lambda^n_k.
\]

\emph{Start of induction.} $W_0$ is included in every subspace
$\delta_i(\Delta^n)$, $i\neq k$, so there exists a sequence of
vertices $\phi^i\in\Delta^{n-1}$, $i\neq k$, such that
$\delta_i\phi^i=v_k$. Let $(b_i)_{i\neq k}$ be the sequence in
$Z_{n-1}$ given by $b_i=z_i\phi^i$. Take $i<j$, $i\not=k\not=j$,
then there exists a vertex $\psi\in\Delta^{n-2}$ such that
$\delta_j\delta_i(\psi)=v_k$, hence $\phi^i=\delta_{j-1}(\psi)$, and
$\phi^j=\delta_i(\psi)$. Now we can check that $b_i$ and $b_j$
satisfy
\[
d_ib_j=d_iz_j\phi^j=d_iz_j\delta_i(\psi)
=d_{j-1}z_i\delta_{j-1}(\psi)=d_{j-1}z_i\phi^i=d_{j-1}b_i.
\]
Thus, $(b_i)_{i\neq k}$ lies in $\Lambda^n_k(Z_\bullet)$.
By assumption the map $c^n_k$ is surjective, so we can find an element
$a\in Z_{n}$ such that $d_i(a)=b_i$ for all $i\neq k$. We define $y$ on $W_0=\{v_k\}$
by $y(v_k)=a\in Z_{n}$.

\emph{Induction step.} Assume that we have defined $y\colon
W_{r-1}\to Z_{n+1}$ such that for each $i\neq k$ we have a
commutative diagram:
\[
\xymatrix@C=2cm{ W_{r-1}\cap \delta_i(\Delta^{n-1})
\ar@{^{(}->}[r] \ar@{^{(}->}[d]& W_{r-1}\ar[r]^y&Z_{n}\ar[d]^{d_i}\\
\Delta^{n-1}\ar[rr]^{z_i}&&Z_{n-1}.\\
}
\]
We want to extend $y$ to a map from $W_r$ with the same property. To
define this extension, it is sufficient to define it on each
$r$-simplex in $W_r$. Such a simplex is determined by an injective
order preserving map $\tau\co [r]\hookrightarrow [n]$ with the
property that $\tau(k')=k$ for some $k'$ in $[r]$. Let $I\subset
[n]$ be the complement of the image of $\tau$. Let us now translate
the problem of extending the map $y$ to the simplex $\tau\subset
W_r$, into a problem of constructing a map $y^\tau \colon
\Delta_r\to Z_n$ with certain properties. Since $\tau$ is a general
$r$-simplex of $W_r$ the collection of all these maps $y^\tau$
defines the extension of $y$ to $W_r$. We have to be careful about
what the transformed problem is.

First, note that $\tau(\Delta^r)\cap W_{r-1}=\tau(\Lambda^r_{k'})$.
Thus $y^\tau$ is already defined on $\Lambda^r_{k'}$ by induction.

There are two cases for the maps $z_i$.
If $i\not\in I$, then $\tau(\Delta^r)\cap\delta_i(\Delta^{n-1})$
is contained in $W_{r-1}$. So the $z_i$'s corresponding to these $i$'s do not
give extra conditions on $y^\tau$, beyond that it has to be an extension of
the already defined map $y^\tau \colon \Lambda^r_{k'}\rightarrow Z_n$. On the other
hand, if $i\in I$, we do get an extra
condition.

The image $\tau(\Delta^r)$ is included in every subspace
$\delta_i(\Delta^n)$, $i\in I$, so there exists a sequence of
injective order preserving maps $\phi^i\co [r]\hookrightarrow[n-1]$,
$i\in I$, such that $\delta_i\phi^i=\tau$. Let $(b_i)_{i\in I}$ be
the sequence of maps $\Delta^r\rightarrow Z_{n-1}$ given by
$b_i=z_i\phi^i$.

Take $i<j$, $i\not\in I\not\ni j$, then there exists an injective
order preserving map $\psi\co [r]\hookrightarrow[n-2]$ such that
$\delta_j\delta_i\psi=\tau$, hence $\phi^i=\delta_{j-1}\psi$, and
$\phi^j=\delta_i\psi$. As before $b_i$ and $b_j$ satisfy
\[
d_ib_j=d_iz_j\phi^j=d_iz_j\delta_i\psi
=d_{j-1}z_i\delta_{j-1}\psi=d_{j-1}z_i\phi^i=d_{j-1}b_i.
\]
Thus, the collection $(b_i)_{i\in I}$ gives a map $b\co
\Delta^r\rightarrow\Lambda^n_k(Z_\bullet)$. The extension problem we
have to solve is the following:
\[
\xymatrix@C=1.5cm{
\Lambda^r_{k'}\ar[r]^{y^\tau}\ar@{^{(}->}[d]& Z_{n}\ar[d]_{c^n_{I}}\\
\Delta^r \ar[r]^{b}\ar@{-->}[ru]& \Lambda_I^n(Z_\bullet).\\
}
\]
However, this is solvable, since the map $c_I^n$ is a Serre
fibration by Lemma~\ref{le:Serre}. This completes the induction
step.

\emph{Final step.} We have constructed a continuous map $y\colon
\Lambda^n_k\to Z_{n}$ such that $z_i=d_i \circ y\circ \delta_i$ for
every $i\not=k$. The composition of $y$ with a retraction
$\Delta^n\rightarrow \Lambda^n_k$ fills the $k$-horn in $D_\bullet$
that we considered initially.
\end{proof}

\begin{Rem}
A $k$-horn, $0\leq k\leq 2$, in dimension $2$ in $D_\bullet=\diag(\Sing Z_\bullet)$
can be filled if $Z_\bullet$ satisfies the weaker condition that
the three boundary maps $d_i\colon Z_2\to Z_1$ are Serre
fibrations, and that $Z_\bullet$ satisfies the discrete Kan
extension condition in dimension $2$.

To see this use the Kan extension property at
the vertex $v_k$ to obtain a point in $Z_2$. Then use that $d_i\colon
Z_2\to Z_1$ is a Serre fibration to obtain a continuous map $y\colon
\Lambda^k_2\to Z_2$. Finally extend this to the topological $2$-simplex.

In this case, we do not need Lemma~\ref{le:Serre}.
\end{Rem}

We end this section by giving a relativizing Theorem~\ref{th:Kan} in
a special case:

\begin{Cor}\label{Cor:topolKanfibr}
Let $Z_\bullet$ be a simplicial space, and $p_0\co Z_0\rightarrow B$
a map of topological spaces. Then the induced map $\diag(\Sing
Z_\bullet)\rightarrow \Sing X$, extending $p_0$, is a Kan fibration
provided that $p_0$ is a Serre fibration and $Z_\bullet$ satisfies
the topological Kan condition.
\end{Cor}

\begin{proof}
Consider a simplicial $n$-dimensional $k$-horn in $\diag(\Sing
Z_\bullet)$ together with a matching $n$-simplex in $\Sing B$. We
need to find an extension to an $n$-simplex in $\diag(\Sing
Z_\bullet)$. The main step of the proof of Theorem~\ref{th:Kan}
transforms this problem into a lifting problem of spaces, namely
\[
\xymatrix@R=0.3cm{ \Lambda^n_k \ar[r] \ar[dd] & Z_n \ar[dd] \ar[dr]^s & \\
&& Z_0. \ar[dl]^{p_0} \\ \Delta^n \ar[r] \ar@{.>}[ruu] & B & }
\]
Here, the degeneracy map $s\co Z_n\rightarrow Z_0$ is a surjective
Serre fibration by Lemma~\ref{le:Serre}, and $p_0$ is a Serre
fibration by assumption. Hence, the lifting exists.
\end{proof}

\section{Good and sufficiently fibrant $2$-groupoids}\label{sect:GSF}

In this section we will give criteria on a topological $2$-category
$2\C$ to ensure that its geometric nerve $\Delta2\C$ is good, and in
the case $2\C$ is a topological $2$-groupoid, we will give criteria
that ensure that $\Delta2\C$ satisfies the topological Kan
condition.

We will use the notion of a \emph{locally equiconnected} (abbr. LEC)
space $X$, i.e. the diagonal map $X\rightarrow X\times X$ is a
closed cofibration. All CW-complexes are LEC spaces by E. Dyer and
S. Eilenberg's adjunction theorem, see~\cite{DyerEilenberg:72}. On
the other hand, as a corollary of their Theorem~II.7, we get:

\begin{Thm}[(Dyer, Eilenberg)]
Assume that $X$ is a LEC space, and $i\co A\hookrightarrow X$ admits
a retraction $r\co X\rightarrow A$. Then $i$ is a closed
cofibration, and $A$ is a LEC space.
\end{Thm}

\begin{proof}
View $A$ as a subset of $X$. A \emph{haloing function} for $A$ is a
function $q\co X\rightarrow I$ such that $q^{-1}(0)=A$. A \emph{halo
retract} $A$ of $X$ is a retract which admits a haloing function.
By~\cite[Theorem~II.7]{DyerEilenberg:72} the result above follows if
we can find a haloing function for our retract $A$.

By Dyer and Eilenberg's characterization of LEC spaces $X$, their
Theorem~II.1, there exists a haloing function $k\co X\times
X\rightarrow I$ for the diagonal in $X\times X$. Now define $q\co
X\rightarrow I$ by the formula
\[
q(x)=k(x,ir(x)).
\]
Since $q(x)=0$ if and only if $x\in A$, we have constructed our
haloing function for $A$.
\end{proof}

\begin{Cor}\label{Cor:LECgood}
A simplicial space $Z_{\bullet}$, where $Z_q$ is a LEC space for infinitely
many $q$, is a good simplicial space.
\end{Cor}

\begin{proof}
By downward induction on $q$ we assume that $Z_q$ is LEC. Any degeneracy $s_i\co
Z_{q-1}\rightarrow Z_{q}$ has a retract, namely $d_i$. Hence, $s_i$
is a closed cofibration and $Z_{q-1}$ is LEC.
\end{proof}

One basic tool to prove that a space is LEC is:

\begin{Thm}[(Heath, {\cite{Heath:86}})]\label{Thm:Heath}
Assume that $X$, $B$ and $E$ are LEC spaces. Then the pullback of
$X\xrightarrow{f}B\xleftarrow{p} E$ is a LEC space provided that $p$
is a Hurewicz fibration.
\end{Thm}

For nerves of topological categories we now get the following:

\begin{Cor}\label{Cor:GoodNerves}
Let $\C$ be a topological category. If the total space of morphisms
$\C_1$ is LEC and the source map $s\co \C_1\rightarrow \C_0$ is a
Hurewicz fibration, then $N_\bullet \C$ is a good simplicial space.
\end{Cor}

\begin{proof}
The space of objects, $\C_0$, is LEC, since it is a retract of
$\C_1$. We will prove by induction that $N_q\C$ is LEC, and the
result then follows by Corollary~\ref{Cor:LECgood}. Observe that
$N_{q+1}\C$ is given as the pullback of
$N_q\C\xrightarrow{v_q}\C_0\xleftarrow{s}\C_1$, where $v_q$ is the
last vertex map, $N_q\C$ is LEC by induction, and $s$ is a Hurewicz
fibration by assumption. By Heath's theorem $N_{q+1}\C$ is also LEC.
\end{proof}

\begin{Def}\label{def:good}
A topological $2$-category $2\C$ is called \emph{good} if its
geometric nerve $\Delta2\C$ is a good simplicial space.
\end{Def}

Let us now list the criteria on a topological $2$-category $2\C$
that together imply that $2\C$ is good:
\begin{description}
\item[\textbf{G1}]\label{assG:s1} The source map $s_1\co 2\C_1\rightarrow2\C_0$
from $1$-morphisms to objects is a Hurewicz fibration.
\item[\textbf{G2}]\label{assG:t2} The target map $t_2\co 2\C_2\rightarrow2\C_1$
from $2$-morphisms to $1$-morphisms is a Hurewicz fibration.
\item[\textbf{G3}]\label{assG:LEC} The total space of $2$-morphisms, $2\C_2$, is LEC.
\item[\textbf{G4}]\label{assG:inv} All $2$-morphisms are isomorphisms
vertically and the assignment $(-)^{-1}\co 2\C_2\rightarrow2\C_2$
sending a $2$-morphism to its vertical inverse is continuous.
\end{description}

The assumptions~\textbf{G1}--\textbf{G3} exclude pathological
$2$-categories. All examples in Section~\ref{sect:examples} satisfy
assumptions~\textbf{G1}--\textbf{G4}.

\begin{Thm}\label{Thm:good}
A topological $2$-category $2\C$ is good if it satisfies the
assumptions~\textbf{G1}, \textbf{G2}, \textbf{G3}, and~\textbf{G4}.
\end{Thm}

Before proving the theorem, we consider the following lemma which,
under assumption~\textbf{G4}, characterize $n$-simplices of
$\Delta2\C$:

\begin{Lem}
Let $2\C$ be a topological $2$-category which satisfies~\textbf{G4}.
Then the space of $n$-simplexes $(\Delta 2\C)_n$ is homeomorphic to
the space of flags
\[
\begin{pmatrix}
x_0 & f_{01} & \phi_{012} & \phi_{013} & \phi_{014} & \phi_{015} & \cdots & \phi_{01n}\\
& x_1 & f_{12} & \phi_{123} & \phi_{124} & \phi_{125} & \cdots & \phi_{12n} \\
&& x_2 & f_{23} & \phi_{234} & \phi_{235} & \cdots & \phi_{23n} \\
&&& x_3 & f_{34} & \phi_{345} & \cdots & \phi_{34n} \\
&&&& x_4 & f_{45} & \cdots & \phi_{45n} \\
&&&&& x_5 & \cdots & \phi_{56n} \\
&&&&&& \ddots & \vdots \\
&&&&&&& x_n
\end{pmatrix}_,
\]
where $x_0,x_1,\ldots,x_n$ are objects of $2\C$, the $f_{i\,i+1}\co
x_i\rightarrow x_{i+1}$ are $1$-morphisms, and $\phi_{i\,i+1\,j}\co
f_{ij}\Rightarrow f_{i+1\,j}*f_{i\,i+1}$, $i+1<j$, are
$2$-morphisms.
\end{Lem}

\begin{proof}
The $1$-morphisms $f_{ij}$ with $i+1<j$ are implicitly defined as
the source of the $2$-morphism $\phi_{i\,i+1\,j}$.

Given a flag of the form above, the only data of an $n$-simplex in
$\Delta2\C$ missing are the $\phi_{ijk}$ with $0\leq i<j<k\leq n$
and $j\neq i+1$. We claim that these may be reconstructed
continuously from the flag and the coherence condition of the
geometric nerve.

Define the \emph{valence} of the symbol $\phi_{ijk}$, $0<i<j<k<n$,
to be the difference $j-i$. Observe that the flag contains all
$\phi_{ijk}$ of valence $1$. A subtetrahedron, of the $n$-simplex
under consideration, is determined by indices $0\leq i<j<k<l\leq n$.
Rewrite the coherence condition for such a subtetrahedron as the
equation
\[
\phi_{ikl}=\left({f_{kl}}*\phi_{ijk}\right)^{-1}\underline{\alpha}
\left(\phi_{jkl}*{f_{ij}}\right)\phi_{ijl}.
\]
Here $(-)^{-1}$ is the inverse $2$-morphism map, and
$\underline{\alpha}$ is the natural associativity isomorphism.
Observe that the valence of the $2$-morphisms on the right hand side
is strictly less than the valence of $\phi_{ikl}$. Inductively, we
use this equation as a definition of $\phi_{ikl}$. A calculation is
needed to check that $\phi_{ikl}$ does not depend on a particular
choice of $j$. For $i=0$, $k=3$, and $l=4$, inspect the following
diagram:
\[
{\xymatrix@C=2cm{
&& {\scriptstyle f_\smind{04}} \ar@{=>}[d]_{\scriptscriptstyle \phi{\,}_\smind{014}}\\
&& {\scriptstyle f_\smind{14}*f_\smind{01} } \ar@{=>}[d]_{\scriptscriptstyle \phi{\,}_\smind{014}*f_\smind{01}}\\
{\scriptstyle f_\smind{24}*f_\smind{02}}
\ar@{=>}[d]^{\scriptscriptstyle \phi{\,}_\smind{234}*f_\smind{02}}
\ar@{=>}[r]^-{\scriptscriptstyle f_\smind{24}*\phi{\,}_\smind{012}}
    & {\scriptstyle f_\smind{24}\!*(f_\smind{12}*f_\smind{01}) } \ar@{=>}[d]_{\scriptscriptstyle \phi{\,}_\smind{234}*(f_\smind{12}*f_\smind{01})} \ar@{=>}[r]
    & {\scriptstyle (f_\smind{24}*f_\smind{12})*f_\smind{01} } \ar@{=>}[d]_{\scriptscriptstyle (\phi{\,}_\smind{234}*f_\smind{12})*f_\smind{01}}\\
{\scriptstyle (f_\smind{34}*f_\smind{23})*f_\smind{02}}
\ar@{=>}[r]^-{\scriptscriptstyle
(f_\smind{34}*f_\smind{23})*\phi{\,}_\smind{012}}
    & {\scriptstyle (f_\smind{34}*f_\smind{23})*(f_\smind{12}*f_\smind{01}) } \ar@{=>}[r]
    & {\scriptstyle ((f_\smind{34}*f_\smind{23})*f_\smind{12})*f_\smind{01} } \\
{\scriptstyle f_\smind{34}*(f_\smind{23}*f_\smind{02}) } \ar@{=>}[u]
\ar@{=>}[r]^-{\scriptscriptstyle
f_\smind{34}*(f_\smind{23}*\phi{\,}_\smind{012})}
    & {\scriptstyle f_\smind{34}*(f_\smind{23}*(f_\smind{12}*f_\smind{01})) } \ar@{=>}[u] \ar@{=>}[d]\\
& {\scriptstyle
f_\smind{34}*((f_\smind{23}*f_\smind{12})*f_\smind{01}) }
\ar@{=>}[r]
    & {\scriptstyle (f_\smind{34}*(f_\smind{23}*f_\smind{12}))*f_\smind{01} } \ar@{=>}[uu]\\
{\scriptstyle f_\smind{34}*f_\smind{03} }
\ar@{=>}[r]^-{\scriptscriptstyle f_\smind{34}*\phi{\,}_\smind{013}}
    & {\scriptstyle f_\smind{34}*(f_\smind{13}*f_\smind{01}) } \ar@{=>}[r] \ar@{=>}[u]^{\scriptscriptstyle f_\smind{34}*(\phi{\,}_\smind{123}*f_\smind{01})}
    & {\scriptstyle (f_\smind{34}*f_\smind{13})*f_\smind{01} }. \ar@{=>}[u]^{\scriptscriptstyle (f_\smind{34}*\phi{\,}_\smind{123})*f_\smind{01}}
}}
\]
The unmarked arrows are natural associativity isomorphisms. Going
from $f_{04}$ to $f_{34}*f_{03}$ on the left side corresponds to
defining $\phi_{034}$ with $j=2$, while the composition on the right
side corresponds to $j=1$. Since the diagram commutes, the two
possible definition agree.
\end{proof}

\begin{proof}[of theorem~\ref{Thm:good}]
It is enough to show that for all $n$ the space of $n$-simplices,
$(\Delta2\C)_n$ is a LEC space.

By assumption~\textbf{G3} it follows that $2\C_2$ is LEC. Observe
that both $2\C_0$ and $2\C_1$ are retracts of $2\C_2$, where the
inclusions are provided by identity maps and the retractions by
source maps. Hence, $2\C_0$ and $2\C_1$ are also LEC spaces.

This provides the start of an induction on the simplicial degree
$n$. Assume that $(\Delta2\C)_{n-1}$ is LEC. We can construct
$(\Delta2\C)_n$ from the $(n-1)$-simplices by iterated pullback,
i.e. we have a tower
\[
(\Delta2\C)_n=P_n\twoheadrightarrow\cdots \twoheadrightarrow
P_1\twoheadrightarrow P_0=(\Delta2\C)_{n-1}.
\]
The steps of the tower corresponds to the entries in the last column
of the flag describing $n$-simplices. $P_1$ is the pullback of
\[
P_0\xrightarrow{x_{n-1}}2\C_0\overset{s_1}{\twoheadleftarrow} 2\C_1.
\]
Here we adjoin the $1$-morphism $f_{n-1\,n}$ along $x_{n-1}$. The
map $s_1$ is a Hurewicz fibration by assumption~\textbf{G1}. Hence
$P_1$ is LEC by Theorem~\ref{Thm:Heath}. The space $P_k$, $k>1$, is
the pullback of
\[
P_{k-1}\xrightarrow{f_{n-k+1,\,n}*f_{n-k,\,n-1+1}}2\C_1\overset{t_2}{\twoheadleftarrow}2\C_2.
\]
Here we adjoin the $2$-morphisms $\phi_{n-k,\,n-k+1,\,n}$ along
$f_{n-k+1,\,n}*f_{n-k,\,n-1+1}$. By assumption~\textbf{G2}, $t_2$ is
a Hurewicz fibration. So $P_k$ is LEC by Heath's theorem.
\end{proof}

Let $2\C$ be a strict discrete $2$-category. For each pair of
objects $x,y$ in $2\C_0$ we have a morphism category $2\C(x,y)$.
Applying the nerve to each of these morphism categories, we get a
simplicial category $N_\bullet2\C$. Apply the nerve construction
again to produce a bisimplicial set $N_\bullet N_\bullet 2\C$.
Bullejos and Cegarra~\cite{BullejosCegarra:03} prove that the
geometric nerve $\Delta 2\C$ is naturally weakly equivalent to the
double nerve $N_\bullet N_\bullet 2\C$.

Inspired by this, we now want to study the geometric nerve of
a strict topological $2$-category $2\C$ via
the nerves $N_\bullet 2\C(x,y)$ of its morphism categories. To do this
we introduce $\operatorname{Ar}^{\parallel}(2\C)$, the subspace of
$2\C_1\times2\C_1$ consisting of a pair $(f,g)$ of parallel
$1$-morphisms, i.e. the source of $f$ is also the source of $g$ and
the target of $f$ is also the target of $g$. Observe that we have a
pullback diagram
\[
\xymatrix{ \operatorname{Ar}^{\parallel}(2\C) \ar[r]^f \ar[d]^g &
2\C_1 \ar[d]^{(s_1,t_1)} \\ 2\C_1 \ar[r]^{(s_1,t_1)} & 2\C_0\times
2\C_0.}
\]
Next we strengthen the assumptions~\textbf{G1} and~\textbf{G2} as follows:

\begin{description}
\item[G1'] The source-target map of $1$-morphisms,
$2\C_1\xrightarrow{(s_1,t_1)}2\C_0\times2\C_0$, is a Hurewicz
fibration.
\item[G2']\label{ass:2st} The source-target map of
$2$-morphisms,
$2\C_2\xrightarrow{(s_2,t_2)}\operatorname{Ar}^{\parallel}(2\C)$, is
a Hurewicz fibration.
\end{description}

Together with~\textbf{G3} and~\textbf{G4}, these assumptions
simplifies comparison of two geometric nerves of strict $2$-categories:

\begin{Lem}\label{lem:DeltaSimplified}
Let $2\C$ and $2\D$ be strict topological $2$-categories that
satisty the assumptions~\textbf{G1'}, \textbf{G2'}, \textbf{G3}
and~\textbf{G4}. Furthemore, let $F\co 2\C\rightarrow 2\D$ be a
continuous strict $2$-functor. If the map on objects $F_0\co
2\C_0\rightarrow 2\D_0$ is a weak equivalence, and for each pair of
objects $x,y\in2\C_0$ the induced map between the nerves of morphism
categories $N_\bullet2\C(x,y)\rightarrow N_\bullet2\D(F(x),F(y))$ is
a weak equivalence, then $F$ induces a weak equivalence of the
geometric nerves, $\Delta2\C\xrightarrow{\simeq}\Delta2\D$.
\end{Lem}

\begin{proof}
We will prove the lemma via a series of claims:

\emph{First claim.} There is a natural weak equivalence
$|\Delta2\C|\simeq |\Sing N_\bullet N_\bullet 2\C|$.

Since $2\C$ is strict, we may form the double nerve $N_\bullet
N_\bullet 2\C$, as well as the geometric nerve $\Delta2\C$. Because
$\Delta2\C$ is a good simplicial space, by Theorem~\ref{Thm:good},
we have a weak equivalence $|\Sing\Delta 2\C|\simeq |\Delta 2\C|$.
Notice that the functor $\Sing$ commutes with nerve constructions.
Fix a simplicial degree $k$ for the $\Sing$-direction. Observe that
$\operatorname{Sing}_k2\C$ is a strict discrete $2$-category. By
Bullejos and Cegarra's theorem we get a natural weak equivalence
\[
|\Delta\operatorname{Sing}_k2\C|\simeq|N_\bullet N_\bullet
\operatorname{Sing}_k2\C|.
\]
This proves the claim since a map between simplicial spaces, which
in each simplicial degree is a weak equivalence, is itself a weak
equivalence.

\emph{Second claim.} The simplicial space $N_\bullet 2\C$ is good,
and moreover, for each pair of objects $x$, $y$ in $2\C$ the
simplicial space $N_\bullet 2\C(x,y)$ is good.

The simplicial space $N_\bullet 2\C$ is the nerve construction
applied to the vertical composition of $2$-morphisms. Thus $N_k2\C$
is topologized as a subspace of the $k$-fold product
$2\C_2\times\cdots\times2\C_2$. Taking the source and target object,
we get a map $N_\bullet 2\C\rightarrow 2\C_0\times 2\C_0$, and
$N_\bullet 2\C(x,y)$ denotes the fiber over the pair of objects
$x,y$.

By Corollary~\ref{Cor:GoodNerves} it is sufficient to check that
$2\C_2$ is LEC and the source map $s_2\co 2\C_2\rightarrow 2\C_1$ is
a Hurewicz fibration. The former is assumption~\textbf{G3}, while
the latter holds because of~\textbf{G1'} and~\textbf{G2'}.

Since the source-target map $2\C_1\rightarrow 2\C_0\times 2\C_0$ is
a Hurewicz fibration,~\textbf{G1'}, it follows by Heath's theorem
that the source map $2\C_2(x,y)\rightarrow 2\C_1(x,y)$ is a Hurewicz
fibration between LEC spaces, for every pair of objects $x,y$ in
$2\C_0$. Thus $N_\bullet2\C(x,y)$ also is good.

\emph{Third claim.} The map $|\Sing N_\bullet 2\C|\rightarrow
|\Sing(2\C_0\times 2\C_0)|$ is a Serre fibration with fiber $|\Sing
N_\bullet2\C(x,y)|$ over the path component corresponding to $(x,y)$
in $2\C_0\times2\C_0$.

Given that the map is a Serre fibration, the statement about the
fiber is obvious. Since the geometric realization of a Kan fibration
is a Serre fibration, we can use Corollary~\ref{Cor:topolKanfibr}.
Hence, we need to show that the map $p_0\co N_0 2\C\rightarrow
2\C_0\times 2\C_0$ is Serre fibration, and that $N_\bullet2\C$
satisfies the topological Kan condition. Now observe that $p_0$ is
the source-target map
$2\C_1\xrightarrow{(s_1,t_1)}2\C_0\times2\C_0$, whence a Serre
fibration by~\textbf{G1'}.

The maps $c^1_k\co N_12\C\rightarrow\Lambda^1_k(N_\bullet 2\C)$ are
for $k=0$ and $k=1$ the target and source maps $t_2,s_2\co
2\C_2\rightarrow 2\C_1$. These maps are Serre fibrations
by~\textbf{G1'} and~\textbf{G2'}, and they are surjective by the
vertical identity map $2\C_1\rightarrow2\C_2$.

The maps $c^2_k\co N_22\C\rightarrow\Lambda^2_k(N_\bullet 2\C)$ are
homeomorphisms by~\textbf{G4}, while the maps $c^n_k$, for $n\geq3$,
are homeomorphisms because the nerve is $2$-coskeletal.

\emph{Fourth claim.} For each $k$ the map $|\Sing N_k N_\bullet
2\C|\rightarrow |\Sing 2\C_0^{\times(k+1)}|$ is a Serre fibration
with fiber $|\Sing N_\bullet2\C(x_0,x_1)|\times\cdots\times|\Sing
N_\bullet2\C(x_{k-1},x_{k})|$ over the path component corresponding
to $(x_0,x_1,\ldots,x_k)$ in $2\C_0^{\times (k+1)}$.

Under consideration here is a category which has $2\C_0$ as its
space of objects and where the collection of $1$-morphisms forms a
simplicial space $N_\bullet 2\C$. We denote this category by its
morphisms, namely $N_\bullet 2\C$. Observe that the horizontal
composition in $2\C$ induces the composition in $N_\bullet2\C$.

Since both $\Sing$ and $|-|$ commutes with finite limits, the
following diagram is pullback
\[
\xymatrix{ |\Sing N_{k+1} N_\bullet 2\C| \ar[r] \ar@{->>}[d] &
|\Sing N_k N_\bullet 2\C|\times |\Sing N_\bullet 2\C| \ar@{->>}[d]
\\
|\Sing 2\C_0^{\times(k+2)}|\ar[r] & |\Sing
2\C_0^{\times(k+1)}|\times |\Sing(2\C_0\times2\C_0)|.}
\]
The claim now follows by induction on $k$.

\emph{Final claim.} The map $|\Sing N_\bullet N_\bullet
2\C|\rightarrow |\Sing N_\bullet N_\bullet 2\D|$ is a weak
equivalence.

Observe that the previous claims applies to $2\D$ as well. Fix a
simplicial degree $k$ for the nerve corresponding to horizontal
composition. Let $x_0,x_1,\ldots,x_k$ be any choice of objects in
$2\C$, and set $y_i=F(x_i)\in2\D_0$. We get a map of fiber sequences
\[
\xymatrix{ |\Sing N_\bullet\left( 2\C(x_0,x_1)\times\cdots\times
2\C(x_{k-1},x_{k})\right) | \ar[r] \ar[d] &
|\Sing N_\bullet\left( 2\D(y_0,y_1)\times\cdots \times 2\D(y_{k-1},y_{k})\right)| \ar[d] \\
|\Sing N_k N_\bullet 2\C| \ar[d] \ar[r] & |\Sing N_k N_\bullet 2\D|
\ar[d] \\ 2\C_0^{\times(k+1)} \ar[r] & 2\D_0^{\times(k+1)}.}
\]
By the assumptions of the lemma, and the second claim, the top and
bottom maps are weak equivalences, and it follows that the middle
map is a weak equivalence. Since this holds for all $k$, the final
claim holds.

Together the first and the final claim proves the lemma.
\end{proof}

We will now give criteria on $2\C$ to ensure that $\Delta2\C$
satisfies the topological Kan condition. Recall from
\cite{Duskin:02} that a $2$-groupoid is a $2$-category whose
$2$-morphisms are isomorphisms, and whose $1$-morphisms
$x_0\xrightarrow{f} x_1$ induce equivalences of categories,
internally in $2\C$:
\begin{align*}
-*f\colon
2\C(x_1,y) \to 2\C(x_0,y),\\
f*-\colon 2\C(y,x_0)\to 2\C(y,x_1).
\end{align*}
We need to make the equivalence of categories continuous. We first
remark that an equivalence of categories is the same as a pair of
adjoint equivalences. For this and for detail about adjoint
functors, see~\cite[Chapter~IV]{MacLane:98}.

Our approach is to define a space $\AdEq(2\C)$ whose points are
tuples $(f,g,\eta,\epsilon)$ consisting of $1$-morphisms $f\co
x_0\rightarrow x_1$ and $g\co x_1\rightarrow x_0$ together with
$2$-morphisms $\eta\co \id_{x_0}\Rightarrow g*f$ and $\epsilon\co
f*g\Rightarrow \id_{x_1}$ such that $\eta$ and $\epsilon$ are
isomorphisms which satisfy
\[
(\epsilon*f)(f*\eta)=\id_f\quad\text{and}\quad(g*\epsilon)(\eta*g)=\id_g.
\]
We give $\AdEq(2\C)$ the subspace topology with respect to the
inclusion into $2\C_1\times2\C_1\times2\C_2\times2\C_2$. Thus the
map $\pi\co \AdEq(2\C)\rightarrow 2\C_1$ sending
$(f,g,\eta,\epsilon)$ to $f$ is continuous.

We arrive at the following reasonable set of fibrancy conditions:
\begin{description}
\item[SF1]\label{assSF:groupoid} $2\C$ is a topological $2$-groupoid.
\item[SF2] The source and target maps for the $1$-morphisms are
Serre fibrations. \label{cond:Serre1}
\item[SF3] The maps from $2\C_2$ to $2\C_1$ given by the source
and by the target of $2$-morphisms are both Serre fibrations.
\label{cond:Segal}
\item[SF4] The assignment sending a $2$-morphism to its vertical inverse is
continuous. \label{cond:Serre2}
\item[SF5] The map $\pi\co \AdEq(2\C)\rightarrow 2\C_1$ is a Serre
fibration. \label{cond:adeq}
\end{description}

\begin{Lem}\label{le:equation}
Assume that $2\C$ is a topological $2$-category satisfying
\textbf{SF1}, \textbf{SF4}, and \textbf{SF5}. Let $\psi\co
D^n\rightarrow 2\C_2$ and $f,h_s,h_t\co D^n\rightarrow 2\C_1$ be
maps such that the target of $f$ is the source of both $h_s$ and
$h_t$, and $\psi\co h_s*f\Rightarrow h_t*f$. Then there exists a
unique map $\phi\co D^n\rightarrow2\C_2$ such that
\[
\phi*f=\psi
\]
for all points in $D^n$. Similarly, if $\psi\co D^n\rightarrow2\C_2$
and $f',h'_s,h'_t\co D^n\rightarrow 2\C_1$ are such that the source
of $f$ is equal to the target of both $h'_s$ and $h'_t$, and
$\psi\co f'*h'_s\Rightarrow f'*h'_t$, then there is a unique
$\phi'\co D^n\rightarrow2\C_2$ such that $f'*\phi'=\psi$.
\end{Lem}

\begin{proof}
We prove the first part, the last statement is dual. We view
$\phi*f=\psi$ as an equation parameterized over $D^n$, where $f$ and
$\psi$ are known, and $\phi$ is the unknown to be solved for.

By the assumptions \textbf{SF1} and \textbf{SF5} there are maps
$g\co D^n\rightarrow2\C_1$ and $\eta,\epsilon\co
D^n\rightarrow2\C_2$ such that $(f,g,\eta,\epsilon)$ is an adjoint
equivalence for all points in $D^n$. The following commutative
diagram shows that $\phi$ can be expressed uniquely in terms of
$\psi*g=(\phi*f)*g$, the adjoint equivalence, and other known
quantities:
\[
\xymatrix@C=1.5cm{ (h_s*f)*g \ar@{=>}[d]_-{\psi*g} & h_s*(f*g)
\ar@{=>}[d]^{\phi*(f*g)} \ar@{=>}[l]_-{\underline{\alpha}}
\ar@{=>}[r]^-{h_s*\epsilon} & h_s*\id \ar@{=>}[d]^{\phi*\id}
\ar@{=>}[r]^-{\underline{\rho}}
& h_s \ar@{=>}[d]^\phi\\
(h_t*f)*g & h_t*(f*g) \ar@{=>}[l]_-{\underline{\alpha}}
\ar@{=>}[r]^-{h_t*\epsilon} & h_t*\id
\ar@{=>}[r]^-{\underline{\rho}} & h_t. }
\]
Here $\underline{\alpha}$ and $\underline{\rho}$ denote the natural
associativity isomorphism and natural right unit isomorphism of
$2\C$ respectively. Also, when writing out the expression for
$\phi$, we observe that assumption \textbf{SF4} is used to invert
some $2$-morphisms. This completes the proof.
\end{proof}

We made the definition above because of the following theorem.

\begin{Thm}\label{Thm:SF}
If \textbf{SF1}, \textbf{SF2}, \textbf{SF3}, \textbf{SF4},
and~\textbf{SF5} holds for a topological $2$-category $2\C$, then
the geometric nerve $\Delta2\C$ satisfies the topological Kan
condition.
\end{Thm}

\begin{proof}
We have to show that all fill-horn-maps
\[
c_k\co (\Delta2\C)_n\rightarrow \Lambda^k_n(\Delta2\C)
\]
are surjective Serre fibrations. By~\cite[Theorem~8.6]{Duskin:02}
all maps $c_k$ are surjective. It remains to show that they are
Serre fibrations. Since $\Delta2\C$ is $3$-coskeletal, the $c_k$s
are homeomorphisms in dimensions $n>3$. For $n=1$ the fill-horn-maps
$c_0$ and $c_1$ are the source and target maps
$2\C_1\rightarrow2\C_0$, which are Serre fibrations by
condition~\textbf{SF2}.

For $n=3$ there are essentially two different cases to consider. For
$k=1$ and $k=2$ the fill-horn-map $c_k$ asks us to solve the
coherence equation
\[
(\phi_{123}*f_{01})\phi_{013}=\underline{\alpha}(f_{23}*\phi_{012})\phi_{023}
\]
with respect to $\phi_{023}$ and $\phi_{013}$ given all other $1$-
and $2$-morphisms. This can be done continuously using assumption
\textbf{SF4} alone. Hence $c_1$ and $c_2$ are homeomorphisms. Next
consider the case $k=0$. In order to show that $c_0$ is a Serre
fibration we view the coherence condition above as an equation,
parameterized over some disk $D^m\times I$, and with $\phi_{123}$ as
the unknown. We want to extend a solution $\phi_{123}\co
D^m\times\{0\}\rightarrow2\C_2$ to all parameters $D^m\times I$.
Using the continuous inverse we reduce to a parameterized equation
of the form $\phi*f=\psi$. This has a unique solution by
Lemma~\ref{le:equation}. The case $k=3$ is similar.

The rest of the proof concerns the case $n=2$, i.e. we are trying to
find lifts for all diagrams of the form
\[
\xymatrix{ D^m\times\{0\} \ar[r] \ar@{^(->}[d] & \Delta(2\C)_2 \ar[d]^{c_k} \\
D^m\times I \ar[r] \ar@{.>}[ur] & \Lambda^k_2(\Delta2\C).}
\]
In the case $k=1$ observe that the following diagram is pullback:
\[
\xymatrix@C=2cm{ \Delta(2\C)_2 \ar[d]^{c_1} \ar[r]^{\phi_{012}} & 2\C_2 \ar[d]^{\text{target}} \\
\Lambda^1_2(\Delta2\C) \ar[r]^{f_{12}*f_{01}} & 2\C_1. }
\]
Thus $c_1$ is a Serre fibration by assumption \textbf{SF3}.

Next we consider the case $k=0$. The map $D^m\times I\rightarrow
\Lambda^0_2(\Delta2\C)$ corresponds to families of $1$-morphisms
$f_{01}$ and $f_{02}$, parametrized over $D^m\times I$, such that
the source of $f_{01}$ equals the source of $f_{02}$. The map
$D^m\rightarrow\Delta2\C_2$ corresponds to two additional maps
$f_{12}\co D^m\rightarrow2\C_1$ and $\phi_{012}\co D^m\rightarrow
2\C_2$ such that the following diagram commutes for all parameters
in $D^m$:
\[
\xymatrix{
&\bullet\ar[rd]^{f_{12}} &\\
\bullet\ar[ru]^{f_{01}}\ar[rr]_{f_{02}}&\ar@{=>}[u]_(.3){\phi_{012}}&\bullet.\\
}
\]
We can choose $g\co D^m\times I\rightarrow2\C_1$ and
$\eta,\epsilon\co D^m\times I\rightarrow2\C_2$ such that
$(f_{01},g,\eta,\epsilon)$ is an adjoint equivalence for all
parameters. By assumption \textbf{SF3} we may solve the lifting
problem
\[
\xymatrix@C=3.5cm{ D^m\times\{0\}
\ar[r]^{(f_{12}*\epsilon)(\phi_{012}*g)}
\ar@{^(->}[d] & 2\C_2 \ar[d]^{\text{source}} \\
D^m\times I \ar[r]^{f_{02}*g} \ar@{.>}[ur]^{\psi} & 2\C_1.}
\]
Furthermore, by Lemma~\ref{le:equation} and assumption~\textbf{SF4}
the equation
\[
(f_{12}*\epsilon)(\phi_{012}*g)=\psi
\]
has a unique solution for $\phi_{012}$ for all parameters $D^m\times
I$, and by uniqueness this solution extends the given $\phi_{012}$
defined for parameters $D^m\times\{0\}$.

The case $k=2$ is dual to the case $k=0$. This finishes the proof.
\end{proof}

\section{Concordance theory}\label{sect:concordance}

Fix a simplicial space $Z_{\bullet}$. In this section we will define
the notion of a $Z_{\bullet}$-bundle over a topological space $X$,
and study its homotopy properties up to concordance. In the case
$Z_{\bullet}=\Delta 2\C$ this recovers the notion of principal
$2\C$-bundles.

\begin{Def}
A \emph{$Z_{\bullet}$-bundle} $\E$ over $X$ consists of an ordered
open cover $\U$ together with a simplicial map $\phi\co
U_{\bullet}\rightarrow Z_{\bullet}$, where $U_{\bullet}$ is the
ordered \v{C}ech complex of $\U$.
\end{Def}

We have to be precise about our conventions and definitions
regarding covers. For every cover $\U$ we will assume that each
$U_\alpha\in \U$ is non-empty. We say that $\U'$ is a
\emph{refinement} of $\U$ if there exists a \emph{carrier function}
$c\co \I'\rightarrow \I$ between the respective indexing sets such
that for each $\alpha\in\I'$ the set $U'_\alpha$ is contained in
$U_{c(\alpha)}$. Whenever the indexing sets $\I$ and $\I'$ are
partially ordered, we will demand that any carrier function is
order-preserving. Observe that to any refinement together with a
choice of carrier function there is a canonical way of associating a
simplicial map $U'_\bullet\rightarrow U_\bullet$ between the
corresponding (ordered) \v{C}ech complexes. Specializing the notion
of a refinement with carrier function, we say that a cover $\U'$ is
a \emph{shrinking} of $\U$, and we write $\U'\subseteq\U$, if the
indexing set $\I'$ is a subset of $\I$ and for each $\alpha\in\I'$
we have $U'_\alpha\subseteq U_\alpha$. A cover $\U$ is \emph{good}
if each finite intersection
$U_{\alpha_0\cdots\alpha_k}=\bigcap_{i=0}^k U_{\alpha_i}$ is either
empty or contractible. Given a map $f\co X\rightarrow Y$, we may
pull back a $Z_\bullet$-bundle $\E$ over $Y$ to $X$ as follows: let
$\U'$ be the ordered open cover of $X$ consisting of the non-empty
open sets having the form $f^{-1}(U_\alpha)$ where $U_\alpha\in \U$.
There is a natural map of ordered \v{C}ech complexes
$U'_\bullet\rightarrow U_\bullet$. Define $f^*\E$ to be $\U'$
together with the composition $U'_\bullet\rightarrow
U_\bullet\xrightarrow{\phi} Z_\bullet$. The pullback to a subspace
$A\subseteq X$ is called a \emph{restriction} and will be denoted
$\E|_A$.

\begin{Def}
We say that $\E_0$ and $\E_1$ over $X$ are \emph{concordant} if
there exists a $Z_{\bullet}$-bundle $\E$ over $X\times I$ such that
the restrictions to $X\times\{0\}$ and $X\times\{1\}$ give $\E_0$
and $\E_1$ respectively. Let $\Con(X,Z_{\bullet})$ be the set of all
concordance classes of $Z_{\bullet}$-bundles over $X$. The class in
$\Con(X,Z_{\bullet})$ corresponding to $\E$ will be written as
$\left[\E\right]$.
\end{Def}

Observe that concordance is an equivalence relation.

Before proceeding with building the concordance theory of
$Z_\bullet$-bundles, it is nice to present an interesting example of
$Z_\bullet$-bundles where $Z_\bullet$ is not the geometric nerve of
a $2$-category:

\begin{Exa}
Wirth and Stasheff,~\cite{WirthStasheff:06}, describe locally
homotopy trivial fibrations with fiber $F$. It seems plausible that
we can define a simplicial space with $k$-simplexes, $Z_k$, the
space of all homotopy coherent functors from $[k]$ to the
topological monoid of homotopy equivalences of $F$, $H(F)$. In that
case a $Z_{\bullet}$-bundle should be the same as a homotopy
transition cocycle, see~\cite[Definition~2.5]{WirthStasheff:06}.
Moreover, if $Z_{\bullet}$ is good, then our
Theorem~\ref{thm:generalizationofmain} shows that these fibrations
are classified by $|Z_{\bullet}|$. Furthermore, $|Z_{\bullet}|$
should be weakly equivalent to $BH(F)$ for reasonable $F$.
\end{Exa}

In the theory that follows we will also consider classes of
$Z_\bullet$-bundles over $X$ fixed on some subspace $A$. We say that
two $Z_\bullet$-bundles $\E_0$ and $\E_1$ over $X$ are \emph{equal
on $A$} if the restriction $\E|_A$ equals $\E'|_A$. We will use the
notation $X\times I/A$ for the quotient space of $X\times I$ where
we identify $(x,t)$ and $(x',t')$ whenever $x=x'$ lies in $A$. Let
$i_0$ and $i_1$ denote the inclusions of $X$ into $X\times I/A$
given by sending $x\in X$ to $(x,0)$ and $(x,1)$ respectively. If
$\E_0$ and $\E_1$ over $A$ are equal on $A$, then they are
\emph{concordant relative to $A$} if there exists a
$Z_\bullet$-bundle $\E'$ over $X\times I/A$ such that
$i_0^*\E'=\E_0$ and $i_1^*\E'=\E_1$. Fixing a $Z_\bullet$-bundle
$\E_A$ over $A$, we denote by
\[
\Con_{Z_{\bullet}}(X,A;\E_A)
\]
the set of all $\E$ over $X$ whose restriction to $A$ equals $\E_A$
modulo the equivalence relation of concordance relative to $A$.

Taking a more categorical approach to concordance, we will now
consider \emph{presheaves on spaces}, i.e. functors $\F\co
\Top^{\op}\rightarrow \Set$. Sending a space $X$ to the set of
$Z_\bullet$-bundles over $X$ defines our first example of a presheaf
on spaces, and we let $\F_{Z_\bullet}$ denote this presheaf. More
generally, we will consider presheaves, $\F\co
\T^{\op}\rightarrow\Set$, defined on categories $\T$ satisfing the
following conditiuons:
\begin{itemize}
\item[i)] There is a faithful functor $V\co \T\rightarrow \Top$ such
that
\item[ii)] The functor $V$ creates pushouts in $\T$, i.e. whenever
given a span $B\leftarrow A\rightarrow X$ in $\T$, the pushout $Y'$
of
\[
V(B)\leftarrow V(A)\rightarrow V(X)
\]
in $\Top$ has a unique lifting to a pushout $Y$ of the original span
in $\T$.
\item[iii)] Moreover, if $A\hookrightarrow X$ in $\T$ maps to an
inclusion $V(A)\hookrightarrow V(X)$ in $\Top$, then the diagram
\[
\xymatrix{ V(A) \ar@{^{(}->}[r] \ar@{^{(}->}[d] & V(X) \ar[d]^{i_0}\\
V(X) \ar[r]^-{i_1} & V(X)\times I/V(A)}
\]
lifts to $\T$, i.e. the space $X\times I/A$ represents a
well-defined object in $\T$ and the inclusions $i_0,i_1\co
X\rightarrow X\times I/A$ are also well-defined in $\T$.
\end{itemize}
Since $V$ is a faithful functor, we will usually omit it from our
notation. Thus we think about $\T$ as a subcategory of $\Top$. For
later use, here are three examples of categories $\T$ satisfying
conditions i)--iii):

\begin{Exa}\label{Exa:Ts}
\begin{itemize}
\item[i)] Fix a space $A$. Let $\T$ be the
subcategory of $\Top$ consisting of spaces $X$ containing $A$ and
maps $X\rightarrow X'$ restricting to the identity on $A$. We call this
$\T$ the category of \emph{spaces containing $A$}.
\item[ii)] Fix a space $A$. Let $\T$ be the category having as objects spaces $X$ containing $A$ together
with the structure of a relative CW-complex on $(X,A)$. The
morphisms of $\T$ are defined to the inclusion of subcomplexes, i.e.
maps $X\hookrightarrow X'$ where the image of $X$ is a union of $A$
and some of the cells of $(X',A)$. We call this $\T$ the category of
\emph{CW-complexes relative to $A$}.
\item[iii)] Let $\T$ be the full subcategory of \emph{compact
spaces}. Condition~ii) holds since any quotient of a compact space
is compact.
\end{itemize}
\end{Exa}

If $\F$ is a presheaf on $\T$, we define concordance as follows: let
$A\hookrightarrow X$ in $\T$ be an inclusion in $\Top$, and fix some
$s_A\in\F(A)$. Let $\F[X,A;s_A]$ denote the set of all $s$ in
$\F(X)$ that restrict to $s_A$ on $A$ modulo the relation that
$s_0\sim s_1$ if there exists an element $s'\in\F(X\times I/A)$ with
$i_0^*s'=s_0$ and $i_1^*s'=s_1$. Let $[s]$ denote the class of $s$
in $\F[X,A;s_A]$.

Suppose that $\F$ is a presheaf on spaces and $\T$ is a category
satisfying the conditions above. We may restrict $\F$ to a presheaf
on $\T$. Observe that condition iii) ensures that concordance of
$\F|_\T$ equals concordance on $\F$, and we may write
\[
\F[X,A;s_A]=\F|_\T[X,A;s_A]
\]
for all inclusions $A\hookrightarrow X$ and $s_A\in\F(A)=\F|_\T(A)$.

\begin{Def}
We say that a presheaf $\F$ on $\T$ is \emph{excisive} if for all
pushout diagrams
\[
\xymatrix{ A \ar@{^{(}->}[r]^{i} \ar[d]_f & X \ar[d]^g\\ B
\ar@{^{(}->}[r]^{j} & Y}
\]
in $\T$ where $i$ is a closed cofibration on $\Top$, and elements
$s_B\in\F(B)$ and $s_X\in\F(X)$ that restrict to the same element
$s_A=f^*s_B=i^*s_X$ in $\F(A)$, there exists a unique element $s_Y$
of $\F(Y)$ such that $g^*s_Y=s_X$ and $j^*s_Y=s_B$. We say that $\F$
is \emph{strongly excisive} if such an element $s_Y$ exists uniquely
whenever $i$ is an inclusion. We call $\F$ \emph{weakly excisive} if
$s_Y$ exists uniquely whenever $(X,A)$ is a relative CW-complex.
\end{Def}

Observe that condition ii) on $\T$ ensures that restriction
of an excisive presheaf $\F$ on spaces to a presheaf $\F|_\T$ on $\T$
also is excisive.

A map between presheaves on $\T$ is a natural transformation $\nu\co
\F\rightarrow \F'$. Clearly $\nu$ induces a maps of relative
concordance classes
$\nu_{[X,A;s_A]}\co\F[X,A;s_A]\rightarrow\F'[X,A;\nu(s_A)]$. A
useful insight in the theory of concordance is that surjectivity
often implies injectivity,
compare~\cite[Proposition~2.18]{MadsenWeiss:07}. We prove a similar
result:

\begin{Prop}\label{Prop:surjimplyinj}
Let $\F$ and $\F'$ be excisive presheaves on $\T$. If $\nu$ induces
a surjective map
\[
\nu_{[X,C;s_C]}\co \F[X,C;s_C]\rightarrow\F'[X,C;\nu(s_C)].
\]
for all closed cofibrations $C\hookrightarrow X$ and $s_C$ in
$\F(C)$, then all these $\nu_{[X,C;s_C]}$ are bijections. If $\F$
and $\F'$ are strongly excisive and $\nu_{[X,C;s_C]}$ are surjective
for all inclusions $C\hookrightarrow X$, then all $\nu_{[X,C;s_C]}$
are bijections.
\end{Prop}

\begin{proof}
Let $s_0$ and $s_1$ be elements of $\F(X)$ whose restrictions to $C$
both equal $s_C$. Since $\F$ is excisive there exists an element
$s_{0,1}$ in $\F(X\amalg_C X)$ that restricts to $s_0$ and $s_1$ on
the two copies of $X$ in $X\amalg_C X$. If $\nu(s_0)$ and $\nu(s_1)$
are concordant relative to $C$ there exists an element
$s'\in\F'(X\times I/C)$ with $i_0^*s'=\nu(s_0)$ and
$i_1^*s'=\nu(s_1)$. Observe that the uniqueness part of excision for
$\F'$ implies that $s'|_{X\amalg_C X}=\nu(s_{0,1})$. Since the map
\[
\F[X\times I/C,X\amalg_C X;s_{0,1}]\xrightarrow{\nu_*} \F[X\times
I/C,X\amalg_C X;\nu(s_{0,1})]
\]
is surjective, we may lift $s'$ to an element $\bar{s}$ in
$\F(X\times I/C)$. This element $\bar{s}$ is a concordance between
$s_0$ and $s_1$ relative to $C$.

If $C\hookrightarrow X$ is a closed cofibration, then $X\amalg_C
X\hookrightarrow X\times I/C$ is also a closed cofibration.
\end{proof}

Let us now verify that the set of $Z_\bullet$-bundles satisfy strong
excision:

\begin{Prop}\label{prop:gluing}
Let $Y$ be the pushout of $f\co A\rightarrow B$ along an inclusion
$A\subseteq X$. Given $Z_\bullet$-bundles $\E$ over $X$ and $\E_B$
over $B$ with $f^*\E_B=\E|_A$, there exists a unique
$Z_\bullet$-bundle $\E'$ over $Y$ such that $\E'|_B=\E_B$ and $\E'$
pulls back to $\E$ over $X$.
\end{Prop}

\begin{proof}
We explicitly construct $\E'$. Let $\I$ and $\I_B$ denote the
indexing sets of the ordered open covers $\U$ and $\U_B$ associated
to $\E$ and $\E_B$. Define $\I'$ as the union $\I\cup\I_B$, and
extend the notation of $U_\alpha\in\U$ and $U^B_\alpha\in\U_B$ so
that $U_\alpha=\emptyset$ for $\alpha\not\in \I$ and
$U^B_\alpha=\emptyset$ for $\alpha\not\in \I_B$. The equation
$f^{-1}(U^B_\alpha)=U_\alpha\cap A$ is satisfied for all $\alpha\in
\I'$. Hence, we may define an ordered open cover
$\U'=\{U'_\alpha\}_{\alpha\in \I'}$ of $Y$ by declaring $U'_\alpha$
to be the image of $U_\alpha\amalg U^B_\alpha$ under the quotient
map $X\amalg B \rightarrow Y$. We see that the corresponding ordered
\v{C}ech complex $U'_\bullet$ is the pushout of
$U^B_\bullet\leftarrow U_\bullet\cap A\hookrightarrow U_\bullet$,
and we define the simplicial map $\phi^Y\co U^Y_\bullet\rightarrow
Z_\bullet$ as the uniquely defined extension of the simplicial maps
corresponding to $\E_B$ and $\E$.
\end{proof}

\begin{Cor}\label{cor:gluing}
Let $Y$ be the pushout of $f\co A\rightarrow B$ along an inclusion
$A\subseteq X$. Fix a $Z_{\bullet}$-bundle $\E_B$ over $B$, and let
$\E_A$ be $f^*\E_B$. Pullback of $Z_\bullet$-bundles induces a
bijection
\[
\Con_{Z_{\bullet}}(Y,B;\E_B)\xrightarrow{\cong}
\Con_{Z_{\bullet}}(X,A;\E_A).
\]
\end{Cor}

\begin{proof}
Keep $A$ fixed and let $\T$ be spaces containing $A$. There is a map
$\nu\co \F\rightarrow \F'$ of strongly excisive presheaves on $\T$
inducing $\Con_{Z_{\bullet}}(Y,B;\E_B)\rightarrow
\Con_{Z_{\bullet}}(X,A;\E_A)$. Here $\F(X)$ and $\F'(X)$ are the
sets of $Z_\bullet$-bundles over $X\amalg_f B$ and $X$ respectively,
and $\nu$ is pullback. By Proposition~\ref{Prop:surjimplyinj} all
$\nu_*$ are bijections if $\nu$ induces a surjection
$\F[X,C;\E_{C\amalg_f B}]\rightarrow \F'[X,C;\nu(\E_{C\amalg_f B})]$
for all inclusions $C\subseteq X$ in $\T$ and $Z_\bullet$-bundles
$\E_{C\amalg_f B}$ over $C\amalg_f B$. This follows immediately from
Proposition~\ref{prop:gluing}. Taking $C=A$ the conclusion follows.
\end{proof}

\begin{Lem}
Let $A$ be a subspace of $X$, and fix a $Z_\bullet$-bundle $\E_A$
over $A$. Then the inclusion $i_1\co X \hookrightarrow X\times I/A$
induces a bijection
\[
i_1^*\co \Con_{Z_\bullet}(X\times I/A,A;\E_A)
\xrightarrow{\cong}\Con_{Z_\bullet}(X,A;\E_A).
\]
\end{Lem}

\begin{proof}
Keep $A$ fixed and let $\T$ be spaces containing $A$. Two strongly
excisive presheaves on $\T$ are defined by letting $\F(X)$ and
$\F'(X)$ be the sets of $Z_\bullet$-bundles over $X\times I/ A$ and
$X$ respectively. Pullback by the inclusion $i_1\co X\hookrightarrow
X\times I/ A$ induces a map of presheaves $i_1^*\co \F\rightarrow
\F'$ inducing $\Con_{Z_{\bullet}}(X\times I/A,A;\E_A)\rightarrow
\Con_{Z_{\bullet}}(X,A;\E_A)$. By
Proposition~\ref{Prop:surjimplyinj} all $i_1^*$ are bijections if
$i_1^*$ induces a surjection $\F[X,C;\E_{C\amalg_f B}]\rightarrow
\F'[X,C;\nu(\E_{C\amalg_f B})]$ for all inclusions $C\subseteq X$ in
$\T$ and $Z_\bullet$-bundles $\E_{C\amalg_f B}$ over $C\amalg_f B$.
The projection map $\pi\co X\times I/A\rightarrow X$ induces on
concordance classes a map $\pi^*$ such that $i_1^*\pi^*$ is the
identity. Consequently $i_1^*$ is surjective. Taking $C=A$ completes
the proof.
\end{proof}

\begin{Cor}\label{Cor:htpyfunct}
$\Con_{Z_{\bullet}}(X,A;\E_A)$ is a homotopy functor of $X$ relative
to $A$.
\end{Cor}

\begin{Lem}\label{Lem:CEP}
Let $A\hookrightarrow B\hookrightarrow X$ be closed cofibrations.
Given $Z_\bullet$-bundles $\E$ over $X$ and $\E_B$ over $B$ that are
equal over $A$, and such that the restriction $\E|_B$ is concordant
to $\E_B$ relative to $A$, then there exists a $Z_\bullet$-bundle
$\E'$ over $X$ concordant to $\E$ relative to $A$ such that the
restriction $\E'|_B$ equals $\E_B$.
\end{Lem}

\begin{proof}
Represent the concordance between $\E|_B$ and $\E_B$ by a
$Z_\bullet$-bundle $\E''$ over $B\times I/A$.
Gluing $\E$ and $\E''$ we get a $Z_{\bullet}$-bundle $\bar{\E}$ over
$\left(X\amalg_{B} B\times I\right)/A$. By the
homotopy extension lifting property, there exists a map $j\co
X\rightarrow \left(X\amalg_{B} B\times I\right)/A$
such that in the following diagram
\[
\xymatrix@C=2cm{ B \ar[r]^-{i_1} \ar@{^(->}[d] &
\left(X\amalg_{B} B\times I\right)/A \ar[d]^{\text{retraction}} \\
X \ar@{=}[r] \ar@{.>}[ur]^j & X }
\]
the upper triangle commutes and the lower diagram commutes up to
homotopy relative to $B$. Let $\E'=j^*\bar{\E}$. Clearly
$\E'_B=\E_B$. Moreover, $j$ is homotopic to the inclusion $i_X\co
X\hookrightarrow \left(X\amalg_{B} B\times I\right)/A$ relative to
$A$. Hence,
\[
\left[\E'\right]=\left[j^*\bar{\E}\right]=\left[i_X^*\bar{\E}\right]=\left[\E\right]
\]
in $\Con_{Z_{\bullet}}(X,A;\E|_A)$.
\end{proof}

\begin{Cor}\label{Cor:HE}
The functor $\Con_{Z_{\bullet}}$ is half exact in the following
sense: let $A\hookrightarrow B\hookrightarrow X$ be closed
cofibrations, and fix some $Z_{\bullet}$-bundle $\E_B$ over $B$. Let
$\E_A=\E_B|_A$. In the sequence
\[
\Con_{Z_{\bullet}}(X,B;\E_B) \rightarrow
\Con_{Z_{\bullet}}(X,A;\E_A) \xrightarrow{\text{restriction}}
\Con_{Z_{\bullet}}(B,A;\E_A)
\]
the image of the first map is precisely the classes in
$\Con_{Z_{\bullet}}(X,A;\E_A)$ that restrict to $[\E_B]$ in
$\Con_{Z_{\bullet}}(B,A;\E_A)$.
\end{Cor}

\begin{Prop}\label{Prop:conlim}
For any relative CW-pair $(X,A)$ and $Z_\bullet$-bundle $\E_A$ over
$A$, the map
\[
\Con_{Z_{\bullet}}(X,A;\E_A)\xrightarrow{\cong}\lim_k\Con_{Z_{\bullet}}(X^k,A;\E_A)
\]
is a bijection. Here $X^k$ is the $k$-skeleton of $X$.
\end{Prop}

\begin{proof}
Keep $A$ fixed and let $\T$ be CW-complexes relative to $A$. There
is a map $\nu\co\F\rightarrow\F'$ of excisive presheaves on $\T$
inducing $\Con_{Z_{\bullet}}(X,A;\E_A)\rightarrow
\lim_k\Con_{Z_{\bullet}}(X^k,A;\E_A)$. Here $\F$ is the presheaf of
$Z_\bullet$-bundles restricted to $\T$, while $\F'(X)$ is the set of
sequences $(\E_0,\E_1,\ldots)$ where $\E_k$ is a $Z_\bullet$-bundle
over the $k$-skeleton $X^k$ such that the restriction of $\E_k$ to
$X^{k-1}$ is concordant to $\E_{k-1}$ relative to $A$. The map $\nu$
is given by sending $\E\in\F(X)$ to $(\E|_{X^0},\E|_{X^1},\ldots)$
in $\F'(X)$. By Proposition~\ref{Prop:surjimplyinj} all $\nu_*$ are
bijections if $\nu$ induces a surjection $\F[X,C;\E_{C}]\rightarrow
\F'[X,C;\nu(\E_{C})]$ for all inclusions $C\hookrightarrow X$ in
$\T$ and $Z_\bullet$-bundles $\E_{C}$ over $C$.

Fix $C\hookrightarrow X$ and $\E_C$, and take a sequence
$(\E_0,\E_1,\ldots)$ in $\F'(X)$ that restricts to $\nu(\E_C)$ in
$\F'(C)$. By definition of $\T$, see Example~\ref{Exa:Ts},
$C\hookrightarrow X$ is the inclusion of a subcomplex. Hence
$C^k\hookrightarrow X^{k-1}\cup C^k \hookrightarrow X^k$ are closed
cofibrations. By induction we construct a sequence $\E'_k$ over
$X^k$ such that the restriction of $\E'_k$ to $X^{k-1}$ is equal to
$\E'_{k-1}$, and $\E'_k$ is concordant to $\E_k$ relative to $C^k$.
We start by taking $\E'_0=\E_0$, and the inductive step uses
Lemma~\ref{Lem:CEP} with $\E_C|_{C^k}$ over $C^k$, the union of
$\E'_{k-1}$ and $\E_C|_{C^k}$ over $X^{k-1}\cup C^k$ and $\E_k$ over
$X^k$. This constructs an honest $Z_\bullet$-bundle $\colim_k\E'_k$
over $X$ mapping to the concordance class represented by the
original sequence $(\E_0,\E_1,\ldots)$, i.e. all $\nu_*$ are
surjective.
\end{proof}

So far we have studied how $\Con_{Z_\bullet}(X,A;\E_A)$ depend on the space $X$.
The next step is to understand what happens when $\E_A$ varies. To see this we
now define a category as follows:

\begin{Def}
Define $\mathscr{C}_{Z_\bullet}(A)$ to be the category with objects
all $Z_\bullet$-bundles over $A$, and with morphisms from $\E_A$ to
$\E'_A$ being the set $\Con_{Z_\bullet}(A\times I,A\amalg
A;\E_A\amalg\E'_A)$. The source and target of a morphism represented
by $\E$ are given as the restrictions $i_0^*\E$ and $i_1^*\E$
respectively. Composition is defined by gluing.
\end{Def}

\begin{Lem}\label{lem:CAZgroupoid}
The category $\mathscr{C}_{Z_\bullet}(A)$ is a groupoid.
\end{Lem}

\begin{proof}
We will show that a morphism represented by an arbitrary concordance
$\E$ has an inverse. Define $\E^{-1}$ by pulling back $\E$ over the
flip $A\times I\rightarrow A\times I$ sending $(a,t)$ to $(a,1-t)$.
The composition of $\E$ with $\E^{-1}$ equals the pullback of $\E$
over the fold $A\times I\rightarrow A\times I$ sending $(a,t)$ to
$(a,|2t-1|)$. A contraction $\theta$ from the fold to the map
$(a,t)\mapsto(a,1)$ is given by $\theta(a,s,t)=(a,(1-s)+s|2t-1|)$.
The pullback $\theta^*\E$ is a concordance relative to $A\amalg A$
from the identity on $i_1^*\E$ to the composition of $\E$ and
$\E^{-1}$.
\end{proof}

\begin{Prop}\label{Prop:changeofcollar}
For closed cofibrations $A\hookrightarrow X$ the expression
$\Con_{Z_\bullet}(X,A;-)$ defines a functor from
$\mathscr{C}_{Z_\bullet}(A)$ to the category of sets.
\end{Prop}

\begin{proof}
By the homotopy extension lifting property there exists a map
$j\co X\rightarrow X\amalg_{A\times\{0\}} A\times I$ such that
$j(a)=(a,1)$ for every $a\in A\subseteq X$ and the composition
\[
X\xrightarrow{j} X\amalg_{A\times\{0\}} A\times I \rightarrow X
\]
is homotopic to the identity on $X$ relative to $A$.

Given a $Z_\bullet$-bundle $\E$ representing a morphism from $\E_A$
to $\E'_A$ in $\mathscr{C}_{Z_\bullet}(A)$ and a $Z_\bullet$-bundle
$\E'$ representing a class of $\Con_{Z_\bullet}(X,A;\E_A)$ we can
glue them together and produce a $Z_\bullet$-bundle $\bar{\E}$ over
$X\amalg_{A\times\{0\}} A\times I$. The functor is defined by
declaring $[j^*\bar{\E}]$ in $\Con_{Z_\bullet}(X,A;\E'_A)$ to
represent the image of $[\E']$ under the map induced by the morphism
$[\E]$.
\end{proof}

Let $[\E]$ be a morphism from $\E_A$ to $\E'_A$ in
$\mathscr{C}_{Z_\bullet}(A)$. It follows from
Lemma~\ref{lem:CAZgroupoid} that $\E$ induces a bijection
\[
[\E]_*\co \Con_{Z_\bullet}(X,A;\E_A)\rightarrow\Con_{Z_\bullet}(X,A;\E'_A).
\]

An open cover $\U$ is \emph{finite and totally ordered} if its
indexing set $\I$ is finite and totally ordered. We call a
$Z_\bullet$-bundle \emph{finite} if its associated open cover is
finite and totally ordered. The finite $Z_\bullet$-bundles form a
strongly excisive presheaf $\F^f_{Z_\bullet}$ on spaces.
Consequently there is a notion of \emph{finite concordance}, and we
denote by
\[
\Con^f_{Z_\bullet}(X)\quad\text{and}\quad\Con^f_{Z_\bullet}(X,A;\E_A)
\]
the classes of finite $Z_\bullet$-bundles over $X$ modulo finite
concordance and finite $Z_\bullet$-bundles over $X$ equal to $\E_A$
over a subset $A$ modulo finite concordance relative to $A$
respectively. Over compact spaces there is no difference between
$\Con$ and $\Con^f$:

\begin{Lem}\label{Lem:finitecon}
Let $A$ be a compact subset of a compact space $X$. Then the natural
map
$\Con^{f}_{Z_\bullet}(X,A;\E_A)\rightarrow\Con_{Z_\bullet}(X,A;\E_A)$
is a bijection for all finite $Z_\bullet$-bundles $\E_A$ over $A$.
\end{Lem}

\begin{proof}
Let $\T$ be the full subcategory of compact space. Inclusion of
finite $Z_\bullet$-bundles into all $Z_\bullet$-bundles gives a map
of strongly excisive presheaves $\nu\co \F^f_{Z_\bullet}\rightarrow
\F_{Z_\bullet}$ on $\T$. By Proposition~\ref{Prop:surjimplyinj} it
is sufficient to show that the map
$\Con^{f}_{Z_\bullet}(X,A;\E_A)\rightarrow\Con_{Z_\bullet}(X,A;\E_A)$
is surjective for all inclusions $A\hookrightarrow X$ in $\T$ and
all finite $\E_A$ over $A$.

Let $\E$ be a $Z_\bullet$-bundle over $X$ equal to $\E_A$ on $A$. By
compactness of $X$, the ordered open cover $\U$ associated to $\E$
has a finite subcover $\U^f$. Without loss of generality we may
assume that the restriction of $\U^f$ to $A$ equals the cover
associated to $\E_A$. Choose a total ordering on $\U^f$ that extends
the partial ordering. Let $\E^f$ be the restriction of $\E$ to the
smaller covering $\U^f$. Observe that the finite $Z_\bullet$-bundle
$\E^f$ is concordant to $\E$ relative to $A$.
\end{proof}

\begin{Prop}\label{Prop:finite}
Let $q\co Z_\bullet\rightarrow W_\bullet$ be a map of simplicial
spaces. If
\[
\Con^{f}_{Z_\bullet}(D^n,S^{n-1};\E_{S^{n-1}}) \rightarrow
\Con^{f}_{W_\bullet}(D^n,S^{n-1};q_*\E_{S^{n-1}})
\]
is surjective for all $n\geq0$ and finite $Z_\bullet$-bundles
$\E_{S^{n-1}}$ over $S^{n-1}$, then
\[
\Con_{Z_\bullet}(X,A;\E_A)\rightarrow \Con_{W_\bullet}(X,A;q_*\E_A)
\]
is a bijection for all pairs $(X,A)$ of the homotopy type of a
CW-pair and all $Z_\bullet$-bundles $\E_A$ over $A$.
\end{Prop}

\begin{proof}
By homotopy invariance, Corollary~\ref{Cor:htpyfunct}, it is
sufficient to consider relative CW-complexes $(X,A)$. Fix $A$ and
let $\T$ be the category of CW-complexes relative to $A$. The
simplicial map $q$ gives a map
$q_*\co\F_{Z_\bullet}\rightarrow\F_{W_\bullet}$ of strongly excisive
presheaves on $\T$ inducing $\Con_{Z_\bullet}(X,A;\E_A) \rightarrow
\Con_{W_\bullet}(X,A;q_*\E_A)$. By
Proposition~\ref{Prop:surjimplyinj} this map of concordance classes
is a bijection if
$\F_{Z_\bullet}[X,C;\E_C]\rightarrow\F_{W_\bullet}[X,C;q_*\E_C]$ is
surjective for all subcomplex inclusions $C\subseteq X$ and
$Z_\bullet$-bundles $\E_C$ over $C$.

By Proposition~\ref{Prop:conlim} we can represent an element $\E'$
of $\Con_{W_\bullet}(X,C;q_*\E_C)$ as a sequence
$(\E'_0,\E'_1,\ldots)$ where $\E'_k$ is a $W_\bullet$-bundle over
$X^k$, such that $\E'_k|_{X^{k-1}}$ is concordant to $\E'_{k-1}$
relative to $C^{k-1}$ and $\E'_k|_{C^k}=(q_*\E_C)|_{C^k}$. By
induction on $k$ we will construct a sequence $(\E_0,\E_1,\ldots)$
where $\E_k$ is a $Z_\bullet$-bundle over $X^k$, such that
$\E_{k}|_{X^{k-1}}=\E_{k-1}$, $\E_k|_{C^k}=\E_C|_{C^k}$, and
$q_*\E_k$ is concordant to $\E'_k$ relative to $C^k$. This new
sequence will represent an element $\E$ of
$\Con_{Z_\bullet}(X,C;\E_C)$ mapping to the class $[\E']$ under
$q_*$.

We start the induction at the $(-1)$-skeleton, i.e.
$X^{-1}=C^{-1}=A$, $\E_{-1}=\E_A$, and $\E'_{-1}=q_*\E_A$. By the
induction hypothesis there is a concordance between
$\E'_k|_{X^{k-1}\cup C^k}$ and $q_*(\E_{k-1}\cup \E_C|_{C^k})$
relative to $C^k$. By Lemma~\ref{Lem:CEP} we can find a
$W_\bullet$-bundle $\tilde{\E}'_k$ over $X^k$ concordant to $\E'_k$
relative to $C^k$ such that $\tilde{\E}'_k|_{X^{k-1}\cup
C^k}=q_*(\E_{k-1}\cup \E_C|_{C^k})$. We claim that
\[
q_*\co \Con_{Z_\bullet}(X^k,X^{k-1}\cup C^k;\E_{k-1}\cup
\E_C|_{C^k})\rightarrow \Con_{W_\bullet}(X^k,X^{k-1}\cup
C^k;q_*(\E_{k-1}\cup \E_C|_{C^k}))
\]
is surjective. If the claim holds, then choose $\E_k$ to be any
$Z_\bullet$-bundle over $X^k$ equal to $\E_{k-1}\cup \E_C|_{C^k}$ on
$X^{k-1}\cup{C^k}$ such that $q_*[\E_k]=[\tilde{\E}'_k]$. This
completes the inductive step.

To prove the claim observe that $X^k$ can be constructed from
$X^{k-1}\cup C^k$ by attaching $k$-cells. By
Corollary~\ref{cor:gluing} and the fact that $\Con$ turns a disjoint
union of spaces into a product, it is sufficient to prove that
\[
\Con_{Z_\bullet}(D^k,S^{k-1};\E_{S^{k-1}})\rightarrow
\Con_{W_\bullet}(D^k,S^{k-1};q_*\E_{S^{k-1}})
\]
is surjective for all $k$ and $Z_\bullet$-bundles $\E_{S^{k-1}}$
over $S^{k-1}$. By Lemma~\ref{Lem:finitecon} $\E_{S^{k-1}}$ is
concordant to a finite $Z_\bullet$-bundle. By
Proposition~\ref{Prop:changeofcollar} nothing is lost by assuming
that $\E_{S^{k-1}}$ actually is finite. Applying
Lemma~\ref{Lem:finitecon} we see that the claim holds if
\[
\Con^f_{Z_\bullet}(D^k,S^{k-1};\E_{S^{k-1}})\rightarrow
\Con^f_{W_\bullet}(D^k,S^{k-1};q_*\E_{S^{k-1}})
\]
is surjective, which is precisely the hypothesis of the proposition
we are proving.
\end{proof}

We now simplify by noticing that any ordered \v{C}ech complex
$U_\bullet$ has \emph{free degeneracies},
see~\cite[Definition~A.4]{DuggerIsaksen:04}. This means that
$U_\bullet$ is the image of a presimplicial space $N_\bullet$ under
the left adjoint to the forgetful functor from simplicial spaces to
presimplicial spaces. It is easy to see what $N_\bullet$ must be; we
define the $k$-simplices by
\[
N_k=\coprod_{\alpha_0<\cdots<\alpha_k} U_{\alpha_0\cdots\alpha_k},
\]
and we call $N_\bullet$ the \emph{non-degenerate part} of
$U_\bullet$. Define $L^kN_\bullet$ to be subspace of $N_k$
consisting of points in the image of some face map $d_i$, i.e.
$L^kN_\bullet=\bigcup_i d_i(N_{k+1})$

Because of the proposition above we are lead to study simplicial
maps whose domains are ordered \v{C}ech complexes of finite and
totally ordered covers $\U$. Let us make the necessary preparations
for the construction of maps $N_\bullet\rightarrow Z_\bullet$ by
downward induction on simplicial degree:

\begin{Lem}\label{le:CW}
Let $N_\bullet$ be the non-degenerate part of the ordered \v{C}ech
complex corresponding to a finite and totally ordered open covering.
The following properties are satisfied:
\begin{itemize}
\item[i)] There exists an integer $K$ so that $N_k=\emptyset$ for $k > K$.
\item[ii)] If $d_i(x)=d_j(y)\in N_k$, $i<j$ there is a $z\in N_{k+2}$
such that $y=d_i(z)$ and $x=d_{j+1}(z)$.
\item[iii)] The obvious map
\[
\coprod_{0\leq i \leq k+1}  N_{k+1} \xrightarrow{\coprod d_i}
L^kN_\bullet
\]
is a quotient map.
\end{itemize}
\end{Lem}

\begin{proof}
Let $K$ be the number of elements in the cover. The next two
statements are general properties of the non-degenerate part of an
ordered \v{C}ech complex.

Addressing the last statement, we consider a subset $V$ of $N_k$
such that all $d_i^{-1}(V)$ are open in $N_{k+1}$. We have to show
that $V\cap L^kN_\bullet$ is open. We can assume that $V$ is a
subset of some $U_{\alpha_0\cdots\alpha_k}$. Notice that all $V\cap
U_\beta$ are open for $\beta\neq\alpha_i$. The result follows since
$V\cap L^kN_\bullet$ equals the finite union
$\bigcup_{\beta\neq\alpha_i}\left(V\cap U_\beta\right)$.
\end{proof}

\begin{Rem}
In the Reedy formalism, see~\cite[Chapter~15]{Hirschhorn:03}, one
considers the category of diagrams, $\mathbf{X}$, indexed by a Reedy
category $\EuScript{C}$. Each object in $\EuScript{C}$ has a degree,
and the usual way of constructing a map $f\co
\mathbf{X}\rightarrow\mathbf{Y}$ between diagrams is by increasing
induction on the degree. Associated to a diagram $\mathbf{X}$ and an
index $\alpha\in\EuScript{C}$ there are the so called latching and
matching objects, $L_\alpha \mathbf{X}\rightarrow
\mathbf{X}_\alpha\rightarrow M_\alpha \mathbf{X}$. If $f$ is already
defined in degrees below $\alpha$, this fixes maps of latching and
matching objects, and an extension of $f$ to $\mathbf{X}_\alpha$
must respect this.

In our case the presimplicial indexing category clearly is a Reedy
category. However, the union $L^kN_\bullet\subseteq N_k$ is neither
a latching nor a matching object. This is not so strange; we are
preparing for downward induction on degree, and $L^kN_\bullet$ plays
the role analogous of a latching object. The analogy for the
matching object is the terminal space.
\end{Rem}

The following lemma plays the role of ``a homotopy extension
property'' for the inclusion $L^kN_\bullet\hookrightarrow N_k$. In a
vague sense this is dual to Segal's definition of a good simplicial
space,~\cite{Segal:74}.

\begin{Lem}\label{Lem:shrink}
Let $N_\bullet$ be the non-degenerate part of the ordered \v{C}ech
complex corresponding to a finite and totally ordered open covering
$\U$ of a normal space $X$ and let $Z$ be a topological space. Given
$g\co N_k\rightarrow Z$ and $G\co L^kN_\bullet\times I\rightarrow Z$
with $g(x)=G(x,0)$ for all $x$ in $L^kN_\bullet$, there exists a
shrinking $\U^\epsilon$ of $\U$, inducing an inclusion
$N^\epsilon_\bullet\hookrightarrow N_\bullet$, and a homotopy
$G^\epsilon\co N^\epsilon_k\times I\rightarrow Z$ such that
$g(x)=G^\epsilon(x,0)$ for all $x\in N^\epsilon_k$ and
$G(x,t)=G^\epsilon(x,t)$ for all $x\in L^kN^{\epsilon}_\bullet$ and
$t\in I$. Moreover, if $t\mapsto G(x,t)$ is constant for some $x\in
N^{\epsilon}_k\cap L^kN_\bullet$, then $t\mapsto G^\epsilon(x,t)$ is
also constant.
\end{Lem}

\begin{proof}
Let $K$ be the number of elements in the cover. Choose a partition
of unity $\{\psi_\alpha\}$ subordinate to $\U$. For $\epsilon>0$
define $\U^\epsilon$ be the family of open sets
$\{\psi_{\alpha}^{-1}\left(\epsilon,1\right]\}$. If $\epsilon$ is
small, $\epsilon<\frac{1}{K}$, then $\U^\epsilon$ is also a cover.
Let $\psi\co N_k\rightarrow I$ be the map given on each intersection
$U_{\alpha_0\cdots\alpha_k}$ as the sum
$\sum_{\beta\neq\alpha_i}\psi_\beta$. Notice that the support of
$\psi$ is contained in $L^kN_\bullet$, and that
$\psi^{-1}\left(\epsilon,1\right]$ equals $L^kN_\bullet^\epsilon$.

Define a homotopy $\hat{G}^\epsilon\co N_k\times I\rightarrow Z$ by
the formula:
\[
\hat{G}^\epsilon(x,t)=\begin{cases} g(x) &\text{ if $x\not\in
L^kN_\bullet$,}\\
G(x,t) &\text{ if $x\in L^kN_\bullet$ and
$\psi(x)\geq\epsilon$, and}\\
G(x,\frac{\psi(x)}{\epsilon}t) &\text{ if $x\in L^kN_\bullet$ and
$\psi(x)\leq\epsilon$.}
\end{cases}
\]
Let $G^\epsilon$ be the restriction of $\hat{G}^\epsilon$ to
$N_k^\epsilon\times I$.
\end{proof}

By shrinking the cover we also get a ``lifting property'' up to
homotopy:

\begin{Prop}\label{Prop:he}
Let $h\co U_\bullet\rightarrow W_\bullet$ and $q\co Z_\bullet \to
W_\bullet$ be maps of simplicial spaces. Suppose that each $q_k\co
Z_k \to W_k$ is the inclusion of a strong deformation retract, and
that $U_\bullet$ is the \v{C}ech complex of a finite and totally
ordered open cover $\U$ of a normal space $X$. There exist a
shrinking $\U^\epsilon$ of $\U$, inducing an inclusion $i\co
U^{\epsilon}_\bullet\hookrightarrow U_\bullet$, together with
simplicial maps $g\co U^\epsilon_\bullet \to Z_\bullet$ and $F\co
U^\epsilon_\bullet\times I\to W_\bullet$ such that
\[
F_k(x,0)=h_k(i_k(x)) \quad\text{and}\quad F_k(x,1)=q_k(g_k(x))
\]
for all $k$ and $x\in U^{\epsilon}_k$. Moreover, if $h_k(x)$ is
contained in the image of $q_k$, then $t\mapsto F_k(x,t)$ is
constant.
\end{Prop}

\begin{proof}
It is enough to consider the non-degerenate part $N_\bullet$ of
$U_\bullet$. We construct $g$ and $F$ by downward induction on $n$.
Assume that there exists shrinking $\U^{\epsilon}$ of $\U$, inducing
$i\co N^{\epsilon}_\bullet\rightarrow N_\bullet$ on the
non-degenerate part of the ordered \v{C}ech complexes, and a
presimplicial map
\[
F^n\co  N^{\epsilon}_\bullet\times I\to W_\bullet
\]
such that $F^n_k(x,0)=h_k(i_k(x))$ and $F^n_k(x,1)\in q_k(Z_k)$ for
all $k>n$ and $x\in N^{\epsilon}_k$. To complete the inductive step
and define $F^{n-1}$, it suffices to construct a shrinking
$\U^{\epsilon'}$ of $\U^{\epsilon}$ and a presimplicial map $G\co
N^{\epsilon'}_\bullet \times I \to W_\bullet$ with the following
properties:
\begin{itemize}
\item[] $G_k(x,0)=F^n_k(x,1)$ for all $k$ and $x\in N^{\epsilon'}_k$,
\item[] $G_k(x,t)=F^n_k(x,1)$ if $k>n$, $x\in N^{\epsilon'}_k$, and
$t\in I$, and
\item[] $G_n(x,1)\in q_n(Z_n) \subset W_n$ for all $x\in
N^{\epsilon'}_k$.
\end{itemize}
Then we obtain $F^{n-1}$ by composing the homotopy $F^n$ with the
homotopy $G$.

We construct $G$ by downward induction on $k$. For $k>n$, the
conditions already define $G_k$. We start the induction at $k=n$ and
take $\U^{\epsilon'}=\U^{\epsilon}$. Let $H\co W_n \times I \to W_n$
be the deformation retraction into $q_n(Z_n)$, and define $G_n$ by
the formula $G_n(x,t)=H(F^n_n(x,1),t)$. Observe that the image of
$G_{n+1}$ is contained in $q_{n+1}(Z_{n+1})$, whence the definition
of $G_n$ ensures commutativity of the diagram
\[
\xymatrix@C=1.5cm{ N^{\epsilon'}_{n+1}\times
I\ar[rr]^{d_i}\ar[d]^{G_{n+1}}
&&N^{\epsilon'}_n\times I \ar[d]^{G_n}\\
q_{n+1}(Z_{n+1})\ar[r]^{\subseteq} &W_{n+1}\ar[r]^{d_i} &W_n. }
\]
Thus $G$ satisfies the presimplicial relations in simplicial degrees
$k\geq n$. It remains to defined $G_k$ for $k<n$.

For each step in the inductive construction of $G$ we will shrink
the cover $\U^{\epsilon'}$ to a smaller cover $\U^{\epsilon''}$. So
we implicitly restrict the domains of the $G_k$'s already
constructed. Assuming that $G_{k+1}\co N^{\epsilon'}_{k+1}\times
I\rightarrow W_{k+1}$ already has been defined we proceed to define
$G_k$.

Because of Lemma~\ref{le:CW}.ii), there is a well defined set map
${\tilde G_k\co L^k N^{\epsilon'}_{\bullet}\times I\to W_k}$ such
that the following diagram commutes
\[
\xymatrix@C=1.5cm{ N^{\epsilon'}_{k+1}\times
I\ar[r]^-{d_i}\ar[d]^{G_{k+1}}&
L^k N^{\epsilon'}_{\bullet}\times I \ar[d]^{\tilde G_{k}}\\
W_{k+1}\ar[r]^{d_i}& W_{k} }
\]
for all $d_i$. According to Lemma~\ref{le:CW}.iii), $\tilde G_k$ is
continuous. By Lemma~\ref{Lem:shrink}, we can find a shrinking
$\U^{\epsilon''}$ of the cover $\U^{\epsilon'}$, inducing an
inclusion $N^{\epsilon''}_\bullet\subseteq N^{\epsilon'}_\bullet$,
together with a map $G_k\co N_k^{\epsilon''}\times I\rightarrow W_k$
such that $G_k(x,0)=F^n_k(x,1)$ for $x\in N^{\epsilon''}_k$, and
$G_k(x,t)={\tilde G_{k}}(x,t)$ for $x\in L^kN^{\epsilon''}_\bullet$
and $t\in I$.

This completes the inductive step, and the proof of everything
except the last part. However, tracing the proof it is easily seen
that all homotopies fix points $x\in N^{\epsilon}$ having $h_k(x)\in
q_k(Z_k)$.
\end{proof}

\begin{Cor}\label{Cor:consdr}
Let $X$ be a space and $A$ a subspace such that $(X,A)$ has the
homotopy type of a CW pair. Let $q\co Z_\bullet\hookrightarrow
W_\bullet$ be a simplicial map which in each simplicial degree is
the inclusion of a strong deformation retract. Let $\E_A$ be any
$Z_\bullet$-bundle over $A$. Then
\[
\Con_{Z_\bullet}(X,A;\E_A)\rightarrow \Con_{W_\bullet}(X,A;q_*\E_A)
\]
is a bijection.
\end{Cor}

\begin{proof}
According to Proposition~\ref{Prop:finite} it is sufficient to show
that the induced map
\[
q_*\co \Con^{f}_{Z_\bullet}(D^n,S^{n-1};\E_{S^{n-1}})\rightarrow
\Con^{f}_{W_\bullet}(D^n,S^{n-1};q_*\E_{S^{n-1}})
\]
is surjective when fixing an arbitrary integer $n$ and a finite
$Z_\bullet$-bundle $\E_{S^{n-1}}$ over the boundary sphere. Let
$\E'$ be a finite $W_\bullet$-bundle over $D^n$ whose restriction to
$S^{n-1}$ equals $q_*\E_{S^{n-1}}$. Let $\U$ denote the
corresponding finite and totally ordered open cover, and let $h\co
U_\bullet\rightarrow W_\bullet$ be the simplicial map associated to
$\E'$. Applying Proposition~\ref{Prop:he} we get a shrinking
$\U^{\epsilon}$ of $\U$, an induced inculsion of ordered \v{C}ech
complexes $U_\bullet^{\epsilon}\hookrightarrow U_\bullet$, a
simplicial map $g\co U_\bullet^{\epsilon}\rightarrow Z_\bullet$, and
a homotopy $F\co U_\bullet^\epsilon\times I\rightarrow W_\bullet$
between the restriction of $h$ to $U_\bullet^{\epsilon}$ and $qg$.

The shrinking $\U^\epsilon$ of $\U$ gives us a family of inclusions
$\{U^\epsilon_\alpha\hookrightarrow U_\alpha\}$. For each $\alpha$
define $U''_\alpha$ to be the open subset $U^\epsilon_\alpha\times I
\cup U_\alpha \times \left[0,\frac{1}{2}\right)$ of $D^n\times I$.
Together all the $U''_\alpha$ constitute a cover $\U''$ of
$D^n\times I$. Let $\phi''\co U''_\bullet\rightarrow W_\bullet$ be
the simplicial map given in simplicial degree $k$ as
\[
\phi''(x,t)=\begin{cases} h(x) & \text{for $x\in U_k$ and $0\leq t <
\frac{1}{2}$, and}\\ F(x,2t-1) &\text{for $x\in U_k^\epsilon$ and
$\frac{1}{2}\leq t\leq 1$.}\end{cases}
\]
This simplicial map explicitly defines a $W_\bullet$-bundle $\E''$
over $D^n\times I$. The restriction of $\E''$ to $D^n\times\{0\}$
lifts to a $Z_\bullet$-bundle $\E$ represented by the simplicial map
$g\co U^\epsilon_\bullet\rightarrow Z_\bullet$. Moreover, the
restriction of $\E''$ to $S^{n-1}\times I$ lifts to a concordance
between the restriction $\E|_{S^{n-1}}$ and our fixed
$Z_\bullet$-bundle $\E_{S^{n-1}}$. Denote this concordance by
$\E_{S^{n-1}\times I}$.

Now consider the diagram
\[
\xymatrix{ \Con^f_{Z_\bullet}(D^n,S^{n-1};\E|_{S^{n-1}})
\ar[d]^{[\E_{S^{n-1}\times I}]_*} \ar[r] &
\Con^f_{W_\bullet}(D^n,S^{n-1};q_*\E|_{S^{n-1}}) \ar[d]^{[q_*\E_{S^{n-1}\times I}]_*} \\
\Con^f_{Z_\bullet}(D^n,S^{n-1};\E_{S^{n-1}}) \ar[r]^{q_*} &
\Con^f_{W_\bullet}(D^n,S^{n-1};q_*\E_{S^{n-1}}). }
\]
The vertical maps, see Proposition~\ref{Prop:changeofcollar},
corresponds to gluing the concordances $\E_{S^{n-1}\times I}$ and
$q_*\E_{S^{n-1}\times I}$ along the boundary sphere. They are
bijections by Lemma~\ref{lem:CAZgroupoid}. We started by picking a
$W_\bullet$-bundle $\E'$ representing an element from the lower
right corner. The $Z_\bullet$-bundle $\E$ represents an element in
$\Con^f_{Z_\bullet}(D^n,S^{n-1};\E|_{S^{n-1}})$, and $\E''$ produces
a concordance between the image of $q_*[\E_{S^{n-1}\times I}]_*\E$
and $\E'$. Hence $q_*$ is surjective.
\end{proof}

We say that two simplicial maps $p,q\co Z_\bullet\rightarrow
W_\bullet$ are \emph{topologically homotopic} if there exists a
simplicial map $F\co Z_\bullet\times I\rightarrow W_\bullet$ with
$F_k(x,0)=p_k(x)$ and $F_k(x,1)=q_k(x)$ for each $k$ and all $x\in
Z_k$.

\begin{Prop}\label{Prop:conhtpyfunctinZ}
Let $A\hookrightarrow X$ be a closed cofibration, and suppose that
$p,q\co Z_\bullet \rightarrow W_\bullet$ are topologically
homotopic. Given any $Z_\bullet$-bundle $\E_A$ over $A$ there exists
a $W_\bullet$-bundle $\E'_{A\times I}$ over $A\times I$, restricting
to $p_*\E_A$ and $q_*\E_A$ on $A\times\{0\}$ and $A\times\{1\}$
respectively, making the following triangle commute
\[
\xymatrix{ & \Con_{Z_\bullet}(X,A;\E_A) \ar[dl]_{p_*} \ar[dr]^{q_*}
& \mbox{}\\ \Con_{W_\bullet}(X,A;p_*\E_A) \ar[rr]^{[\E'_{A\times
I}]_*} && \Con_{W_\bullet}(X,A;q_*\E_A).}
\]
\end{Prop}

\begin{proof}
Let $F\co Z_\bullet\times I\rightarrow W_\bullet$ denote a homotopy
between $p$ and $q$.

Given an ordered open cover $\U$ of $X$ we define a cover $\U'$ of
$X\times I$ consisting of the open sets $U'_\alpha=U_\alpha\times I$
where $U_\alpha\in\U$. The corresponding ordered \v{C}ech complexes
satisfy $U'_\bullet=U_\bullet\times I$.

Let $\E$ be a $Z_\bullet$-bundle over $X$ equal to $\E_A$ on $A$,
and consider the associated simplicial map $\phi\co
U_\bullet\rightarrow Z_\bullet$. Subordinate to the cover $\U'$ we
define a $W_\bullet$-bundle $\E'$ over $X\times I$ having as its
simplicial map the composition
\[
U'_\bullet=U_\bullet\times I\xrightarrow{\phi\times\id}
Z_\bullet\times I \xrightarrow{F} W_\bullet.
\]
The restriction of $\E'$ to $A\times I$ defines the concordance
$\E'_{A\times I}$ between $p_*\E_A$ and $q_*\E_A$. Moreover, $\E'$
yields a concordance between $[\E'_{A\times I}]_*p_*\E$ and $q_*\E$.
\end{proof}

\begin{Thm}\label{Thm:ConDegwiseHEQ}
Let $A\hookrightarrow X$ be a closed cofibration such that $(X,A)$
has the homotopy type of a CW pair. Let $f\co Z_\bullet\rightarrow
W_\bullet$ be a simplicial map such that each $f_k$ is a homotopy
equivalence. Let $\E_A$ be any $Z_\bullet$-bundle over $A$. Then $f$
induces a bijection
\[
f_*\co\Con_{Z_\bullet}(X,A;\E_A)\rightarrow
\Con_{W_\bullet}(X,A;f_*\E_A).
\]
\end{Thm}

\begin{proof}
Let $M_\bullet$ be the degreewise mapping cylinder of $f$. There are
inclusions of degreewise strong deformation retracts $i\co Z_\bullet
\hookrightarrow M_\bullet$ and $j\co W_\bullet \hookrightarrow
M_\bullet$. Moreover, the map $jf$ is topologically homotopic to
$i$. According to Proposition~\ref{Prop:conhtpyfunctinZ} there
exists a concordance $\E'_{A\times I}$ of $M_\bullet$-bundles
between $i_*\E_A$ and $j_*f_*\E_A$ making the following diagram
commute
\[
\xymatrix{ \Con_{Z_\bullet}(X,A;\E_A) \ar[d]_{i_*} \ar[r]^{f_*} &
\Con_{W_\bullet}(X,A;f_*\E_A) \ar[d]^{j_*} \\
\Con_{M_\bullet}(X,A;i_*\E_A) \ar[r]^{[\E'_{A\times I}]_*} &
\Con_{M_\bullet}(X,A;j_*f_*\E_A).}
\]
The vertical maps are bijections by Corollary~\ref{Cor:consdr}. The
bottom map is a bijection by Lemma~\ref{lem:CAZgroupoid}.
\end{proof}

In our next line of arguments we will show that filling a horn in
$Z_\bullet$ does not affect the concordance classes of
$Z_\bullet$-bundles.

A shrinking $\U'$ of $\U$ induces an inclusion of ordered \v{C}ech
complexes $i\co U'_\bullet\hookrightarrow U_\bullet$. Given a
simplicial map $U_\bullet\xrightarrow{\phi}Z_\bullet$, we define the
\emph{cover-restriction} of $\phi$ to $\U'$ the composition
$U'_\bullet\xrightarrow{i} U_\bullet\xrightarrow{\phi} Z_\bullet$.
Similarly, if $\E$ is the $Z_\bullet$-bundle associated to $\phi$,
then we denote the $Z_\bullet$-bundle associated to $\phi i$ by
$i^*\E$, and call it a \emph{cover-restriction} of $\E$.

\begin{Prop}\label{Prop:fillhorn}
Let $h\co \Lambda_\bullet^{n,k}\rightarrow Z_\bullet$ be a $k$-horn
of dimension $n$ in a simplicial space $Z_\bullet$, and define $q\co
Z_\bullet\hookrightarrow W_\bullet$ by filling this horn, i.e.
$W_\bullet$ is the pushout of
$\Delta_\bullet^n\hookleftarrow\Lambda_\bullet^{n,k}\xrightarrow{h}
Z_\bullet$. Let $\phi\co U_\bullet\rightarrow W_\bullet$ be a
simplicial map from the ordered \v{C}ech complex associated to a
finite and totally ordered open cover $\U$ of a normal space $X$.
Then there exists a zigzag chain of finite and totally ordered open
covers
\[
\U=\U_0\subseteq \U'_1\supseteq \U_1\subseteq \U'_2\supseteq
\ldots\subseteq \U'_{s-1}\supseteq \U_{s-1}\subseteq \U'_s \supseteq
\U_s
\]
together with simplicial maps $\phi_i$ and $\phi'_i$ from the
various ordered \v{C}ech complexes into $W_\bullet$ such that each
$\phi_i$ is equal to the cover-restriction of both $\phi'_{i-1}$ and
$\phi'_{i}$, and the last simplicial map $\phi_s$ lifts through
$q\co Z_\bullet\rightarrow W_\bullet$. Moreover, if $A$ is a subset
of $X$ such that $\phi$ restricted to $\U\cap A$ already lifts to
$Z_\bullet$, then we may assume that the restriction to $A$ of all
$\phi'_i$ and $\phi_i$ also lifts to $Z_\bullet$.
\end{Prop}

At the center of this long and complicated proof there is an
unpolished gem; the emptying of the horn by modifications of the
cover.

\begin{proof}
It is enough to consider presimplicial maps from the non-degenerate
part of the various ordered \v{C}ech complexes.

Our focus will be to move $\phi\co N_\bullet\rightarrow W_\bullet$
away from the $k$th face of the filled horn. Recall that the
$r$-simplices of $\Delta^n_\bullet$ correspond to order-preserving
maps $\theta\co[r]\rightarrow [n]$. Denote by $w_\theta$ the
$r$-simplex in $W_r$ corresponding to $\theta$ under the simplicial
map $\Delta^n_\bullet\rightarrow W_\bullet$. Let $\Theta$ be the set
of all $\theta$ that factors as a surjective map
$[r]\twoheadrightarrow [n-1]$ followed by the coface map
$\delta_k\co[n-1]\rightarrow[n]$ that omits $k$. Here we allow $r$
to vary. Choose a total ordering on $\Theta$ such that
$\theta_1<\theta_2$ whenever the domain of $\theta_1$ contains fewer
elements than the domain of $\theta_2$. Each $\theta\in\Theta$
corresponds to a degenerate simplex above the $k$th face of
$\Delta^n_\bullet$. Furthermore, for all $\theta\in\Theta$ the
simplex $w_\theta$ is an isolated point of $W_r$.

The cover $\U$ has a finite and totally ordered indexing set $\I$.
An order-preserving injection $\rho\co[r]\hookrightarrow \I$
determines a subset $U_\rho=U_{\rho(0)}\cap\cdots\cap U_{\rho(r)}$.
Moreover, $U_\rho$ is a summand of $N_r$. If the presimplicial map
$\phi\co N_\bullet\rightarrow W_\bullet$ does not factor through
$q\co Z_\bullet\hookrightarrow W_\bullet$, then some degenerate
simplex above the $k$th face of $\Delta^n_\bullet$ lies in the image
of $\phi$. From now on reserve $\theta$ as a symbol for the maximal
element of $\Theta$ such that $w_\theta$ is contained in
$\phi_r(N_r)$ for some $r$.

Define the \emph{$\theta$-valence} of a presimplicial map $\phi\co
N_\bullet\rightarrow W_\bullet$ to be the number of order-preserving
injections $\rho\co[r]\hookrightarrow\I$ such that $\phi_r(U_\rho)$
contains $w_\theta$. We will soon construct a finite and totally
ordered cover $\U'$ containing $\U$ as a subcover together with an
extension $\phi'\co N'_\bullet\rightarrow W_\bullet$ of $\phi$ such
that there exists a shrinking $\U^\epsilon$ of $\U'$ where the
cover-restriction $\phi^\epsilon\co N_\bullet^\epsilon\rightarrow
W_\bullet$ of $\phi'$ has smaller $\theta$-valence than $\phi$.
Repeating this construction we will eventually reach a presimplicial
map of $\theta$-valence zero, i.e. the point $w_\theta$ is no longer
in its image. By downward induction on $\theta\in\Theta$ the
forthcoming construction will establish the lemma.

For such a maximal $\theta$ choose any $\rho$ with
$w_\theta\in\phi_r(U_\rho)$. Let $j$ be the largest integer such
that $\theta(j-1)< k$. Define the indexing set $\I'$ by inserting
into $\I$ a new index $\beta$ as the successor of $\rho(j-1)$, i.e.
we have $\I'=\I\cup\{\beta\}$. Since $w_\theta$ is an isolated point
of $W_r$ the subset of $U_\rho$ mapping to $\theta$ is both closed
and open, and we define $U_\beta=U_\rho\cap \phi^{-1}(w_\theta)$.
Now view $U_\beta$ as a subspace of $X$, and consider the
restriction of the presimplicial map $\phi$ to $U_\beta$. We claim
that there is a factorization of $\phi|_{U_\beta}$ as a
presimplicial map $N_\bullet\cap U_\beta\rightarrow
\Delta^n_\bullet$ followed by the simplicial map
$\Delta^n_\bullet\rightarrow W_\bullet$.

Let us verify this claim. Given an order-preserving injection
$\tau\co [t]\hookrightarrow \I$ let $\hat{t}$ be the smallest
integer such that there exists order-preserving injections
$\hat{\tau}$, $\mu$, and $\nu$ making the diagram
\[
\xymatrix{  & [t] \ar@{^{(}->}[dl]_{\nu} \ar@{^{(}->}[dr]^{\tau} &
\\ [\hat{t}] \ar@{^{(}->}[rr]^{\hat{\tau}} && \I \\ & [r]
\ar@{^{(}->}[ul]^{\mu} \ar@{^{(}->}[ur]_{\rho} }
\]
commute. Observe that $\hat{\tau}$, $\mu$, and $\nu$ are unique.
Informally, $\hat{\tau}$ is the union of $\tau$ and $\rho$.

By the constructions above $\phi_r$ maps any $x$ in $U_\rho\cap
U_\beta$ to $w_\theta$ in $W_r$. Suppose that $U_{\hat{\tau}}$ meets
$U_\beta$ and take some $x$ in $U_{\hat{\tau}}\cap U_\beta$. Observe
that $\phi_{\hat{t}}(x)\in W_{\hat{t}}$ is sent to $w_\theta$ by the
face operator $\mu^*$ corresponding to
$\mu\co[r]\rightarrow[\hat{t}]$. Inspecting the definition of
$W_\bullet$ we see that all simplices mapping to $w_\theta$ by some
face operator lies in the image of $\Delta^n_\bullet\rightarrow
W_\bullet$. Hence $\phi_{\hat{t}}(x)=w_{\hat{\theta}}$ for some
$\hat{\theta}\co[\hat{t}]\rightarrow[n]$. The equation
$\theta=\hat{\theta}\mu$ is satisfied, and for any $i\in[\hat{t}]$
outside the image of $\mu$ we must have $\hat{\theta}(i)=k$.
Otherwise we get a contradiction, i.e. some $\theta_2>\theta$ in
$\Theta$ would be in the image of $\phi$. Thus $\hat{\theta}$ is
uniquely determined. At last we see that $\phi_t$ maps any $x\in
U_\tau\cap U_\beta$ to $w_{\hat{\theta}\nu}$ in $W_t$. This
completes the verification of the claim.

Now define the cover $\U'$ as the union $\U\cap\{U_\beta\}$. Any extension
$\phi'\co N'_\bullet\rightarrow W_\bullet$ is completely determined by the
restrictions of $\phi'_{t+1}$ to $U_{\tau'}$ for the order-preserving injections
$\tau'\co[t+1]\hookrightarrow\I'$ containing $\beta$ in their image.
For any such $\tau'$ define $\tau$ and $\delta_l$ by the pullback diagram
\[
\xymatrix{ [t] \ar@{^{(}->}[r]^{\delta_l} \ar@{^{(}->}[d]_{\tau} & [t+1] \ar@{^{(}->}[d]^{\tau'} \\
\I \ar@{^{(}->}[r] & \I'.}
\]
Applying the notation introduced above, we define $\theta'\co[t+1]\rightarrow[n]$
by $\theta'(l)=k$ and $\theta'\delta_l=\hat{\theta}\nu$. For $x$ in $U_{\tau'}=U_\tau\cap U_\beta$
define $\phi'_{t+1}$ by the formula
\[
\phi'_{t+1}(x)=w_{\theta'}.
\]
Since
$d_l\phi'_{t+1}(x)=w_{\theta'\delta_l}=w_{\hat{\theta}\nu}=\phi_t(x)$
it follows that $\phi'$ satisfies the presimplicial identities.
Furthermore, $\phi$ and $\phi'$ have the same $\theta$-valence since
$w_\theta$ is not contained in any $\phi'_{t+1}(U_{\tau'})$. What
remains to be proved is that we can find a shrinking $\U^\epsilon$
of $\U'$ such that the cover-restriction $\phi^\epsilon$ of $\phi'$
has $\theta$-valence smaller than the $\theta$-valence of $\phi$ and
$\phi'$.

Since $\U$ is a finite open cover of a normal space, we can choose a
partition of unity $\{\psi_\alpha\}$ subordinate to the cover $\U$.
Let $C$ be the closure of $U_\beta$ in $X$. Define $B$ to be the set
of all points $x\in X\smallsetminus U_\beta$ with
$\psi_{\rho(i)}(x)=\frac{1}{r+1}$ for $i=0,\ldots,r$. Because
$U_\rho\smallsetminus U_\beta$ is an open neighborhood of $B$ in
$X$, we see that $B$ and $C$ are disjoint closed subsets in $X$. By
the Tietze extension theorem there exists a continuous function
$\epsilon\co X\rightarrow \left[0,\frac{1}{r+1}\right]$ such that
$\epsilon(x)=0$ for all $x\in C$ and $\epsilon(x)=\frac{1}{r+1}$ for
all $x$ in $B$. Informally, we think about the points $x$ with all
$\psi_{\rho(i)}(x)=\frac{1}{r+1}$ as the bad guys. The role of
$\epsilon$ is to separate the bad points inside $U_\beta$ from those
outside.

We will now define a shrinking $\U^\epsilon$ of $\U'$ by shrinking
the $U_{\rho(i)}$ while leaving the other open sets unchanged.
Define $U^\epsilon_{\rho(i)}$ to be the set of all points $x\in
U_{\rho(i)}$ satisfying the inequality
\[
r\cdot \psi_{\rho(i)}(x)+\epsilon(x) > \psi_{\rho(0)}(x)+ \cdots +
\psi_{\rho(i-1)}(x) + \psi_{\rho(i+1)}(x) + \cdots +
\psi_{\rho(r)}(x).
\]
Explicitly, the cover $\U^\epsilon$ is given as the union
$\{U^\epsilon_{\rho(i)}\}_{i=0}^r \cup \{U_\alpha\}_{\alpha\neq
\rho(i),\beta} \cup \{U_\beta\}$. Let us verify that $\U^\epsilon$
really is a cover of $X$. Take an arbitrary point $x$ in $X$.
\begin{itemize}
\item[i)] If $\sum_{i=0}^r\psi_{\rho(i)}(x)<1$, then $x\in U_\alpha$ for some
$\alpha\in \I$ outside the image of $\rho$.
\item[ii)] If $\psi_{\rho(j)}(x)>\frac{1}{r+1}$ for some $j$, then
\[
(r+1)\cdot \psi_{\rho(j)}(x)>1\geq \sum_{i=0}^r\psi_{\rho(i)}(x).
\]
Consequently, $x$ lies in $U^\epsilon_{\rho(j)}$.
\item[iii)] If $x\in B$, then $\psi_{\rho(i)}(x)=\frac{1}{r+1}$ for
$i=0,\ldots,r$ and $\epsilon(x)=\frac{1}{r+1}$. Hence,
\[
r\cdot\psi_{\rho(0)}(x)+\epsilon(x)=1>\frac{r}{r+1}=\sum_{i=1}^r\psi_{\rho(i)}(x).
\]
So, $x$ lies in $U^\epsilon_{\rho(0)}$.
\item[iv)] Otherwise, all $\psi_{\rho(i)}(x)=\frac{1}{r+1}$ and
$x\not\in B$, whence $x\in U_\beta$.
\end{itemize}
This shows that $\U^\epsilon$ indeed is a cover of $X$.

Define $\phi^\epsilon\co N^\epsilon_\bullet\rightarrow W_\bullet$ to
be the cover-restriction of $\phi'$ over $\U^\epsilon\subseteq \U'$.
We claim that the $\theta$-valence of $\phi^\epsilon$ is less than
the $\theta$-valence of $\phi'$. More precisely, $w_\theta$ is
contained in $\phi'(U'_\rho)$, but not in
$\phi^\epsilon(U^\epsilon_\rho)$. To see this, assume for a
contradiction that $\phi^\epsilon(x)=w_\theta$ for some $x\in
U^\epsilon_\rho$. By definition of $U_\beta$, we see that such $x$
must be contained in $U_\beta$, and consequently $\epsilon(x)=0$.
Moreover $x$ satisfies each of the inequalities defining
$U^\epsilon_{\rho(i)}$ for $i=0,\ldots,r$. Adding all these
inequalities we get
\[
\sum_{j=0}^r( r\cdot \psi_{\rho(j)}(x)) + r\cdot \epsilon(x)
> \sum_{j=0}^r \sum_{i\neq j} \psi_{\rho(i)}(x),
\]
which contradicts $\epsilon(x)=0$.

Addressing the last statement of the lemma, we consider a subset
$A\subseteq X$ such that the restriction of $\phi$ to $N_\bullet\cap
A$ already lifts through $q\co Z_\bullet \hookrightarrow W_\bullet$.
Let $\theta$, $\rho$, and $U_\beta$ be defined as above. Notice that
$A\cap U_\beta$ must be empty, otherwise $w_\theta$ would be in the
image of $\phi|_A$. Consequently, the restriction of
$\U'=\U\cup\{U_\beta\}$ to $A$ equals $\U\cap A$, and similarly
$\phi|_A=\phi'|_A$. It follows that $\phi'|_A$ lifts through
$Z_\bullet$. Since $\phi^\epsilon|_A$ is a cover-restriction of
$\phi'|_A$, also the restriction of $\phi^\epsilon$ to $A$ lifts
through $Z_\bullet$, and we are done.
\end{proof}

\begin{Cor}\label{Cor:confillhorn}
Let $X$ be a space and $A$ a subspace such that $(X,A)$ has the
homotopy type of a CW pair. Let $h\co
\Lambda_\bullet^{n,k}\rightarrow Z_\bullet$ be a horn in a
simplicial space $Z_\bullet$, and define $q\co
Z_\bullet\hookrightarrow W_\bullet$ by filling this horn, i.e. let
$W_\bullet$ be the pushout of
$\Delta_\bullet^n\hookleftarrow\Lambda_\bullet^{n,k}\xrightarrow{h}
Z_\bullet$. Let $\E_A$ be any $Z_\bullet$-bundle over $A$. Then
\[
\Con_{Z_\bullet}(X,A;\E_A)\rightarrow \Con_{W_\bullet}(X,A;q_*\E_A)
\]
is a bijection.
\end{Cor}

\begin{proof}
By Proposition~\ref{Prop:finite} it is enough to show that $q$
induces surjections
\[
q_*\co \Con^{f}_{Z_\bullet}(D^m,S^{m-1};\E_{S^{m-1}})\rightarrow
\Con^{f}_{W_\bullet}(D^m,S^{m-1};q_*\E_{S^{m-1}})
\]
for arbitrary integers $m$ and finite $Z_\bullet$-bundles
$\E_{S^{m-1}}$ over the boundary sphere. Let $\E_0$ be a finite
$W_\bullet$-bundle over $D^m$ whose restriction to $S^{m-1}$ equals
$q_*\E_{S^{m-1}}$. Let $\U_0$ denote the corresponding finite and
totally ordered open cover, and let $\phi_0\co
U^0_\bullet\rightarrow W_\bullet$ be the simplicial map associated
to $\E_0$. By Proposition~\ref{Prop:fillhorn} there exists a zigzag
chain of finite and totally ordered open covers
\[
\U_0\subseteq \U'_1\supseteq \U_1\subseteq \U'_2\supseteq
\ldots\subseteq \U'_{s-1}\supseteq \U_{s-1}\subseteq \U'_s \supseteq
\U_s
\]
together with simplicial maps $\phi_i$ and $\phi'_i$ from the
various ordered \v{C}ech complexes into $W_\bullet$ such that each
$\phi_i$ is equal to the cover-restriction of both $\phi'_{i-1}$ and
$\phi'_{i}$, and the last simplicial map $\phi_s$ lifts through
$q\co Z_\bullet\rightarrow W_\bullet$. Moreover, the restriction to
$A$ of each $\phi'_i$ and $\phi_i$ lifts to $Z_\bullet$.

We will now build a concordance $\E_{D^m\times I}$ between
$W_\bullet$-bundles such that $i_0^*\E_{D^m\times I}=\E_0$ while the
restrictions of $\E_{D^m\times I}$ to $D^m\times\{1\}$ and
$S^{m-1}\times I$ lifts to $Z_\bullet$-bundles. Let $\I'_i$ be the
indexing set of $\U'_i$. The elements of $\U_i$ and $\U'_i$ will be
denoted by $U_{i,\alpha}$ and $U'_{i,\alpha}$ respectively, and we
extend the notation by letting $U_{i,\alpha}=\emptyset$ and
$U'_{i,\alpha}=\emptyset$ whenever the index $\alpha$ lies outside
the respective indexing sets. Now define a cover $\U$ of $X\times I$
indexed by the union $\I=\bigcup_{i=1}^{s}\I'_i$ by defining the
open set $U_\alpha$ to be the union
\[
\left(\bigcup_{i=0}^sU_{i,\alpha}\times\left\{\frac{i}{s}\right\}\right)
\cup \left(\bigcup_{i=1}^s
U'_{i,\alpha}\times\left(\frac{i-1}{s},\frac{i}{s}\right)\right).
\]
We define a simplicial map $\phi$ from the ordered \v{C}ech complex
of $\U$ into $W_\bullet$ by defining each restriction
$\phi|_{X\times\left\{\frac{i}{s}\right\} }$ to be equal to $\phi_i$
whereas each restriction $\phi|_{X\times\left\{t\right\} }$ is
defined equal to $\phi'_i$ for all $t$ between $\frac{i-1}{s}$ and
$\frac{i}{s}$. Let $\E_{D^m\times I}$ be the $W_\bullet$-bundle
corresponding to $\phi\co U_\bullet\rightarrow W_\bullet$.

Recall that $D^m\times I/S^{m-1}$ denotes the quotient space of
$D^{m}\times I$ where we have identified $(\mathbf{x},t_0)$ with
$(\mathbf{x},t_1)$ for $\mathbf{x}\in S^{m-1}$ and any two
$t_0,t_1\in I$. Let $H\co D^m\times I/S^{m-1}\rightarrow D^m\times
I$ be a map with $H(\mathbf{x},1)\in S^{m-1}\times I \cup
D^m\times\{1\}$ and $H(\mathbf{x},0)=(\mathbf{x},0)$ for all
$\mathbf{x}\in D^m$. Then $H^*\E_{D^m\times I}$ gives a concordance
relative to $q_*\E_{S^{m-1}}$ between the $W_\bullet$-bundle $\E_0$
and another $W_\bullet$-bundle $\E'$ such that $\E'=q_*\E$ for some
$Z_\bullet$-bundle $\E$ with $\E|_{S^{m-1}}=\E_{S^{m-1}}$. Thus we
have shown that the map $ \Con^{f}_{Z_\bullet}(D^m, S^{m-1};
\E_{S^{m-1}}) \rightarrow \Con^{f}_{W_\bullet}(D^m, S^{m-1};
q_*\E_{S^{m-1}})$ is surjective.
\end{proof}

We will now show how Quillen's small object argument, see for
example~\cite[Proposition~7.17]{DwyerSpalinski:95}, can be use to
replace an arbitrary simplicial space $Z_\bullet$ with a simplicial
space $\tilde{Z}_\bullet$ satisfying the topological Kan condition,
Definition~\ref{Def:Kanconditions}. First observe that
$\tilde{Z}_\bullet$ satisfies the topological Kan condition if and
only if all horn filling problems of the forms
\[
\xymatrix{ \Lambda^{n,k}_\bullet \ar@{^{(}->}[d] \ar[r] & \tilde{Z}_\bullet \\
\Delta^n_\bullet \ar@{.>}[ur] } \quad\text{and}\quad\xymatrix{
D^l\times I\times\Lambda^{n,k}_\bullet\cup
D^l\times\{0\}\times\Delta^n_\bullet
 \ar@{^{(}->}[d] \ar[r] & \tilde{Z}_\bullet \\
D^l\times I\times \Delta^n_\bullet \ar@{.>}[ur] }
\]
can be solved. Therefore we define $J$ as to be the union of the
sets of simplicial maps
\[
\{ \Lambda^{n,k}_\bullet\hookrightarrow\Delta^n_\bullet \} \cup \{
D^l\times I\times\Lambda^{n,k}_\bullet\cup
D^l\times\{0\}\times\Delta^n_\bullet \hookrightarrow D^l\times
I\times \Delta^n_\bullet \}
\]
where the integers $k$, $l$ and $n$ vary freely.

A \emph{relative $J$-cell complex} is a map $q\co
Z_\bullet\rightarrow W_\bullet$ that can be constructed as a
transfinite composition of pushouts of elements of $J$,
see~\cite[Definition~10.5.8]{Hirschhorn:03}. A \emph{presentation}
of a relative $J$-cell complex $q\co Z_\bullet\rightarrow W_\bullet$
is a particular choice of how to construct $W_\bullet$ from
$Z_\bullet$, see~\cite[Section~10.6]{Hirschhorn:03}. Moreover, a
\emph{subcomplex} of a presented relative $J$-cell complex $q\co
Z_\bullet\rightarrow W_\bullet$ is a presented relative $J$-cell
complex $q'\co Z_\bullet\rightarrow W'_\bullet$ together with a map
$W'_\bullet\rightarrow W_\bullet$ sending each $J$-cell of
$W'_\bullet$ to a $J$-cell of $W_\bullet$. A relative $J$-cell
complex is \emph{finite} if it can be constructed by attaching only
finitely many $J$-cells, i.e. if the transfinite composition is
actually finite.

\begin{Lem}\label{lem:compact}
Let $C_\bullet$ be a presimplicial space such that each $C_k$ is
compact and there exists an integer $K$ such that $C_k=\emptyset$
for all $k>K$. Let $q\co Z_\bullet\rightarrow W_\bullet$ be a
relative $J$-cell complex. Then for any presimplicial map $f\co
C_\bullet\rightarrow W_\bullet$ there is a finite subcomplex $q'\co
Z_\bullet\rightarrow W'_\bullet$ such that $f$ factors as a
presimplicial map $f'\co C_\bullet\rightarrow W'_\bullet$ followed
by the subcomplex inclusion $W'_\bullet\rightarrow W_\bullet$.
\end{Lem}

\begin{proof}
Fix a presentation of $q\co Z_\bullet\rightarrow W_\bullet$, and
observe that each $q_k\co Z_k\rightarrow W_k$ is a relative cell
complex of topological spaces.
By~\cite[Proposition~10.8.7]{Hirschhorn:03} it follows that $f_k\co
C_k\rightarrow W_k$ intersects only interiors of finitely many
cells. Using that $C_k$ is non-empty in only finitely many
simplicial degrees, we see that the image of $f$ intersects only
interiors of finitely many $J$-cells of the given presentation of
$Z_\bullet\rightarrow W_\bullet$.

To finish the proof it is enough to show that each $J$-cell of
$W_\bullet$ is contained in a finite subcomplex. The argument uses a
transfinite induction over the given presentation of the relative
$J$-cell complex. Observe that the non-degenerate part of the domain
of a $J$-cell satisfies the same assumptions as the presimplicial
space $C_\bullet$. It follows from the first part of this proof that
any attaching map of a $J$-cell intersects only interiors of
finitely many $J$-cells. By the induction hypothesis it follows that
the attaching map of a $J$-cell has image contained in a finite
subcomplex, and we are done.
\end{proof}

\begin{Prop}\label{Prop:relJcellisconequiv}
Let $X$ be a space and $A$ a subspace such that $(X,A)$ has the
homotopy type of a CW pair. Let $q\co Z_\bullet\rightarrow
W_\bullet$ be a relative $J$-cell complex. Then for any
$Z_\bullet$-bundle $\E_A$ over $A$ the induced map
\[
\Con_{Z_\bullet}(X,Z;\E_A)\rightarrow \Con_{W_\bullet}(X,A;q_*\E_A)
\]
is a bijection.
\end{Prop}

\begin{proof}
Using Proposition~\ref{Prop:finite} we reduce to showing that $q$
induces surjections
\[
q_*\co \Con^{f}_{Z_\bullet}(D^m,S^{m-1};\E_{S^{m-1}})\rightarrow
\Con^{f}_{W_\bullet}(D^m,S^{m-1};q_*\E_{S^{m-1}})
\]
for any $m$ and finite $Z_\bullet$-bundle $\E_{S^{m-1}}$ over the
boundary sphere. Let $\E_0$ be a $W_\bullet$-bundle that equals
$q_*\E_{S^{m-1}}$ over $S^{m-1}$, and let $\U$ be the finite and
totally ordered cover associated to $\E_0$. Choose a partition of
unity $\{\psi_\alpha\}$ subordinate to $\U$. Denote by $\phi\co
U_\bullet\rightarrow W_\bullet$ the simplicial map associated to
$\E$.

Since $\U$ is finite there is an $\epsilon>0$ such that
$\U^{\epsilon}$ defined as
$\{\phi_{\alpha}^{-1}\left(\epsilon,1\right]\}$ is a cover of $D^m$.
Between $\U$ and its shrinking $\U^{\epsilon}$ there is a cover
$\{\phi_{\alpha}^{-1}\left[\epsilon,1\right]\}$ consisting of
compact subsets. Associated to these three covers we have the
non-degenerate part of the ordered \v{C}ech complexes and natural
inclusions
\[
N^{\epsilon}_\bullet\rightarrow C_\bullet\rightarrow N_\bullet.
\]
We observe that $C_\bullet$, associated to the compact cover,
satisfies the condition of Lemma~\ref{lem:compact}. Consequently,
the restriction of $\phi$ to $C_\bullet$ factors through some finite
subcomplex $q'\co Z_\bullet\rightarrow W'_\bullet$ of the relative
$J$-cell complex $q\co Z_\bullet\rightarrow W_\bullet$. Hence, the
restriction of $\phi$ to $U_\bullet^{\epsilon}$ factors as a map
$U_\bullet^{\epsilon}\rightarrow W'_\bullet$ followed by the
subcomplex inclusion $W'_\bullet\rightarrow W_\bullet$. Using the
technique from the proof of Corollary~\ref{Cor:confillhorn} we can
now construct a finite $W'_\bullet$-bundle $\E'$ with
$\E'|_{S^{m-1}}=q'_*\E_{S^{m-1}}$ such that the associated
$W_\bullet$-bundle is concordant to $\E_0$ relative to $S^{m-1}$.
Thus $[\E_0]$ is in the image of
\[
\Con^f_{W'_\bullet}(D^m,S^{m-1};q'_*\E_{S^{m-1}})\rightarrow
\Con^f_{W_\bullet}(D^m,S^{m-1};q_*\E_{S^{m-1}}).
\]

Since $q'\co Z_\bullet\rightarrow W'_\bullet$ is a finite relative
$J$-cell complex there is a finite chain of simplicial spaces
\[
Z_\bullet=W_\bullet^{0}\rightarrow W_\bullet^{1}\rightarrow
W_\bullet^{2} \rightarrow \cdots \rightarrow
W_\bullet^{s-1}\rightarrow W_\bullet^{s}= W'_\bullet
\]
such that each $W_\bullet^{i}\rightarrow W_\bullet^{i+1}$ is the
pushout of some element from $J$. Let $q^i$ denote the composition
$Z_\bullet\rightarrow W^i_\bullet$. We claim that for each $i$ the induced map of
concordance classes
\[
\Con^{f}_{W^{i}_\bullet}(D^m,S^{m-1};q^{i}_*\E_{S^{m-1}})
\xrightarrow{\cong}
\Con^{f}_{W^{i+1}_\bullet}(D^m,S^{m-1};q^{i+1}_*\E_{S^{m-1}})
\]
is a bijection. There are two cases to verify. First, if
$W^{i+1}_\bullet$ is the pushout of $D^l\times I\times
\Delta^n_\bullet \hookleftarrow D^l\times
I\times\Lambda^{n,k}_\bullet\cup
D^l\times\{0\}\times\Delta^n_\bullet \rightarrow W^{i}_\bullet$,
then $W^i_\bullet\rightarrow W^{i+1}_\bullet$ is in each simplicial
degree a homotopy equivalence, whence the map above is a bijection
by Theorem~\ref{Thm:ConDegwiseHEQ}. Secondly, if $W^{i+1}$ is the
pushout of $\Delta^n_\bullet\hookleftarrow
\Lambda_\bullet^{n,k}\rightarrow W^i_\bullet$, the the map above is
a bijection by Corollary~\ref{Cor:confillhorn}.

Composing these bijections, we get a bijection
\[
\Con^{f}_{Z_\bullet}(D^m,S^{m-1};\E_{S^{m-1}})
\xrightarrow{\cong}
\Con^{f}_{W'_\bullet}(D^m,S^{m-1};q'_*\E_{S^{m-1}}).
\]
Thus there is a $Z_\bullet$-bundle $\E$ whose restriction to $S^{m-1}$ equals $\E_{S^{m-1}}$
and such that $[q_*\E]=[\E_0]$. This finish the proof.
\end{proof}

\begin{Lem}\label{lem:relJcelliswe}
If $q\co Z_\bullet\rightarrow W_\bullet$ is a relative $J$-cell
complex, then its geometric realization $|Z_\bullet|\rightarrow
|W_\bullet|$ is a weak equivalence.
\end{Lem}

\begin{proof}
Pushouts commutes with geometric realization, and the geometric
realization of any element of $J$ is both a relative CW complex and
a weak equivalence. The conclusion follows.
\end{proof}

By Quillen's small object argument we now get:

\begin{Thm}\label{Thm:topolKanreplacement}
Let $Z_\bullet$ be a simplicial space, let $X$ be a space, and $A$ a
subspace such that $(X,A)$ has the homotopy type of a CW pair. Then
there exists a simplicial space $\tilde{Z}_\bullet$ together with a
simplicial map $q\co Z_\bullet\rightarrow \tilde{Z}_\bullet$ such
that
\begin{itemize}
\item[i)] $\tilde{Z}_\bullet$ satisfies the topological Kan condition,
\item[ii)] for any $Z_\bullet$-bundle $\E_A$ over $A$ the induced map
\[
\Con_{Z_\bullet}(X,A;\E_A)\xrightarrow{q_*}\Con_{\tilde{Z}_\bullet}(X,A;q_*\E_A)
\]
is a bijection, and
\item[iii)] the geometric realization $|q|\co|Z_\bullet|\rightarrow |\tilde{Z}_\bullet|$ is
a weak equivalence.
\end{itemize}
Moreover, if $Z_\bullet$ is a good simplicial space, then
$\tilde{Z}$ can also be taken to be a good simplicial space.
\end{Thm}

\begin{proof}
The set $J$ of simplicial maps permits the small object argument. By
the small object argument, see for
example~\cite[Proposition~10.5.16]{Hirschhorn:03}, it is possible to
construct a simplicial map $q\co
Z_\bullet\rightarrow\tilde{Z}_\bullet$ such that $q$ is a relative
$J$-cell complex and $\tilde{Z}_\bullet\rightarrow *$ has the right
lifting property with respect to any elemeny of $J$. The last
statement immediately implies that $\tilde{Z}_\bullet$ satisfies the
topological Kan condition. By Lemma~\ref{lem:relJcelliswe} the
geometric realization $|q|$ is a weak equivalence. Property iii) is
implied by Proposition~\ref{Prop:relJcellisconequiv}.
\end{proof}

\begin{Rem}
It is natural to ask for a model structure on simplicial spaces such
that the weak equivalences are the maps $Z_\bullet\rightarrow
W_\bullet$ that induce bijections $q\co
\Con_{Z_\bullet}(D^m,S^{m-1};\E_{S^{m-1}})\rightarrow\Con_{W_\bullet}(D^m,S^{m-1};q_*\E_{S^{m-1}})$
for all $m$ and $\E_{S^{m-1}}$, and fibrant objects are the
simplicial spaces satisfying the topological Kan condition. Such a
model structure would not be identical to the realization model
structure of~\cite{RezkSchwedeShipley:01}. We would like more
fibrant objects, compare~\cite[Lemma~8.10]{RezkSchwedeShipley:01}.
However, answering this question is beyond the scope of the present
paper.
\end{Rem}

\section{Proof of Theorem~\ref{thm:main}}\label{sect:proof}

We want to relate concordance classes of $Z_\bullet$-bundles to
homotopy classes of maps into the geometric realization
$|Z_\bullet|$. Intuitively, a $Z_\bullet$-bundle over $X$ is a
simplicial map $\phi\co U_\bullet\rightarrow Z_\bullet$, and
realizing gives a map $|U_\bullet|\rightarrow|Z_\bullet|$. As
explained in~\cite{DuggerIsaksen:04} there is a natural weak
equivalence $X\xleftarrow{\simeq}|U_{\bullet}|$ for any ordered open
cover $\U$. Whenever $X$ has the homotopy type of a CW-complex each
$Z_\bullet$-bundle yields a well-defined homotopy class of maps from
$X$ into $|Z_\bullet|$. Actually, we will not use Dugger and
Isaksens weak equivalence, but rather we compare homotopy classes
and concordance classes by more direct means.

Given a topological space $Y$ there are several ways to functorially
associate a simplicial space. We have the \emph{constant simplicial
space} $C_\bullet(Y)$ defined in each simplicial degree as
$C_k(Y)=Y$ and with all face and degeneracy maps equal to the
identity map on $Y$. On the other hand we have the \emph{continuous
singular simplicial space} $S_\bullet(Y)$ whose space of
$k$-simplices is the mapping space $\Top(\Delta^k,Y)$. Observe that
$C_\bullet(-)$ is the left adjoint to the space of $0$-simplices
functor $Z_\bullet\mapsto Z_0$ and that $S_\bullet(-)$ is the right
adjoint to geometric realization of simplicial spaces. Since we have
a homeomorphism $Y\xrightarrow{\cong} S_0(Y)$ there is a natural map
$C_\bullet(Y)\rightarrow S_\bullet(Y)$. Inspecting the definitions,
we see that this natural map assigns to a point $y$ in $C_k(Y)=Y$
the constant map sending all of $\Delta^k$ to $y\in Y$. We may also
consider the discrete topology on the set of maps
$\Delta^k\rightarrow Y$. This defines the \emph{singular simplicial
set} $\Sing Y$. There is a natural simplicial map $\kappa\co\Sing
Y\rightarrow S_\bullet(Y)$ given as the identity on the underlying
sets.

\begin{Lem}\label{Lem:SYtoYwe}
The simplicial spaces $C_\bullet(Y)$, $S_\bullet(Y)$, and $\Sing Y$ are good,
the geometric realization $|C_\bullet(Y)|$ is naturally homeomorphic to $Y$, and
the natural maps induces weak equivalences
\[
Y\cong |C_\bullet(Y)|\xrightarrow{\simeq} |S_\bullet(Y)| \xleftarrow{\simeq} |\Sing Y|.
\]
\end{Lem}

\begin{proof}
Clearly $C_\bullet(Y)$ is a good simplicial space, i.e. all degeneracies are closed cofibrations, and
we see that the geometric realization of $C_\bullet(Y)$ is
naturally homeomorphic to $Y$. Furthermore, it is not
difficult to show explicitly that the map
$\Top(\Delta^{k-1},Y)\rightarrow \Top(\Delta^k,Y)$, induced by
$\sigma_i\co \Delta^k\rightarrow\Delta^{k-1}$, is the inclusion of a
deformation retract. Consequently, also $S_\bullet(Y)$ is good. The singular simplicial
set is in each simplicial degree discrete, hence $\Sing Y$ is good.

Let $q\co Z_\bullet \rightarrow W_\bullet$ be simplicial map between
good simplicial spaces. It is a classical result that the geometric
realization $|q|$ is a weak equivalence if each $q_k\co
Z_k\rightarrow W_k$ is a weak equivalence, see for
example~\cite[Section~2.2]{Madsen:94}. This applies to the
simplicial map $C_\bullet(Y)\rightarrow S_\bullet(Y)$, i.e. we have
a weak equivalence
\[
Y\cong |C_\bullet(Y)|\xrightarrow{\simeq} |S_\bullet(Y)|.
\]

It is well-known that the counit of the adjunction between $|-|$ and
$\Sing$ is a natural weak equivelence $|\Sing Y|\xrightarrow{\simeq}
Y$. Now observe that the counit of the adjunction between geometric
realization and the functor $S_\bullet(-)$ is a retraction to the
weak equivelence $|C_\bullet(Y)|\xrightarrow{\simeq}|S_\bullet(Y)|$.
The result now follows by commutativity of the diagram
\[
\xymatrix{ |S_\bullet(Y)| \ar[dr]^{\simeq} && |\Sing Y|
\ar[ll]_{\kappa} \ar[dl]_{\simeq} \\ & Y.}
\]
\end{proof}

Let $A$ be a subspace of $X$, and fix a continuous map $g\co
A\rightarrow Y$. Denote by $[X,A;Y,g]$ the set of maps $f\co
X\rightarrow Y$ whose restriction to $A$ equals the fixed map $f$
modulo homotopy relative to $A$, i.e. two such maps $f_0$ and $f_1$
represent the same class if there is a homotpy between them that is
constant equal to $g$ on $A$.

Let $\phi\co U_\bullet\rightarrow C_\bullet(Y)$ be the simplicial
map of a $C_\bullet(Y)$-bundle $\E$ over $X$ subordinate to some
ordered open cover $\U$. At the level of $0$-simplices $\phi$ maps
each open set $U_\alpha\in\U$ into $Y$. Moreover, these maps agree
on each intersection $U_\alpha\cap U_\beta$. Hence $\phi$ comes from
a continuous map $X\rightarrow Y$, and we call this map the
\emph{underlying map} of the $C_\bullet(Y)$-bundle. We now have the
following important, but easy observation:

\begin{Prop}\label{prop:conhtpy}
Let $X$ and $Y$ be topological spaces, and let $A$ be a subspace of
$X$. Let $g\co A\rightarrow Y$ be the underlying map of a
$C_\bullet(Y)$-bundle $\E_A$ over $A$. Then the underling map of
$C_\bullet(Y)$-bundles defines a natural bijection
\[
\Con_{C_\bullet(Y)}(X,A;\E_A)\xrightarrow{\cong}[X,A;Y,g].
\]
\end{Prop}

\begin{proof}
We compare the notion of a map $X\rightarrow Y$ to the notion of a
$C_{\bullet}(Y)$-bundle over $X$. The only additional piece of data
contained in the latter definition is the choice of an ordered open
cover $\U$ over $X$. Clearly any two ways of making this choice is
equivalent up to concordance. The result follows.
\end{proof}

Let us compare $C_\bullet(Y)$-bundles and $S_\bullet(Y)$-bundles.

\begin{Prop}\label{prop:conCS}
Let $A\hookrightarrow X$ be a closed cofibration such that $(X,A)$
has the homotopy type of a CW pair. Let $Y$ be a topological space,
and let $\E_A$ be any $C_\bullet(Y)$-bundle over $A$. Then the
natural map $i\co C_\bullet(Y)\rightarrow S_\bullet(Y)$ induces a
bijection
\[
\Con_{C_\bullet(Y)}(X,A;\E_A)\xrightarrow{\cong}\Con_{S_\bullet(Y)}(X,A;i_*\E_A).
\]
\end{Prop}

\begin{proof}
This follows immediately from Theorem~\ref{Thm:ConDegwiseHEQ} since
the simplicial map $i$ in each simplicial degree $k$ is a homotopy
equivalence $Y=C_k(Y)\rightarrow S_k(Y)=\Top(\Delta^k,Y)$.
\end{proof}

Let us also compare $S_\bullet(-)$-bundles for weakly equivalent
spaces $Y$ and $Y'$.

\begin{Cor}\label{cor:conSSwe}
Let $A\hookrightarrow X$ be a closed cofibration such that $A$ has
the homotopy type of a CW complex and $(X,A)$ has the homotopy type
of a CW pair. Let $f\co Y\rightarrow Y'$ be a weak equivalence, and
let $\E_A$ be any $S_\bullet(Y)$-bundle over $A$. Then $f$ induces a
bijection
\[
\Con_{S_\bullet(Y)}(X,A;\E_A)\xrightarrow{\cong}\Con_{S_\bullet(Y')}(X,A;f_*\E_A).
\]
\end{Cor}

\begin{proof}
By Proposition~\ref{prop:conCS} there is a $C_\bullet(Y)$-bundle
$\E'_A$ over $A$ and a concordance $\E_{A\times I}$ between the
$S_\bullet(Y)$-bundles $i_*\E'_A$ and $\E_A$. By
Lemma~\ref{lem:CAZgroupoid} and
Proposition~\ref{Prop:changeofcollar} gluing of $\E_{A\times I}$
gives bijections. Hence, we have a diagram
\[
\xymatrix{
\Con_{C_\bullet(Y)}(X,A;\E'_A) \ar[d]_{\text{\ding{172}}} \ar[r]^{\text{\ding{173}}} &
\Con_{S_\bullet(Y)}(X,A;i_*\E'_A) \ar[d] \ar[r]^{\text{\ding{175}}} &
\Con_{S_\bullet(Y)}(X,A;\E_A) \ar[d]^{\text{\ding{177}}}\\
\Con_{C_\bullet(Y')}(X,A;f_*\E'_A) \ar[r]^{\text{\ding{174}}} &
\Con_{S_\bullet(Y')}(X,A;f_*i_*\E'_A) \ar[r]^{\text{\ding{176}}} &
\Con_{S_\bullet(Y')}(X,A;f_*\E_A) }
\]
where \ding{175} and \ding{176} are bijections. Furthermore,
\ding{172} is a bijection by Proposition~\ref{prop:conhtpy}, and
\ding{173} and \ding{174} are bijections by
Proposition~\ref{prop:conCS}. It follows that \ding{177} is a
bijection.
\end{proof}

Let $Z_\bullet$ be a simplicial space. Since the diagram $C_\bullet(-)\rightarrow S_\bullet(-)
\leftarrow \Sing(-)$ is natural, we can insert $Z_\bullet$ and produce a diagram
of bisimplicial spaces
\[
C_\bullet(Z_\bullet)\rightarrow S_\bullet(Z_\bullet)
\leftarrow \Sing(Z_\bullet).
\]
Applying the diagonal we recover from $C_\bullet(Z_\bullet)$ the
simplicial space $Z_\bullet$. Define the functors $T$ and $D$ by
sending $Z_\bullet$ to $\diag S_\bullet(Z_\bullet)$ and $\diag
\Sing(Z_\bullet)$ respectively. The underlying sets of
$D(Z_\bullet)$ and $T(Z_\bullet)$ are identical. In both cases a
$k$-simplex is represented by a continuous map $\Delta^k\rightarrow
Z_k$. The difference is the topology we place on this set. In the
case of $D(Z_\bullet)$ we take the discrete topology. In the case of
$T(Z_\bullet)$ we take the compact-open topology on each of the
mapping spaces. Let $\kappa\co D(Z_\bullet)\rightarrow T(Z_\bullet)$
denote the natural map which is the identity on the underlying sets.

\begin{Lem}\label{Lem:DTZagree}
Assume that $Z_\bullet$ is a good simplicial space.
The geometric realization of $Z_\bullet\rightarrow T(Z_\bullet)\leftarrow D(Z_\bullet)$
gives a diagram of natural weak equivalences
\[
|Z_\bullet|\xrightarrow{\simeq} |T(Z_\bullet)|\xleftarrow{\simeq} |D(Z_\bullet)|.
\]
\end{Lem}

\begin{proof}
Suppose that $W_{\bullet,\bullet}$ is a bisimplicial space. We can
take the geometric realization in two steps; first by defining a
simplicial space $[k]\mapsto |W_{k,\bullet}|$ and then realizing
again to get $\left|[k]\mapsto |W_{k,\bullet}|\right|$, or we can
take the geometric realization of the diagonal, $|\diag
W_{\bullet,\bullet}|$. It is a classical result that these two
construction produce naturally homeomorphic spaces.

By Lemma~\ref{Lem:SYtoYwe} we get for each $k$ weak equivalences
\[
Z_k\xrightarrow{\simeq} |S_\bullet(Z_k)| \xleftarrow{\simeq} |\Sing Z_k|.
\]
The simplicial space $Z_\bullet$ is good by assumption, and it is
not hard to see that $[k]\mapsto|\Sing Z_k|$ is good. Observing that
each degeneracy $s_i\co Z_{k-1}\rightarrow Z_k$ is a closed
cofibration, we use~\cite[Lemma~4]{Strom:72} to conclude that the
induced map $\Top(\Delta^{l},Z_{k-1})\rightarrow
\Top(\Delta^{l},Z_k)$ also is a closed cofibration. It follows that
$[k]\mapsto |S_\bullet(Z_k)|$ is a good simplicial space.

The statement of the lemma follows since $|Z_\bullet|\rightarrow
|T(Z_\bullet)|\leftarrow |D(Z_\bullet)|$ is homeomorphic to the
geometric realization of simplicial maps which are weak equivalences
in each simplicial degree.
\end{proof}

Given a simplicial space $Z_\bullet$, let $\eta\co
Z_\bullet\rightarrow S_\bullet(|Z_\bullet|)$ denote the unit of the
adjunction between geometric realization and the continuous singular
simplicial space. We use $\eta$ to relate $Z_\bullet$-bundles to
homotopy classes of maps into $|Z_\bullet|$. To be precise, if $\E$
is a $Z_\bullet$-bundle over a CW-complex $X$ represented by a
simplicial map $U_\bullet\xrightarrow{\phi}Z_\bullet$, then
composing with $\eta$ yields $U_\bullet\rightarrow
S_\bullet(|Z_\bullet|)$, and this simplicial map corresponds to a
well-defined homotopy class of maps $X\rightarrow|Z_\bullet|$ by
Proposition~\ref{prop:conhtpy} and Proposition~\ref{prop:conCS}. Our
aim is to prove that that $\eta$, under reasonable assumptions,
induces a bijection between concordance classes. To handle a
technical point in this proof we need:

\begin{Thm}[(Whitney)]\label{thm:Whitney}
Every open subset of $S^{n-1}$ is a CW-complex.
\end{Thm}

\begin{proof}
Every proper subset $U$ of $S^{n-1}$ embeds as an open subset of
$\R^{n-1}$. Consider meshes $\mathscr{M}_s$ in $\R^{n-1}$ of cubes
with side length $2^{-s}$, such that $\mathscr{M}_{s+1}$ is a
subdivision of $\mathscr{M}_s$.
Following~\cite[Section~8]{Whitney:34}, we can write $U$ as a union
$\bigcup_s K_s$ such that $K_s$ consists of cubes in
$\mathscr{M}_s$, the interiors of the $K_s$'s are mutually disjoint,
and $K_s$ meets only $K_{s-1}$ and $K_{s+1}$. This is a CW-structure
on $U$.
\end{proof}

\begin{Lem}\label{Lem:surjcon}
Assume that $Z_{\bullet}$ satisfies the topological Kan condition
and $\U$ is a good ordered open cover of $S^{n-1}$. For every
simplicial map $\phi\co U_\bullet\rightarrow T(Z_\bullet)$ there
exists a simplicial map $\Phi\co U_\bullet\times I\rightarrow
T(Z_\bullet)$ such that $\phi=\Phi(-,1)$ and $\Phi(-,0)$ factors as
$U_{\bullet}\rightarrow D(Z_{\bullet})\xrightarrow{\kappa}
T(Z_\bullet)$.
\end{Lem}

\begin{proof}
Let $N_\bullet$ be the non-degenerate part of the ordered \v{C}ech
complex $U_\bullet$ associated to the good ordered open cover $\U$.
It is enough to define $\Phi$ as a presimplicial map
$N_\bullet\times I\rightarrow T(Z_\bullet)$. Inspecting the
definition of $T(Z_\bullet)$, we see that maps $N_\bullet\rightarrow
T(Z_\bullet)$ can be described by continuous maps $\phi_k\co
N_k\times\Delta^k\rightarrow Z_k$ such that the following
\emph{presimplicial coherence diagrams} commute for all face maps
$d_i$:
\[
\xymatrix@C=1.5cm{ & N_k\times \Delta^k \ar[r]^-{\phi_k} & Z_k
\ar[dd]^{d_i}
\\ N_k\times \Delta^{k-1} \ar[ur]^{\id\times\delta_i}
\ar[dr]_{d_i\times \id} && \\ & N_{k-1}\times\Delta^{k-1}
\ar[r]^-{\phi_{k-1}} & Z_{k-1}.}
\]
Observe that a map $N_\bullet\rightarrow T(Z_\bullet)$ factors
through $\kappa\co D(Z_\bullet)\rightarrow T(Z_\bullet)$ if and only
if for each $k$ and $\mathbf{t}\in \Delta^k$ the function sending
$x\in N_k$ to $\phi_k(x,\mathbf{t})$ is a locally constant on $N_k$.

By definition the space of non-degenerate $k$-simplices is
$N_k=\coprod_{\alpha_0<\cdots<\alpha_k} U_{\alpha_0\cdots\alpha_k}$.
Since $\U$ is a good cover, each connected component of $N_k$ is
contractible. Hence, there exists local contractions $H^k\co
N_k\times I\rightarrow N_k$, i.e. homotopies such that for $x\in
N_k$ we have $H^k(x,1)=x$, whereas $H^k(x,0)$ is locally constant.

Let $\phi_k\co N_k\times\Delta^k\rightarrow Z_k$ be the collection
of maps that represents the given simplicial map $\phi\co
U_\bullet\rightarrow T(Z_\bullet)$. By induction on $k$ we will
construct maps $\Phi_k\co N_k\times I\times \Delta^k\rightarrow Z_k$
such that
\begin{itemize}
\item[i)] $\Phi_{k}(x,1,\mathbf{t})=\phi_{k}(x,\mathbf{t})$ for all
$x\in N_k$ and $\mathbf{t}\in \Delta^k$,
\item[ii)] for $\mathbf{t}\in \Delta^k$ fixed, the expression $\Phi_{k}(x,s,\mathbf{t})$ is locally constant as
a function $N_k\times \left[ 0,\frac{1}{k+2}\right]\rightarrow Z_k$,
and
\item[iii)] all presimplicial coherence diagrams for $\Phi$ commute.
\end{itemize}
We begin by defining $\Phi_0$. For $x\in N_0$, $s\in I$ and
$1\in\Delta^0$ we define
\[
\Phi_0(x,s,1)=\begin{cases}\phi_0(H^0(x,2s-1),1)&\text{if
$s\geq\frac{1}{2}$}\\ \phi_0(H^0(x,0),1)&\text{if
$s\leq\frac{1}{2}$.}\end{cases}
\]
Next, we construct $\Phi_k$ from $\Phi_i$, $i<k$. This is done in
four steps:

\emph{First step.} In this step we will extend $\phi_k\co
N_k\times\{1\}\times \Delta^k\rightarrow Z_k$ to a map
\[
\Phi_k\co \big( N_k\times\left[ \frac{1}{k+1},1\right]\times
\partial\Delta^k\big)\cup\big(N_k\times\{1\}\times
\Delta^k\big)\rightarrow Z_k.
\]
We will do this by induction over the dimension of proper faces
$\sigma$ of $\Delta^k$. Given a face $\sigma$, let $I$ be the subset
of $[k]$ such that $n\in I$ if and only if $\sigma$ is contained in
the $n$'th $(k-1)$-face of $\Delta^k$. We can, by the induction
hypothesis, construct a diagram
\[
\xymatrix{ N_k\times \left( \left[ \frac{1}{k+1},1\right]\times
\partial\sigma \cup \{1\}\times \sigma
\right) \ar[r] \ar[d] & Z_k \ar[d]^{c^I_k}\\
N_k\times \left[ \frac{1}{k+1},1\right]\times \sigma \ar[r] &
C\Lambda^I_k(Z_\bullet), }
\]
where the lower map comes from the $\Phi_i$, $i<k$. The vertical map
on the left side is a cellular inclusion between CW-complexes, and a
homotopy equivalence. The topological Kan condition says that
$c^I_k$ is a Serre fibration. Hence, a lift exists, and we use this
to define $\Phi_k$ on $N_k\times \left[ \frac{1}{k+1},1\right]\times
\sigma$.

\emph{Second step.} Let
\[
r\co \left[ \frac{1}{k+1},1\right]\times \Delta^k \rightarrow \big(
\left[ \frac{1}{k+1},1\right]\times
\partial\Delta^k\big) \cup \big(\{1\}\times \Delta^k \big)
\]
be a retraction. For $x\in N_k$, $\frac{1}{k+1}\leq s\leq 1$,
$\mathbf{t}\in\Delta^k$, we define
\[
\Phi_k(x,s,\mathbf{t})=\Phi_k(r(x,s,\mathbf{t})),
\]
where the $\Phi_k$ on the right side is defined by the previous
step.

\emph{Third step.} Recall that $H^k$ denoted a local contraction of
$N_k$. Choose a function $f\co
\left[\frac{1}{k+2},\frac{1}{k+1}\right]\rightarrow I$ such that
$f(\frac{1}{k+2})=0$ and $f(\frac{1}{k+1})=1$. For $x\in N_k$,
$\frac{1}{k+2}\leq s\leq 1$, $\mathbf{t}\in\Delta^k$, we now define
\[
\Phi_k(x,s,\mathbf{t})=\begin{cases} \Phi_k(x,s,\mathbf{t})
&\text{if $\frac{1}{k+1}\leq s$,}\\
\Phi_k(H^k(x,f(s)),\frac{1}{k+1},\mathbf{t}) &\text{if
$\frac{1}{k+2}\leq s\leq\frac{1}{k+1}$.}\end{cases}
\]
Again, the $\Phi_k$ on the right comes from the previous step. At
this stage, we should verify that the presimplicial coherence
diagrams commute for $s\in\left[\frac{1}{k+2},\frac{1}{k+1}\right]$,
i.e. we need to check the equality
\[
d_i\Phi_k(x,s,\delta_i\mathbf{t})=\Phi_{k-1}(d_ix,s,\mathbf{t})
\]
for all face maps $d_i$, $x\in N_k$, $s$ as above and
$t\in\Delta^{k-1}$. This is an easy exercise using that equality
holds for $s=\frac{1}{k+1}$, and that when fixing $\mathbf{t}$ the
function $\Phi_{k-1}(y,s,\mathbf{t})$ is locally constant on
$N_{k-1}\times\left[ 0,\frac{1}{k+1}\right]$.

\emph{Final step.} The steps above define $\Phi_k$ on
$N_k\times\left[\frac{1}{k+2},1\right]\times\Delta^k$. Now extend
the domain to all of $N_k\times I\times \Delta^k$ by
\[
\Phi_k(x,s,\mathbf{t})=\Phi_k(x,\frac{1}{k+2},\mathbf{t}),\quad\text{for
$s\leq\frac{1}{k+2}$.}
\]
This completes the construction of all $\Phi_k$s.

Since $\Phi_k(x,0,\mathbf{t})$ is locally constant for fixed
$\mathbf{t}\in\Delta^k$, it follows that $\Phi(-,0)\co
N_\bullet\rightarrow T(Z_\bullet)$ factors through $D(Z_\bullet)$.
This completes the proof of the lemma.
\end{proof}

\begin{Lem}\label{Lem:corelemma}
Let $Z_\bullet$ be a good simplicial space satisfying the
topological Kan condition. Let $\U$ be a good ordered open cover of
$D^n$ such that the restriction to $S^{n-1}$ also is good. Let
$U_\bullet$ and $U^{\partial}_\bullet$ denote the ordered \v{C}ech
complex associated to $\U$ and its restriction to $S^{n-1}$
respectively. Suppose that the pair
$(|U_\bullet|,|U^\partial_\bullet|)$ is a relative CW-complex.  Let
$\phi^\partial$ be a simplicial map $U^\partial_\bullet\rightarrow
D(Z_\bullet)$, and let $\tilde{\phi}$ be a continuous map
$|U_\bullet|\rightarrow |T(Z_\bullet)|$ such that the diagram
\[
\xymatrix{ U^\partial_\bullet \ar[d] \ar[r]^{\phi^\partial} &
D(Z_\bullet) \ar[d]^{\eta\kappa} \\ U_\bullet \ar[r]^-{\phi} &
S_\bullet(|T(Z_\bullet)|)}
\]
commutes. Here $\phi$ is the adjoint of $\tilde{\phi}$. Then there
exists a simplicial map $\phi'\co U_\bullet\rightarrow D(Z_\bullet)$
extending $\phi^{\partial}$ such that the composition
$|U_\bullet|\xrightarrow{|\phi'|}|D(Z_\bullet)|\xrightarrow{|\kappa|}|T(Z_\bullet)|$
is homotopic to $\tilde{\phi}$ relative to $|U^\partial_\bullet|$.
\end{Lem}

\begin{proof}
Since $|\kappa|\co|D(Z_\bullet)|\rightarrow |T(Z_\bullet)|$ is a
weak equivalence it follows by classical techniques,
see~\cite[Theorem~7.6.22]{Spanier:91}, that there exists a map
$\tilde{\phi_1}\co |U_\bullet|\rightarrow |D(Z_\bullet)|$ such that
the composition
\[
|U_\bullet|\xrightarrow{\tilde{\phi_1}}|D(Z_\bullet)|\xrightarrow{|\kappa|}|T(Z_\bullet)|
\]
is homotopic to $\tilde{\phi}$ relative to $|U^\partial_\bullet|$.

A key idea in what follows is to replace the \v{C}ech complexes
$U_\bullet$ and $U^\partial_\bullet$ by their simplicial sets of
path components $V_\bullet$ and $V^\partial_\bullet$. To be precise,
we define the $k$-simplices of $V_\bullet$ to be the set of path
components of $U_k$, i.e. we have $V_k=\pi_0 U_k$. Because $\U$ is
good the natural quotient map $q\co U_\bullet\rightarrow V_\bullet$
is a weak equivalence in each simplicial degree. Consequently, we
get a weak equivalence $|U_\bullet|\xrightarrow{\simeq}
|V_\bullet|$. Similarly for $V^\partial_\bullet$.

Inside the mapping cylinder $|U_\bullet|\times I\cup_{|q|}
|V_\bullet|$ we find $|U_\bullet|\times\{0\}\cup_0
|U^\partial_\bullet|\times I\cup_{|q^\partial|}
|V^\partial_\bullet|$ lying as a deformation retract. The
restriction of $\tilde{\phi_1}$ to $|U^\partial_\bullet|$ factors
through $|V_\bullet^\partial|$ since it is the realization of the
simplicial map $\phi^{\partial}\co U^\partial_\bullet\rightarrow
D(Z_\bullet)$ and $D(Z_\bullet)$ is a simplicial set. It follows
that $\tilde{\phi_1}$ extends to a map
\[
|U_\bullet|\times\{0\}\cup_0 |U^\partial_\bullet|\times
I\cup_{|q^\partial|} |V^\partial_\bullet|\rightarrow |D(Z_\bullet)|
\]
being constant in the $I$-coordinate. Using the deformation
retraction, we extend to a map from the mapping cylinder of $q$,
\[
|U_\bullet|\times I\cup_{|q|} |V_\bullet|\rightarrow |D(Z_\bullet)|.
\]
Consequently we get a continuous map $\tilde{\phi_2}\co
|V_\bullet|\rightarrow |D(Z_\bullet)|$ whose restriction to
$|V^\partial_\bullet|$ agrees with the realization of
$\phi^\partial$ and such that there is a homotopy relative to
$|U_\bullet^\partial|$ between $\tilde{\phi_1}$ and the composition
\[
|U_\bullet|\xrightarrow{|q|}|V_\bullet|\xrightarrow{\tilde{\phi_2}}
|D(Z_\bullet)|.
\]

The adjoint of $\tilde{\phi_2}$ fits into the following solid
commutative diagram of simplicial sets,
\[
\xymatrix{ V^\partial_\bullet \ar[r] \ar[d] & D(Z_\bullet)
\ar[d]^{\eta}\\ V_\bullet \ar[r]^-{\phi_2} \ar@{.>}[ur] &
\Sing|D(Z_\bullet)|.}
\]
Since $Z_\bullet$ satisfies the topological Kan condition, we know
that $D(Z_\bullet)$ is a Kan complex by Theorem~\ref{th:Kan}. It
follows that $\eta$ is a weak equivalence between fibrant and
cofibrant objects. It is well-known, see for
example~\cite[Lemma~4.6]{Jardine:04}, that under these circumstances
the dotted map exists such that the upper triangle commutes, while
the lower triangle commutes up to homotopy relative to
$V_\bullet^\partial$. This implies that there is a simplicial map
$\phi_3\co V_\bullet\rightarrow D(Z_\bullet)$ whose restriction to
$V^\partial_\bullet$ is $\phi^\partial$ and such that there is a
homotopy relative to $|U_\bullet^\partial|$ between
$\tilde{\phi_2}$. This concludes the proof.
\end{proof}

\begin{Lem}\label{Lem:cleverrefinement}
Let $\U'$ be any ordered open cover of $D^n$. Then there exists a
refinement $\U$ of $\U'$ together with a carrier function, such that
$\U$ and its restriction to $S^{n-1}$ both are good covers and the
pair $(|U_\bullet|,|U^\partial_\bullet|)$ is a relative CW-complex.
Here $U_\bullet$ and $U^\partial_\bullet$ are the ordered \v{C}ech
complexes associated to $\U$ and the restriction $\U^\partial$ of
$\U$ to $S^{n-1}$.
\end{Lem}

\begin{proof}
To prove the lemma we introduce \emph{convex fragments of spherical
shells}. We call a subset $V$ of $S^{n-1}$ \emph{convex} if
\begin{itemize}
\item[i)] $V$ contains no pair of antipodal point, and
\item[ii)] whenever $\mathbf{x}$ and $\mathbf{y}$ lies in $V$ then the
shortest segment of the great arc between $\mathbf{x}$ and
$\mathbf{y}$ is completely contained in $V$.
\end{itemize}
We call a subset $U$ of $D^n$ a \emph{convex fragment of a spherical
shell} if there is a convex subset $V$ of $S^{n-1}$ and an interval
$J\subseteq(0,1]$ such that
\begin{itemize}
\item[i)] the set $U$ does not contain the center of $D^n$, and
\item[ii)] a point $\mathbf{x}$ in
$D^n$ is contained in $U$ if and only if there exists $\mathbf{y}\in
V$ and $r\in J$ such that $\mathbf{x}=r\mathbf{y}$.
\end{itemize}
A subset $U$ of $D^n$ is a ball if there exists an interval
$J\subseteq[0,1]$ such that $0\in J$ and $\mathbf{x}$ lies in $U$ if
and only if $\|x\|\in J$. Observe that the intersection between a
convex fragment of a spherical shell and either a ball or another
convex fragment of a spherical shell is again a convex fragment of a
spherical shell.

Given the ordered open cover $\U'$ we let $\U$ be an ordered open
cover such that
\begin{itemize}
\item[i)] each $U_\alpha$ in $\U$ is
either an open ball of radius less than $1$ or a open convex
fragment of a spherical shell, and
\item[ii)] each $U_\alpha$ in $\U$ is contained in some $U'_\beta$ in
$\U'$.
\end{itemize}
The existence of such a cover $\U$ follows from the fact that the
collection of all open balls and open convex fragments of spherical
shells is a basis for the topology on $D^n$. Choose any carrier
function $c\co \I\rightarrow \I'$ between the corresponding indexing
sets.

Any non-empty finite intersection $U_{\alpha_0\ldots\alpha_k}$ of
subsets in $\U$ is contractible since it is either a ball or a
convex fragment of a spherical shell. Consequently $\U$ is a good
cover. Moreover, the restriction of $\U$ to $S^{n-1}$ is a cover
$\U^\partial$ consisting of open convex subsets. Hence $\U^\partial$
is also good. To check that $(|U_\bullet|,|U_\bullet^\partial|)$ is
a relative CW complex it is enough to see that for each finite
intersection the pair
$(U_{\alpha_0\ldots\alpha_k},U_{\alpha_0\ldots\alpha_k}\cap
S^{n-1})$ is a relative CW complex. Such a pair is either an open
ball relative to the empty set or an open convex fragment of a
spherical shell relative to an open convex subset of $S^{n-1}$. In
both cases an argument using the technique of
Theorem~\ref{thm:Whitney} shows that we indeed have a relative CW
complex. This proves the lemma.
\end{proof}

\begin{Prop}\label{Prop:Tsurj}
Assume that $Z_{\bullet}$ is good and satisfies the topological Kan
condition. Let $\E_{S^{n-1}}$ be any $T(Z_\bullet)$-bundle over
$S^{n-1}$. Then the map
\[
\eta_*\co\Con_{T(Z_{\bullet})}(D^n,S^{n-1};\E_{S^{n-1}})\rightarrow
\Con_{S_\bullet(|T(Z_{\bullet})|)}(D^n,S^{n-1};\eta_*\E_{S^{n-1}})
\]
is surjective.
\end{Prop}

\begin{proof}
Let $\E'$ be any $S_\bullet(|T(Z_\bullet)|)$-bundle with
$\E'|_{S^{n-1}}=\eta_*\E_{S^{n-1}}$, and let $\U'$ be the ordered
open cover associated to $\E'$. Use Lemma~\ref{Lem:cleverrefinement}
to choose a refinement $\U$ of $\U'$. The carrier function
$c\co\I\rightarrow \I'$ between the corresponding indexing sets
induces simplicial maps $i\co U_\bullet\rightarrow U'_\bullet$ and
$i_\partial\co U^\partial_\bullet\rightarrow {U'}_\bullet^\partial$
between the ordered \v{C}ech complexes of $\U$ and $\U'$ and their
restriction to $S^{n-1}$ respectively. Using $i$ and $i_\partial$ we
change the associated cover and get a new
$S_{\bullet}(|T(Z_\bullet)|)$-bundle $i^*\E'$ and a new
$T(Z_\bullet)$-bundle $i^*_\partial\E_{S^{n-1}}$.

We claim that there exists a concordance of
$S_\bullet(|T(Z_\bullet)|)$-bundles between $\E'$ and $i^*\E'$ whose
restriction to $S^{n-1}$ lifts to a concordance of
$T(Z_\bullet)$-bundles between $\E_{S^{n-1}}$ and
$i^*_\partial\E_{S^{n-1}}$. Let $\phi'\co U'_\bullet\rightarrow
S_\bullet(|T(Z_\bullet)|)$ be the simplicial map corresponding to
$\E'$. By definition the composition $\phi'i$ corresponds to
$i^*\E'$. To see the existence of these concordances it is
sufficient to construct a new ordered open cover $\U''$ together
with a simplicial map $\phi''\co U''_\bullet\rightarrow Z_\bullet$
such that $\U$ and $\U'$ are subcovers and $\phi'i$ and $\phi'$ are
the restrictions of $\phi''$. Let the indexing set $\I''$ of $\U''$
be the disjoint union $\I\amalg \I'$ and give $\I''$ the smallest
partial ordering such that the inclusions of $\I$ and $\I'$ are
order preserving and $\alpha>\beta$ whenever $c(\alpha)=\beta$. We
define $\U''=\U\cup\U'$. Observe that $c$ extends to a carrier
function $c''\co\I''\rightarrow \I'$ inducing a map $r\co
U''_\bullet\rightarrow U'_\bullet$ of the corresponding ordered
\v{C}ech complexes such that $r$ is a retraction of
$U'_\bullet\rightarrow U''_\bullet$. Define $\phi''=\phi'r$. This
proves the claim.

As a consequence of this claim it follows that there is a
commutative diagram
\[
\xymatrix{ \Con_{T(Z_\bullet)}(D^n,S^{n-1};i_\partial^*\E_{S^{n-1}})
\ar[r]^-{\text{\ding{173}}} \ar[d]^{\cong} &
\Con_{S_\bullet(|T(Z_\bullet)|)}(D^n,S^{n-1};\eta_*i_\partial^*\E_{S^{n-1}})
\ar[d]^{\cong} \\
\Con_{T(Z_\bullet)}(D^n,S^{n-1};\E_{S^{n-1}})
\ar[r]^-{\text{\ding{172}}} &
\Con_{S_\bullet(|T(Z_\bullet)|)}(D^n,S^{n-1};\eta_*\E_{S^{n-1}}), }
\]
where the vertical maps are gluing of concordances over $S^{n-1}$.
These vertical maps are bijections by Lemma~\ref{lem:CAZgroupoid}
and Proposition~\ref{Prop:changeofcollar}. Moreover, the class
represented by $\E'$ corresponds to the class represented by
$i^*\E'$. So it is enough to prove that $i^*\E'$ lies in the image
of \ding{173}.

Since $\U^\partial$ is a good cover, Lemma~\ref{Lem:surjcon} applies
to the simplicial map $\phi_{S^{n-1}} i_\partial\co
U^\partial_\bullet\rightarrow T(Z_\bullet)$ associated to
$i_\partial^*\E_{S^{n-1}}$. Thus we get a $D(Z_\bullet)$-bundle
$\E^\partial$ over $S^{n-1}$ together with simplicial map $\Phi\co
U^\partial_\bullet\times I\rightarrow T(Z_\bullet)$ such that
$\Phi(-,1)=\phi_{S^{n-1}} i_\partial$ and
$\Phi(-,0)=\kappa\phi^\partial$, where $\phi^{\partial}\co
U^\partial_\bullet\rightarrow D(Z_\bullet)$ is associated to
$\E^{\partial}$.

Let $\tilde{\phi}'\co|U_\bullet|\rightarrow |T(Z_\bullet)|$ denote
the adjoint of the simplicial map $\phi'i\co U_\bullet\rightarrow
S_\bullet(|T(Z_\bullet)|)$ associated to $i^*\E'$. We see that the
geometric realization of $\Phi$ match up with $\tilde{\phi}'$ and
produce a continuous map
\[
|U^\partial_\bullet|\times I\cup |U_\bullet|\times\{1\}\rightarrow
|T(Z_\bullet)|.
\]
Since $(|U_\bullet|,|U_\bullet^\partial|)$ is a relative CW complex
it follows that the inclusion $|U^\partial_\bullet|\times I\cup
|U_\bullet|\times\{1\}\hookrightarrow |U_\bullet|\times I$ is a
trivial cofibration, and consequently we may extend the map above to
a continuous map
\[
\tilde{\Phi}\co |U_\bullet|\times I\rightarrow |T(Z_\bullet)|.
\]

The simplicial map $\Phi$ immediately gives a concordance of
$T(Z_\bullet)$-bundles between $i^*_\partial\E_{S^{n-1}}$ and
$\kappa_*\E^\partial$. Furthermore, the adjoint of $\tilde{\Phi}$
gives a concordance of $S_\bullet(|T(Z_\bullet)|)$-bundles from
$i^*\E'$ to a new $S_\bullet(|T(Z_\bullet)|)$-bundle $\E''$. The
restriction of the latter concordance to $S^{n-1}$ is the image of
the first concordance under $\eta_*$. Moreover, $\E''$ has the
following properties:
\begin{itemize}
\item[i)] it is subordinate to the good cover $\U$,
\item[ii)] the restriction of $\E''$ to $S^{n-1}$ corresponds to a
simplicial maps that factors as
\[
U^\partial_\bullet\xrightarrow{\phi^{\partial}}
D(Z)\xrightarrow{\eta\kappa} S_\bullet(|T(Z_\bullet)|).
\]
\end{itemize}
Consequently, the gluing of these concordances gives a commutative
diagram
\[
\xymatrix{ \Con_{T(Z_\bullet)}(D^n,S^{n-1};\kappa_*\E^\partial)
\ar[r]^-{\text{\ding{174}}} \ar[d]^{\cong} &
\Con_{S_\bullet(|T(Z_\bullet)|)}(D^n,S^{n-1};\eta_*\kappa_*\E^\partial) \ar[d]^{\cong} \\
\Con_{T(Z_\bullet)}(D^n,S^{n-1};i_\partial^*\E_{S^{n-1}})
\ar[r]^-{\text{\ding{173}}} &
\Con_{S_\bullet(|T(Z_\bullet)|)}(D^n,S^{n-1};\eta_*i_\partial^*\E_{S^{n-1}}),
}
\]
where $i^*\E'$ in the lower left corner lifts to $\E''$ in the upper
left corner.

Translating Lemma~\ref{Lem:corelemma} into the language of
$?_\bullet$-bundles we get a $D(Z_\bullet)$-bundle $\E$ over $D^n$
whose restriction to $S^{n-1}$ equals $\E^{\partial}$, together with
a concordance relative to $S^{n-1}$ between the
$S_{\bullet}(|T(Z_\bullet)|)$-bundles $\E''$ and $\eta_*\kappa_*\E$.
Thus the class of $\E''$ is in the image of \ding{174}. This implies
that the class of $\E'$ is in the image of \ding{172}, and we are
done.
\end{proof}

\begin{Cor}\label{Cor:Tbijection}
Assume that $Z_{\bullet}$ is good and satisfies the topological Kan
condition. Let $A\hookrightarrow X$ be a closed cofibration such
that $(X,A)$ has the homotopy type of a CW-pair. Let $\E_A$ be any
$Z_\bullet$-bundle over $A$. Then the map
\[
\eta_*\co\Con_{T(Z_{\bullet})}(X,A;\E_A)\rightarrow
\Con_{S_\bullet(|T(Z_{\bullet})|)}(X,A;\eta_*\E_A)
\]
is a bijection.
\end{Cor}

\begin{proof}
This follows directly from the proposition above together with
Proposition~\ref{Prop:finite}.
\end{proof}

We have now almost proved the following theorem:

\begin{Thm}\label{thm:generalizationofmain}
Assume that $X$ has the homotopy type of a CW-complex, and that
$Z_{\bullet}$ is good. Then geometric realization induces a
bijection $\Con_{Z_\bullet}(X)\xrightarrow{\cong}[X,|Z_\bullet|]$.
\end{Thm}

In some sense this theorem is related
to~\cite[Proposition~A.1.1]{MadsenWeiss:07}. Let us now explain how
the our theory of concordances gives a proof:

\begin{proof}
By Theorem~\ref{Thm:topolKanreplacement} there is a good simplicial
space $\tilde{Z}_\bullet$ satisfying the topological Kan condition
together with a simplicial map $q\co
Z_\bullet\rightarrow\tilde{Z}_\bullet$ inducing bijections of
concordance classes and homotopy classes of maps. Now consider the
following diagram:
\[
\xymatrix{
\Con_{Z_\bullet}(X) \ar[r]^{\text{\ding{172}}} \ar[d] &
\Con_{\tilde{Z}_\bullet}(X) \ar[r]^{\text{\ding{173}}} \ar[d] &
\Con_{T(\tilde{Z}_\bullet)}(X) \ar[d]^{\text{\ding{174}}} \\
\Con_{S_\bullet(|Z_\bullet|)}(X) \ar[r] &
\Con_{S_\bullet(|\tilde{Z}_\bullet|)}(X) \ar[r] &
\Con_{S_\bullet(|T(\tilde{Z}_\bullet)|)}(X)\\
\Con_{C_\bullet(|Z_\bullet|)}(X) \ar[u]_{\text{\ding{175}}} \ar[r] \ar[d]^{\text{\ding{176}}} &
\Con_{C_\bullet(|\tilde{Z}_\bullet|)}(X) \ar[u]_{\text{\ding{175}}} \ar[r] \ar[d]^{\text{\ding{176}}} &
\Con_{C_\bullet(|T(\tilde{Z}_\bullet)|)}(X) \ar[u]_{\text{\ding{175}}} \ar[d]^{\text{\ding{176}}} \\
[X,|Z_\bullet|] \ar[r]^{\text{\ding{178}}} &
[X,|\tilde{Z}_\bullet|] \ar[r]^{\text{\ding{177}}} &
[X,|T(\tilde{Z}_\bullet)|].
}
\]
The maps \ding{172} and \ding{178} are induced by $q$ and they are
bijections by Theorem~\ref{Thm:topolKanreplacement}. Since the
natural map $\tilde{Z}_\bullet\rightarrow T(\tilde{Z}_\bullet)$ is a
homotopy equivalence in each simplicial degree, it follows from
Theorem~\ref{Thm:ConDegwiseHEQ} that \ding{173} is a bijection.
Moreover, the map \ding{177} is a bijection for the same reason, see
Lemma~\ref{Lem:DTZagree}. By Corollary~\ref{Cor:Tbijection} the map
\ding{174} is a bijection. All maps marked \ding{175} are bijections
by Proposition~\ref{prop:conCS}, and all maps marked \ding{176} are
bijections by Proposition~\ref{prop:conhtpy}. Consequently, we get
our natural bijection by composing these maps.
\end{proof}

Observe that we easily could generalize the theorem to a CW-complex
$X$ relative to fixed data over some subcomplex $A$.

Of course we have the case $Z_{\bullet}=\Delta 2\C$ in mind, and
when $2\C$ is a good topological $2$-category the condition of
Theorem~\ref{thm:generalizationofmain} is satisfied. Thus we get a
proof of the Theorem~\ref{thm:main}. Reasonable criteria implying
that $2\C$ is good are provided by Theorem~\ref{Thm:good}.





%

\end{document}